\documentclass{tran-l}

\usepackage{amsfonts, amsmath, amsthm}
\usepackage{epsfig, psfrag}

\theoremstyle{definition}

\theoremstyle{remark}

\numberwithin{equation}{section}

\swapnumbers \theoremstyle{plain}

\newtheorem{thm}{Theorem}[section]
\newtheorem*{thm2}{Theorem}
\newtheorem{lem}[thm]{Lemma}

\newtheorem{cor}[thm]{Corollary}

\newtheorem{conj}[thm]{Conjecture}
\newtheorem*{conj2}{Conjecture}

\newtheorem{prop}[thm]{Proposition}

\newtheorem{ques}[thm]{Question}

\theoremstyle{remark}

\theoremstyle{definition}

\newcommand{\T}{\mathcal T}
\newcommand{\C}{\mathcal C}

\renewcommand{\P}{{\mathcal P}}
\renewcommand{\S}{{\mathcal S}}
\renewcommand{\t}{{\tau}}
\newcommand{\bdy}{\partial}
\newcommand{\bbb}{\mathbb}
\newcommand{\rppp}{\mathbb{R}P^3}
\newcommand{\rpp}{\mathbb{R}P^2}

\newcommand{\td}{\tilde}
\newcommand{\open}[1]{\stackrel{\circ}{#1}}

\newcommand{\lra}{\leftrightarrow}
\newcommand{\un}{\underline}

\newcommand{\abs}[1]{\lvert#1\rvert}

\begin{document}

\title{Layered-triangulations of $3$--manifolds}

\author{William Jaco}
\address{Department of Mathematics, Oklahoma State University,
Stillwater, OK 74078}

\email{jaco@math.okstate.edu}
\thanks{The first author was partially supported by NSF Grants
 DMS9971719 and DMS0204707,
The Grayce B. Kerr Foundation, and The American Institute of
Mathematics (AIM)}

\author{J.~Hyam Rubinstein}
\address{Department of Mathematics and Statistics,
University of Melbourne, Parkville, VIC 3052, Australia}
\email{rubin@maths.unimelb.edu.au}
\thanks{The second author was partially supported by The Australian Research Council, The
Grayce B. Kerr Foundation and Stanford University}


\date{\today}


\begin{abstract} A family of one-vertex triangulations of
the genus-$g$-handlebody, called layered-triangulations, is defined.
These triangulations induce a one-vertex triangulation on the
boundary of the handlebody, a genus $g$ surface. Conversely, any
one-vertex triangulation of a genus $g$ surface can be placed on the
boundary of the genus-$g$-handlebody in infinitely many distinct
ways; it is shown that any of these can be extended to a
layered-triangulation of the handlebody. To organize this study, a
graph is constructed, for each genus $g\ge 1$, called the $L_g$
graph; its $0$--cells are in one-one correspondence with equivalence
classes (up to homeomorphism of the handlebody) of one-vertex
triangulations of the genus $g$ surface on the boundary of the
handlebody and its $1$--cells correspond to the operation of a
diagonal flip (or $2\lra 2$ Pachner move) on a one-vertex
triangulation of a surface. A complete and detailed analysis of
layered-triangulations is given in the case of the solid torus ($g
=1$), including the classification of all normal and almost normal
surfaces in these triangulations. An initial investigation of normal
and almost normal surfaces in layered-triangulations of higher
genera handlebodies is discussed.  Using Heegaard splittings,
layered-triangulations of handlebodies can be used to construct
special one-vertex triangulations of $3$-manifolds, also called
layered-triangulations. Minimal layered-triangulations of lens
spaces (genus one manifolds) provide a common setting for new proofs
of the classification of lens spaces admitting an embedded non
orientable surface and the classification of embedded non orientable
surfaces in each such lens space, as well as a new proof of the
uniqueness of Heegaard splittings of lens spaces, including $S^3$
and $S^2\times S^1$. Canonical triangulations of Dehn fillings
called triangulated Dehn fillings are constructed and applied  to
the study of Heegaard splittings and efficient triangulations of
Dehn fillings. It is shown that all closed $3$--manifolds can be
presented in a new way, and with very nice triangulations, using
layered-triangulations of handlebodies that have special one-vertex
triangulations of a closed surface on their boundaries, called
$2$--symmetric triangulations. We provide a quick introduction to a
connection between layered-triangulations and foliations. Numerous
questions remain unanswered, particularly in relation to the
$L_g$-graph, $2$--symmetric triangulations of a closed orientable
surface, minimal layered-triangulations of genus-$g$-handlebodies,
$g\ge 1$ and the relationship of layered-triangulations to
foliations.

\end{abstract}

\maketitle

\tableofcontents

\section{Introduction} This work began with an interest in
managing triangulations of Dehn fillings. Specifically, suppose $X$
is a compact $3$--manifold with boundary a torus; $\T$ is a
one-vertex triangulation of $X$; and $\alpha$ is a slope on $\bdy
X$. We wanted a one-vertex triangulation $\T(\alpha)$ of the Dehn
filled manifold $X(\alpha)$ that restricts to $\T$ on $X$; i.e.,
$\T(\alpha)$ is an extension of $\T$. Furthermore, we needed to
understand  how normal surfaces in the triangulation $\T(\alpha)$
could meet the solid torus of the Dehn filling. It turned out that
what we are now calling layered-triangulations of a solid torus are
exactly what we needed. As we refined our understanding of these
triangulations, we realized that layered-triangulations of the solid
torus lead to nice triangulations of all lens spaces, again called
layered-triangulations, as well as provide an excellent tool to
study Heegaard splittings of Dehn fillings. Moreover,
layered-triangulations of the solid torus and of lens spaces (genus
$1$ Heegaard splittings) raise the issue of extending these ideas to
triangulations of genus-$g$-handlebodies and nice triangulations of
$3$--manifolds having genus $g$ Heegaard splittings. We define
layered-triangulations of higher genus handlebodies analogously to
that for the solid torus and lay a foundation for further study of
these triangulations and their use leading to layered-triangulations
of all $3$--manifolds.

We provide a thorough study of layered-triangulations of the solid
torus and of lens spaces. An understanding of layered-triangulations
of the solid torus is necessary in our work on $1$--efficient
triangulations \cite{jac-rub1}; however, we have carried it much
further. Once we realized that layered-triangulations of lens spaces
provide a common framework in which to recover many of the
interesting results about lens spaces, using quite elementary
techniques, we carried out the various details. This reinforced our
belief that nice triangulations of $3$--manifolds can be used for a
better understanding of the topology of $3$--manifolds and indeed
might provide a combinatorial tool for a better understanding of
their geometry.

In Sections 4 and 5 we consider the mechanics of
layered-triangulations of the solid torus and classify the normal
and almost normal surfaces in these triangulations. A way to
visualize this situation is to think of the infinity of
possibilities of placing a one-vertex (two-triangles, minimal)
triangulation of the torus on the boundary of the solid torus. Two
such triangulations on the boundary of the solid torus are
considered equivalent if there is a homeomorphism of the solid torus
taking one to the other. We have that the equivalence classes of
one-vertex triangulations on the boundary of the solid torus are in
one-one correspondence with the rational numbers $p/q, 0<p<q$, $p$
and $q$ relatively prime, including the forms $0/1, 1/1$. We next
consider extending these triangulations to triangulations of the
solid torus; layered-triangulations provide such extensions. We call
a layered-triangulation extending a $p/q$-triangulation on the
boundary of the solid torus a $p/q$--layered-triangulation (of the
solid torus). We conjecture that the unique minimal
$p/q$--layered-triangulation is the minimal such extension of the
$p/q$--triangulation on the boundary, in the sense of using the
minimal number of tetrahedra. This conjecture is one of the
compelling open problems. A graph having vertices corresponding to
an equivalence class of one-vertex triangulations on the boundary of
the solid torus and having edges corresponding to the elementary
move from one such triangulation to another via a ``diagonal flip"
can be defined, called the $L$-graph. We then have that
layered-triangulations of the solid torus correspond to certain
paths  and minimal layered-triangulations correspond to certain
minimal arcs in the $L$-graph. Finally, we show that a minimal
layered-triangulation of a solid torus has only finitely many
connected normal and almost normal surfaces; we classify them
completely.

In Sections 6 and 7 we define layered-triangulations of lens spaces
and classify the normal and almost normal surfaces in minimal
layered-triangulations of lens spaces. All lens spaces have
layered-triangulations. A very effective method to study them is to
start with a minimal $p/q$--layered-triangulation of the solid torus
and then use any one of the three possible ways to identify the two
triangles on the boundary; each identification gives a lens space.
We provide the combinatorics giving the lens spaces obtained from a
$p/q$--layered-triangulation; and, conversely, for the lens space
$L(X,Y)$, we give the solution as to what layered-triangulation of
the solid torus will give $L(X,Y)$ upon identifying its boundary
faces. What is quite pleasing is the simplicity provided for
understanding the normal and almost normal surfaces in minimal
layered-triangulations of lens spaces. The only orientable embedded,
normal surfaces are small neighborhoods of edges in the
$1$--skeleton; the triangulations have only one-vertex and so all
edges are simple closed curves. There are no almost normal octagonal
surfaces, except in the one-tetrahedron layered-triangulation of
$S^3$, the two-tetrahedra layered-triangulations of $\rppp$ and
$S^2\times S^1$, and one (the ``bad" one) of the two distinct
two-tetrahedra minimal layered-triangulations of $L(3,1)$. Finally,
the only normal non orientable surface embedded in a minimal
layered-triangulation of a lens space is the unique minimal genus
and incompressible one, except in the two-tetrahedra
layered-triangulations of $S^2\times S^1$,  which admits an
embedded, normal Klein bottle, and of $\rppp$, which admits an
embedded normal non orientable surface of genus $3$ in addition to
$\rpp$. We leave unanswered the analogous question to that above
regarding the minimal triangulation extending a $p/q$-triangulation
on the boundary of the solid torus; namely, is the minimal
layered-triangulation of a lens space the minimal triangulation of a
lens space. We conjecture it is (this has been independently
conjectured by S. Matveev  in \cite{matveev2} but using the language
of simple spines and handle decompositions).

In Section 8, we provide applications of these methods. We recover
the earlier work of G. Bredon and J. Wood \cite{bredon-wood}, giving
the classification of lens spaces admitting embedded, non orientable
surfaces and classify such surfaces in these lens spaces. Similarly,
our methods lend themselves to the study of Heegaard splittings. We
recover the earlier work of F. Bonahon and J.P. Otal
\cite{bon-otal}, classifying Heegaard splittings of lens spaces,
including the cases of $S^3$ and $S^2\times S^1$ done by F.
Waldhausen \cite{wald-HS3}. In particular, in a minimal
layered-triangulation of a lens space, distinct from $S^3$, an
almost normal surface is a Heegaard surface if and only if it is the
vertex-linking $2$--sphere with thin edge-linking tubes and an
almost normal tube along a thick edge or the vertex-linking
$2$--sphere with thin edge-linking tubes and an almost normal tube
at the same level as a thin edge-linking tube. There are no normal
Heegaard surfaces in a minimal layered-triangulation of a lens space
distinct from $S^3$. The one-tetrahedron triangulation of $S^3$ has
the vertex-linking $2$--sphere as a normal Heegaard surface and,
necessarily, has an octagonal almost normal $2$--sphere as a
Heegaard surface. Finally, we use these triangulations to give
examples of $0$-- and $1$--efficient triangulations, as well as
examples that are neither. In fact, we classify those
layered-triangulations of lens spaces that are $0$--efficient and
those that are $1$--efficient. All lens spaces, except $\rppp$ and
$S^2\times S^1$, admit infinitely many $0$--efficient triangulations
but a given lens space admits only finitely many $1$--efficient
triangulations.

While these facts about lens spaces are not new, we encourage the
reader to consider the methods and the use of nice-triangulations
for understanding the topology and geometry of $3$--manifolds. These
triangulations also provide examples for students wanting to see
real normal and almost normal surfaces in real triangulations;
however, we believe possibly their most important feature is that
they might guide us in using layered-triangulations in the study of
manifolds of higher genus than one.

In Subsection $8.4$, we define triangulated Dehn fillings. These are
canonical triangulations for Dehn fillings of a knot-manifold. We
use these for consideration of Heegaard splittings of Dehn fillings.
A reader familiar with work done in
\cite{mor-rubHS-curved,shar-HS,rieckHS,rieck-sedgHSfinite,rieck-sedgHSstructure}
will immediately recognize how layered-triangulations of a solid
torus and almost normal, strongly irreducible Heegaard splittings
fit together in a study of Heegaard splittings of Dehn fillings. We
provide a detailed proof of the results in
\cite{mor-rubHS-curved,shar-HS}, using our combinatorial methods to
provide a contrast with the more difficult topological methods. In
fact, we are able to give a complete analysis of the intersection
between a normal or almost normal Heegaard surface and a solid torus
with a layered-triangulation without any restrictions regarding the
nature of the intersection curves on the boundary of the solid
torus. We complete this subsection with conditions as to when a
$0$--efficient triangulation of a knot-manifold extends to a
$0$--efficient triangulated  Dehn-filling. We conclude that minimal
triangulations of knot-manifolds, which must be $0$--efficient,
extend to $0$--efficient triangulated Dehn fillings of that manifold
in all but finitely many cases. A similar result is true for
$1$--efficient triangulations but we do not give that here. It will
appear in \cite{jac-rub1}.

In Sections 9 and 10 we look at layered-triangulations of
handlebodies and of $3$--manifolds in general. These sections, we
believe, give exciting new ways to study $3$--manifolds via Heegaard
splittings with the aid of nice triangulations and may be the most
important sections in this paper.

As in the above situation with the solid torus, we consider
one-vertex (minimal) triangulations of the genus $g$ closed,
orientable surface on the boundary of the genus-$g$-handlebody. We
consider two such triangulations equivalent if there is a
homeomorphism of the handlebody throwing one to the other (including
orientation reversing homeomorphisms). A graph having vertices
corresponding to an equivalence class of one-vertex triangulations
on the boundary of the genus-$g$-handlebody and having edges
corresponding to the elementary move from one such triangulation to
another via a ``diagonal flip" can be defined, called the
$L_g$-graph. This graph admits an extension to a higher dimensional
complex, just as the complex used by J. Harer \cite{harer} and L.
Mosher \cite{mosher-guide} in studying the mapping class group. We
believe the $L_g$ graph is a very important graph to study, as well
as its extension and relationship to the Harer complex.

We define layered-triangulations of a handlebody in a manner
analogous to the definition in the genus one case. We show that any
one-vertex-triangulation on the boundary of the genus-$g$-handlebody
can be extended to a layered-triangulation of the
genus-$g$-handlebody. It follows that layered-triangulations of the
genus-$g$-handlebody correspond to paths in the $L_g$ graph. If
$[\tau]$ is an equivalence class of a one-vertex triangulation on
the boundary of a genus-$g$-handlebody,  we call a
layered-triangulation extending $\tau$, a
$[\tau]$-layered-triangulation. A minimal such layered-triangulation
extending the class $[\tau]$ is called a minimal $[\tau]$-layered
triangulation. Since the $L_g$ graph is not simply connected, it is
suspected that there are many distinct minimal
$[\tau]$-layered-triangulations of the genus-$g$-handlebody. We
examine the normal surfaces in a minimal (smallest number of
tetrahedra possible) four-tetrahedra triangulation of the
genus-$2$-handlebody. According to Ben Burton, using
{\footnotesize\textsc REGINA} \cite{burton-regina}, there are $196$
distinct minimal triangulations of the genus-$2$-handlebody; we do
not know if all are layered. The embedded normal surfaces in higher
genera handlebodies are much richer than what we discovered in genus
one; for example, in the particular four-tetrahedra triangulation of
the genus-$2$-handlebody we consider, there is an infinite family of
distinct normal meridional disks. However, we suspect that while the
task may be much harder, we should be able to understand normal and
almost normal surfaces in these triangulations.

Given a genus $g$ Heegaard splitting of a $3$--manifold, we can
select any one-vertex triangulation on the genus $g$ Heegaard
surface, extend this triangulation to each of the handlebodies in
the splitting, and gain a one-vertex triangulation of the
$3$--manifold, which we call a layered-triangulation. Just as above
for lens spaces, we can modify this view to achieve a new and very
curious presentation for $3$--manifolds. The description we want
also leads to an interesting phenomenon of minimal triangulations of
the genus $g$ closed surface.  In the case of genus one, we started
with a layered-triangulation of the solid torus and observed that if
we use any one of the three possible orientation reversing
identifications of the two faces in the boundary, we get a lens
space and a layered-triangulation of that lens space. Such
identifications of the two faces in the boundary can equally well be
viewed as an identification of the two faces with the one-triangle
M\"obius band. We have three possibilities because the edges in the
one-vertex triangulation of the torus all have a symmetry about them
that allows a simplicial involution, keeping the edge fixed. The
quotient of this involution is the one-triangle M\"obius band. For
genus $g$, we do not have such $2$--symmetry about all edges in a
minimal (one-vertex) triangulation; in fact, some of these
triangulations have no $2$--symmetry at all about an edge. However,
we do have lots of these triangulations with $2$--symmetry, which
can be seen by considering the genus $g$ orientable surface, $S$, as
a two-sheeted branched cover over a compact $2$--manifold, $B$,
having one boundary component and Euler characteristic $1-g$, with
the branching locus the boundary of $B$. Then a minimal (one-vertex)
triangulation of $B$ lifts to a $2$--symmetric minimal (one-vertex)
triangulation of $S$. Hence, every closed $3$--manifold can be
obtained from a layered-triangulation of a handlebody having a
$2$--symmetric minimal triangulation on its boundary by identifying
the faces in the boundary of the triangulation via the associated
simplicial involution. We call such a presentation of a
$3$--manifold a triangulated Heegaard splitting. We suspect that the
most fruitful analysis of $3$--manifolds coming from such
presentations will be in using minimal
$[\tau]$--layered-triangulations of handlebodies, where $\tau$ is
$2$--symmetric. There are numerous open and very compelling
questions regarding these constructions. In Section 10, we examine
many of these questions and explore some possibilities in examples.

Finally, we suggest on first reading, one might skim through
Sections 2-7.  Followed by starting a close reading in Section 8 and
referring back to earlier sections for possible definitions and
clarification. The applications of layered-triangulations in Section
8 give the model for some of our interest in nice triangulations.
Sections 9 and 10 leave more questions than answers. In fact,
layered-triangulations of handlebodies, an understanding of
$L_g$-graphs, and the applications of triangulated Heegaard
splittings to the study of $3$--manifolds is entirely unexplored
ground. We have not had time to carry our investigations further in
this area but find it very compelling.

We have taken a very long time to write up this work, which was done
between 1996-98. Related work, with the knowledge of our work, has
been published in \cite{jac-sedg-dehn} and in \cite{matveev-fom} and
\cite{fom}. Also, Ben Burton has put his particular spin on some of
this work in his thesis \cite{burton-thesis} and used it in the
development of the program {\footnotesize\textsc{REGINA}}
\cite{burton-regina}. His view makes for interesting reading
complementary to this. We would like to thank David Letscher and
Eric Sedgwick, as well as Ben Burton, for many interesting
discussions regarding this material.

\section{Triangulations, normal and almost normal surfaces}

We remark that our triangulations are more likely to be encountered
in the literature under the name pseudo-triangulations. Furthermore,
since the closure of our tetrahedra are not embedded, we need to add
a technical caveat to the definition of both normal and almost
normal surfaces in these triangulations. We feel these slight
modifications are well worth the effort, as we find these
triangulations to be nicely suited for algorithms and the study of
the topology of $3$--manifolds. The reader is referred to our paper
on $0$--efficient triangulations \cite{jac-rub0} for a more complete
discussion of these notions.

Suppose  $\boldsymbol{\Delta}=
\td{\Delta}_1\amalg\td{\Delta}_2\amalg\ldots\amalg\td{\Delta}_n\amalg\ldots$
is a disjoint union of tetrahedra. If $\sigma_i$ and $\sigma_j$ are
distinct faces of $\Delta_{i}$ and  $\Delta_{j}$, respectively, then
we call an affine isomorphism  $\phi_{i,j}:\sigma_i\lra \sigma_j$ a
{\it face pairing} or {\it face identification}. While $\sigma_i$
and $\sigma_j$ are always distinct faces, it is possible $i = j$.
Since a face pairing is an affine isomorphism we think of the
pairing going in either direction.  A family of face pairings
induces an equivalence relation on  $\boldsymbol{\Delta}$. We denote
the identification space determined by this relation
$\boldsymbol{\Delta}/\boldsymbol{\Phi}$ and the associated
projection map by
$\rho:\boldsymbol{\Delta}\rightarrow\boldsymbol{\Delta}/\boldsymbol{\Phi}$.
Typically, the identification space
$\boldsymbol{\Delta}/\boldsymbol{\Phi}$ is not a $3$--manifold. If
we require the family of face pairings to be {\it monogamous}; i.e.,
if a face is in a face pairing, then it is in only one face pairing,
then $\boldsymbol{\Delta}/\boldsymbol{\Phi}$ will be a $3$--manifold
at all points except possibly at the image of certain edges and
certain vertices. If we take a bit more care regarding orientation
of identified edges, then the identification space will be a
$3$--manifold at all points, except possibly at the image of the
vertices. This will then provide us with a satisfactory
combinatorial structure.
 In the case $\boldsymbol{\Delta}/\boldsymbol{\Phi}$
is homeomorphic to a manifold $M$, we set $\T =
(\boldsymbol{\Delta},\boldsymbol{\Phi})$ and say $\T$ is a
triangulation of $M$. If $\boldsymbol{\Delta}/\boldsymbol{\Phi}$
minus the image of the vertices is a $3$--manifold, then it is
homeomorphic to the interior of a compact $3$--manifold $X$; in this
case, we say $\T$ is an ideal triangulation of $\open{X}$, the
interior of $X$. A vertex is then called an {\it ideal vertex}. In
either of the cases that $\T$ is a triangulation or an ideal
triangulation, we have the interiors of the simplices in the
tetrahedra of $\boldsymbol{\Delta}$ embedded. Thus we use the terms
vertex, edge, face and tetrahedron for the images under the
identification map of a vertex, edge, face or tetrahedron.

We remark that the probability of getting a $3$--manifold, even
after avoiding the problems in faces and about edges, goes to zero
as the number of tetrahedra in $\boldsymbol{\Delta}$ increases
without bound. Finally, since the second derived subdivision of such
a triangulation of a $3$--manifold gives a PL-triangulation in the
classical sense, we use classical PL terminology such as ``small
regular neighborhood" to mean a suitable small regular neighborhood
in some derived PL-subdivision.

\begin{figure}[htbp]
\psfrag{a}{\scriptsize{$a$}}\psfrag{b}{\scriptsize{$b$}}\psfrag{c}{\scriptsize{$c$}}
\psfrag{d}{\scriptsize{$d$}}\psfrag{e}{\scriptsize{$e$}}

            \psfrag{1}{\footnotesize{$$}}

            \psfrag{A}{\Large{$(A)$}}
\psfrag{B}{\Large{$(B)$}} \psfrag{S}{\Large{$S^3$}}

        \vspace{0 in}
        \begin{center}
\epsfxsize = 2.5 in \epsfbox{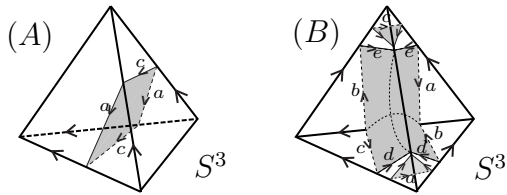}
 \caption{A normal torus (A) in the one-tetrahedron, two-vertex triangulation
of $S^3$ and an almost normal $2$--sphere (B) in the
one-tetrahedron, one-vertex triangulation of $S^3$.}
\label{f-normal-almost-normal}
\end{center}
\end{figure}

We shall assume the reader is familiar with the basics of normal and
almost normal surface theory. However, observe that since our
tetrahedra, faces and edges are not embedded (only their interiors
are embedded), it is not necessarily the case that a normal surface
in $M$ is a surface that meets the tetrahedron in $M$ in a
collection of normal triangles and normal quadrilaterals but rather
a normal surface in $M$ is a surface that meets the tetrahedra in
$M$ so that the lifts  of the intersection of the surface with the
tetrahedra in $M$ are normal triangles and normal quadrilaterals in
the tetrahedra of  $\boldsymbol{\Delta}$. We exhibit this in the two
one-tetrahedron triangulations of $S^3$ in Figure
\ref{f-normal-almost-normal}. However, this is a technicality that
cause no problems and will be implicitly understood.

\section{One-vertex triangulations} A closed surface with non positive
Euler characteristic admits a one-vertex triangulation (a
triangulation of the $2$--sphere must have at least three vertices
while one of $\rpp$ must have at least two vertices).  One of the
consequences of the work in \cite{jac-rub0} is that irreducible,
orientable $3$--manifolds have one-vertex triangulations. Also, in
Section \ref{lay-handle}, we provide a proof that all closed,
orientable $3$--manifolds (including reducible ones) admit a
one-vertex triangulation. (Such results also follow from the
existence of special spines for closed, orientable $3$--manifolds
\cite{matveev1,casler}.)  A one-vertex triangulation of a closed
surface is minimal. On the other hand, by using $2\leftrightarrow 3$
Pachner moves \cite{pachner1}, it is easy to see that any
$3$--manifold that has a one-vertex triangulation has infinitely
many; so, in particular, a fixed $3$--manifold can have a one-vertex
triangulation with an arbitrarily large number of tetrahedra. In
\cite{jac-rub0}, it is shown that a minimal triangulation of an
irreducible, orientable $3$--manifold, except $S^3, \rppp$ and
$L(3,1)$, has only one vertex. ($S^3$ has $2$ one-tetrahedron
triangulations, $1$ with one vertex and $1$ with two vertices;
$\rppp$ has $2$ two-tetrahedra (minimal) triangulations, $1$ with
one vertex and $1$ with two vertices; and $L(3,1)$ has $3$
two-tetrahedra (minimal) triangulations, $2$ with one vertex and $1$
with two vertices. We do not know if minimal triangulations of
reducible $3$--manifolds have just one vertex, but we suspect this
also is true.  Finally, we point out a very interesting problem;
namely, given a triangulation of a closed $3$--manifold, can it be
decided if the triangulation is minimal? An affirmative solution to
the Homeomorphism Problem for $3$--manifolds enables an algorithm to
decide if a given triangulation of a manifold is minimal; the
Homeomorphism Problem follows from the Thurston Geometrization
Theorem, which recently has been  announced by G.Perelman
\cite{perelman}. However, this certainly strikes us as a circuitous
resolution of the problem. From the census given in
\cite{matveev2,matveev3} or the census by Ben Burton given by
{\footnotesize\textsc{REGINA}} \cite{burton-regina}, we have a list
of distinct minimal triangulations of closed, orientable
$3$--manifolds up to seven tetrahedra (\cite{matveev3}), extended to
nine tetrahedra (\cite{italians-3}), and a census extending these to
bounded $3$--manifolds, non orientable $3$--manifolds and ideal
triangulations of the interior of some compact $3$--manifolds, up to
seven tetrahedra (\cite{burton-thesis}).

We have the following observations about one-vertex triangulations
of closed surfaces.

\begin{lem} [\cite{jac-sedg-dehn}] In a one-vertex triangulation of a closed surface every
trivial, normal curve is vertex-linking.\end{lem}

\begin{lem} [\cite{jac-sedg-dehn}] In a one-vertex triangulation of the torus, two normal
curves are normally isotopic if and only if they are isotopic.
\end{lem}

The previous Lemma is not true for normal curves in one-vertex
triangulations of higher genera (genus greater than one) surfaces.
Typically, in a one-vertex triangulation of a closed orientable
surface other than the torus, there are infinitely many distinct
normal curves in an isotopy class. An easy way to see this is to
consider two disjoint non trivial simple closed curves $\ell$ and
$\ell'$ that co-bound an annulus $A$ containing the vertex ($\ell$
and $\ell'$ are isotopic but for an orientable surface different
from the torus, they are {\it not} isotopic via an isotopy missing
the vertex). We also assume $A$ does not separate the surface. For
any simple closed curve meeting $A$ in a single non trivial spanning
arc, we can make ``finger moves" within $A$ over the vertex;
thereby, generating an infinite family of simple closed curves that
are isotopic but no two are isotopic via an isotopy missing the
vertex. These curves can be normalized to normal simple closed
curves all of which are in the same isotopy class but no two are
normally isotopic. On the other hand, we do have the following
partial generalization for higher genera surfaces.

\begin{lem} In a one-vertex triangulation of a closed surface, two curves
are isotopic via an isotopy missing the vertex if and only if they
are normally isotopic.\end{lem}

This can be proved using the same techniques as in
\cite{jac-sedg-dehn}.

If $Q$ is a quadrilateral, possibly with some edge identifications,
then we can subdivide $Q$ into two triangles, without adding
vertices, by adding one of the two possible diagonals as an edge.
The operation of going from one of these triangulations of $Q$ to
the other is called a {\it diagonal flip}. Now, suppose $F$ is a
surface and $\mathcal{P}$ is a triangulation of $F$. If $e$ is an
interior edge of $\mathcal{P}$, then there are triangles $\sigma$
and $\beta$ meeting along $e$. If $\sigma \ne\beta$, the union $Q =
\sigma\cup \beta$ is a quadrilateral, possibly with some
identifications of boundary edges, with diagonal $e$ and
anti-diagonal the edge $e'$. If we swap the diagonal $e$ for the
diagonal $e'$, then we have a new triangulation $\mathcal{P}'$ of
$F$, which we say is obtained from $\mathcal{P}$ by a {\it diagonal
flip}. We have the following well-known lemma. There is a proof in
\cite{mosher-guide}

\begin{lem}\label{diag-flip} Suppose $S$ is a closed surface and
 $\mathcal{P}$ and $\mathcal{P}'$ are one-vertex
triangulations of $S$. Then there is a sequence of triangulations
$\P = \mathcal{P}_0,\mathcal{P}_1,\ldots,\mathcal{P}_n = \P'$ of
$S$, where $\mathcal{P}_{i+1}$ differs from $\mathcal{P}_i$ by a
diagonal flip and (possibly) an isotopy.
\end{lem}

\begin{figure}[htbp]
\psfrag{1}{\small{$1$}}\psfrag{2}{\small{$2$}}\psfrag{3}{\small{$3$}}
\psfrag{4}{\small{$4$}}\psfrag{5}{\small{$5$}}
            \psfrag{a}{\large{$\sigma$}}
            \psfrag{b}{\large{$\beta$}}
            \psfrag{c}{\large{$\td{\sigma}$}}
            \psfrag{d}{\large{$\td{\beta}$}}
            \psfrag{i}{\Large{$\td{\Delta}$}}
            \psfrag{f}{\Large{${\Delta}$}}
            \psfrag{g}{\Large{$\P'$}}
            \psfrag{h}{\Large{$\P$}}
            \psfrag{e}{$e$}
            \psfrag{j}{$e'$}
            \psfrag{k}{$\td{e}$}

        \vspace{0 in}
        \begin{center}
\epsfxsize = 2.5 in \epsfbox{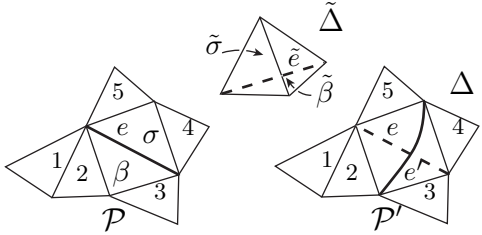}
 \caption{Layering along an edge $e$ results in a diagonal flip taking the triangulation
 $\P$ to the triangulation $\P'$.} \label{f-layer-at-e}
\end{center}
\end{figure}

 Now, suppose $M$ is a compact $3$--manifold with
nonempty boundary, $\T$ is a triangulation of $M$, and $\mathcal{P}$
is the triangulation induced on $\bdy M$ by $\T$.  Furthermore,
suppose $e$ is an edge in $\P$ and there are two distinct triangles
$\sigma$ and $\beta$ in $\P$ meeting along the edge $e$. Let
$\td{\Delta}$ be a tetrahedron distinct from the tetrahedra in $\T$
and let $\td{e}$ be an edge in $\td{\Delta}$. Let $\td{\sigma}$ and
$\td{\beta}$ be the faces of $\td{\Delta}$ that meet along $\td{e}$.
Now, if we identify $\td{e}$ with $e$ and extend this to
 face
identifications from $\td{\sigma}\rightarrow\sigma$ and
$\td{\beta}\rightarrow\beta$, we get a new triangulation $\T'$ of
the $3$--manifold $M$ so that  the triangulation  $\P'$, induced by
$\T'$ on $\bdy M$, is obtained from $\P$ via a diagonal flip in the
quadrilateral $\sigma\cup\beta$. There are two ways to attach
$\td{\Delta}$, depending on the direction we attach $\td{e}$ to $e$
(we have specified that $\td{\sigma}$ goes to $\sigma$ and
$\td{\beta}$ goes to $\beta$); however, the resulting triangulations
are isomorphic. We call this operation a {\it layering on the
triangulation $\T$ along the edge $e$ of $\P$} or simply {\it a
layering on $\T$ along the edge $e$}. The operation on $\P$ is
called a Pachner or bi-stellar move of type $2\leftrightarrow 2$ on
$\P$ along the edge $e$ (as well as a diagonal flip). See Figure
\ref{f-layer-at-e}.

\section{Layered-triangulations of a solid torus}
We consider very special one-vertex triangulations of the solid
torus, which we call ``layered-triangulations." We show that given a
solid torus and a one-vertex triangulation on its boundary, there is
an algorithm to construct a layered-triangulation of the solid torus
that restricts to the given one-vertex triangulation on its
boundary.  In fact, for any positive integer $N$ there is such an
extension to a layered-triangulation of the solid torus having more
than $N$ tetrahedra; however, there is a unique (up to combinatorial
isomorphism) extension of a one-vertex triangulation on the boundary
of a solid torus to a minimal layered-triangulation of the solid
torus. It remains an open question as to whether the minimal
layered-triangulation of a solid torus that extends a one-vertex
triangulation on the boundary of the solid torus is the minimal such
extension. We discuss this and related problems about extending
triangulations from the boundary of a $3$--manifold to
triangulations of the $3$--manifold in later sections.

\subsection{One-vertex triangulations on the boundary of a solid
torus} An isotopy class of a simple, closed, normal curve, $\gamma$,
in a one-vertex triangulation of a torus (which is the same as the
normal isotopy class of $\gamma$) determines a unique unordered
triple of nonnegative integers $\{y_1,y_2,y_3\} = \{\abs{\gamma \cap
e_1},\abs{\gamma \cap e_2},\abs{\gamma \cap e_3}\}$, where $y_i =
\abs{\gamma \cap e_i}$ is the number of times the normal curve
$\gamma$ meets the edge $e_i, i = 1,2,3$ of the triangulation. This
includes the triples $\{0,1,1\}$ and $\{1,1,2\}$. In any case, we
have for one of these numbers, say $y_k$, that $y_k = y_i + y_j$ and
with only the exceptions $\{0,1,1\}$ and $\{1,1,2\}$, the integers
$y_i, y_j$ (hence, also $y_k$) are relatively prime. It follows, in
all cases, that precisely one of the terms $y_i, y_j$ or $y_k$ is
even, as is the sum $y_1 + y_2+ y_3$. Note these are not the normal
coordinates of the normal curve in the one-vertex triangulation;
however, there is an easy translation between normal coordinates and
these triples (see \cite{jac-sedg-dehn}). We use standard
terminology and call an isotopy class of a nontrivial simple closed
curve in a torus a {\it slope}. If $\gamma$ is the unique normal
curve representing a slope, then we will often just use $\gamma$ for
both the slope and the normal curve.

For a solid torus, there is a unique slope in its boundary so that
any curve, with that slope, bounds a disk. We call this slope the
{\it meridional slope} and a disk it bounds is called a {\it
meridional disk}. So, if $\T_{\bdy}$ is a one-vertex triangulation
on the boundary of a solid torus, there is a unique normal curve
representing the meridional slope and thus a unique triple of
nonnegative integers $\{p,q,p+q\}$ determined by the intersections
of the (normal) meridional slope with the edges of the triangulation
$\T_{\bdy}$. We allow the triples $\{0,1,1\}$ and $\{1,1,2\}$; the
first being when the meridional slope is the slope of an edge and
the second is the slope of the anti-diagonal, if we designate the
edge meeting the curve in two points the diagonal. The triple
$\{p,q,p+q\}$ is completely determined once we specify $p$ and $q$
and except for $p=q=1$, we may assume notation is such that $0\le p\
<q$ and, except for $p=0, q = 1$, $p$ and $q$ are relatively prime.
We will call such a one-vertex triangulation on the boundary of a
solid torus a \emph{$p/q$--triangulation} on the boundary of the
solid torus, including $0/1$ and $1/1$ for the triples $\{0,1,1\}$
and $\{1,1,2\}$, respectively. Notice that while there is a unique
one-vertex triangulation of the torus, there are infinitely many
ways that a one-vertex triangulation of the torus can sit upon the
boundary of a solid torus, differing in how the edges of the
triangulation meet the meridional slope. The following lemma shows
that the association of a rational $p/q$ with a triangulation on the
boundary of the solid torus is unique, up to homeomorphism of the
solid torus.

\begin{lem}\label{pq-determine} Suppose $\T_\bdy$  and $\T_\bdy'$ are
$p/q$-- and $p'/q'$--triangulations, respectively, on the boundary
of the solid torus $\bbb{T}$. There is a homeomorphism of $\bbb{T}$
taking the triangulation $\T_\bdy$ to the triangulation $\T_\bdy'$
if and only if  $p/q = p'/q'$.
\end{lem}

\begin{proof} Suppose $\T_\bdy$ and $\T_{\bdy}'$ are one-vertex triangulations on
the boundary of the solid torus $\bbb{T}$. Denote their edges by
$e_1,e_2,e_3$ and $e_{1}',e_{2}',e_{3}'$, respectively. Let $\mu$
denote the meridional curve on the boundary. Since $\mu$ is unique
up to isotopy, it determines a unique normal isotopy class in
$\T_{\bdy}$ and in $\T_{\bdy}'$. As above, $\T_{\bdy}$ determines
a unique triple $\{y_1,y_2,y_3\}$, where $y_i = \abs{e_i\cap\mu}$
and similarly, $\T_{\bdy}'$ determines a unique triple
$\{y_{1}',y_{2}',y_{3}'\}$, where $y_{i}' = \abs{e_{i}'\cap\mu}$.
(These triples are the triples $\{p,q,p+q\}$ and
$\{p',q',p'+q'\}$, respectively, given in the hypothesis. We have
chosen this notation here to avoid ambiguity as we proceed with
the proof.)

Now, if there is a homeomorphism of the solid torus $\bbb{T}$
taking $\T_\bdy$ to $\T_{\bdy}'$, then by the uniqueness of the
normal isotopy class of the meridian, we have $\{y_1,y_2,y_3\}
\equiv \{y_{1}',y_{2}',y_{3}'\}$ as unordered triples.

So, conversely, suppose $\{y_1,y_2,y_3\} \equiv
\{y_{1}',y_{2}',y_{3}'\}$ and notation has been chosen so that
$y_i = y_{i}', i = 1,2,3$ and $y_3 = y_1 + y_2$ ($y_{3}' = y_{1}'
+ y_{2}'$).

Choose orientation on $e_1, e_2$ and $\mu$ so $\mu = y_{2}e_1 +
y_{1}e_2$. Choose a longitude, say $\lambda$, on the boundary of
$\bbb{T}$ so $\lambda = x_{2}e_1 + x_{1}e_2$, where $x_1, x_2$
satisfy $x_{1}y_2 - x_{2}y_1 = 1$. It follows that $$e_{1} =
x_{1}\mu - y_{1}\lambda,$$ $$e_2 = -x_{2}\mu + y_{2}\lambda$$ and
so, $$e_3 = e_2 - e_1 = (-x_1 - x_2)\mu + (y_1 + y_2)\lambda.$$

Similarly, we can choose orientation on $e_{1}'$ and $e_{2}'$ and
$\mu$ so $\mu = y_{2}e_{1}' + y_{1}e_{2}'$ and there exists
$x_{1}', x_{2}'$ such that we can have $\pm\lambda = x_{2}'e_{1}'
+ x_{1}'e_{2}'$ and $x_{1}'y_2 - x_{2}'y_1 = 1$. Notice, we want
$x_{1}'y_2 - x_{2}'y_1 = 1$, therefore, we have the ambiguity
$\pm\lambda$. Hence,
$$e_{1}' = x_{1}'\mu \mp y_{1}\lambda,$$ $$e_{2}' = -x_{2}'\mu \pm y_{2}\lambda$$ and
so, $$e_{3}' = e_{2}' - e_{1}' = (-x_{1}' - x_{2}')\mu \pm (y_1 +
y_2)\lambda.$$

A Dehn twist about $\mu$, repeated $n$ times, takes
$\mu\rightarrow \mu$ and $\lambda\rightarrow n\mu + \lambda.$ So,
$$e_1\rightarrow x_{1}\mu -y_{1}n\mu - y_{1}\lambda = (x_{1} -
y_{1}n)\mu - y_{1}\lambda$$ $$e_2\rightarrow -x_{2}\mu +y_{2}n\mu
+ y_{2}\lambda = (-x_{2} + y_{2}n)\mu + y_{2}\lambda.$$

Recall, we have $x_{1}y_{2} - x_{2}y_{1} = x_{1}'y_2 - x_{2}'y_1$;
so, since $y_1$ and $y_2$ are relatively prime, it follows that
$(x_1 - x_{1}') = y_{1}n$ for some integer $n$. Similarly, we get
$(x_2 - x_{2}') = y_{2}n$ for the same integer $n$. We conclude
that for this $n$, we have $$x_{1}' = x_1 - y_{1}n$$ and $$-x_{2}'
= -x_2 + y_{2}n.$$ Therefore
$$e_1\rightarrow x_{1}'\mu - y_1\lambda,$$ $$e_2\rightarrow -x_{2}'\mu + y_2\lambda,$$
and $$e_3\rightarrow (-x_{1}' - x_{2}')\mu + (y_1 + y_2)\lambda.$$
Hence, either there is a Dehn twist taking  $\T_\bdy$ to
$\T_{\bdy}'$ or a Dehn twist followed by an orientation reversing
homeomorphism of $\bbb{T}$ leaving $\mu$ fixed and reversing
orientation of the longitude $\lambda$ taking $\T_\bdy$ to
$\T_{\bdy}'$.\end{proof}

Since there are different ways in which a one-vertex triangulation
of a torus can sit upon the boundary of a solid torus, it is natural
to wonder if any such triangulation can be extended (without adding
vertices) to a triangulation of the solid torus. Below, we show that
any one-vertex triangulation on the boundary of a solid torus can be
extended to a (one-vertex) triangulation of the solid torus. A proof
of this appears in \cite{jac-sedg-dehn} following our work. As
pointed out above, we leave open the interesting and closely related
question of given a (one-vertex) triangulation on the boundary of a
solid torus, what is the minimal number of tetrahedra needed to
extend this triangulation to a triangulation of the solid torus?
More generally, it can be shown that any triangulation on the
boundary of a compact $3$--manifold can be extended without adding
vertices to a triangulation of the $3$--manifold; so, in particular,
all vertices of such an extension are in the boundary of the
$3$--manifold. In this more general context, we also have the
question of what is the minimal number of tetrahedra needed to
extend a given triangulation on the boundary.

\subsection{Layered-triangulations of a solid torus} If $\T$ is a
one-vertex triangulation of a solid torus $\mathbb{T}$, then the
boundary torus has a triangulation with two faces, three edges and
one vertex. If $e$ is one of the edges in the induced
triangulation on the boundary, then we can add a tetrahedron
$\td{\Delta}$, which is disjoint from the lifts of tetrahedra in
$\T$,  by layering along the edge $e$; the result after layering
is still a solid torus with a triangulation having one more
tetrahedron. We denote the new solid torus by $\mathbb{T}' =
\mathbb{T}\cup_e\Delta$, where $\Delta$ is the image of
$\td{\Delta}$, and denote the new triangulation by $\T' = \T
\cup_e \td{\Delta}$. As we mentioned above, once the edge $e$ is
designated, there are several possibilities for how we attach the
tetrahedron $\td{\Delta}$ along $e$ and its two adjacent faces;
however, in this case, the resulting triangulations are all
isomorphic.

\begin{figure}[htbp]

           \psfrag{D}{\Large{$\td{\Delta}$}}
\psfrag{a}{\small{$a$}} \psfrag{b}{\small{$b$}}
\psfrag{c}{\small{$c$}} \psfrag{d}{\small{$d$}}
\psfrag{1}{\small{$1$}} \psfrag{2}{\small{$2$}}
\psfrag{3}{\small{$3$}}\psfrag{x}{\small{$x$}}\psfrag{y}{\small{$y$}}\psfrag{z}{\small{$z$}}
\psfrag{M}{\begin{tabular}{c}
           $\langle
a,b,c\rangle\leftrightarrow\langle x,y,z\rangle$\\
        $\langle b,c,d\rangle\leftrightarrow\langle
x,y,z\rangle$
            \end{tabular}}

\psfrag{N}{\begin{tabular}{c}$\langle
a,b,c\rangle\leftrightarrow\langle y,z,x\rangle$\\
$\langle d,b,c\rangle\leftrightarrow\langle y,z,x\rangle$
\end{tabular}} \psfrag{O}{\begin{tabular}{c}$\langle
a,b,c\rangle\leftrightarrow\langle z,x,y\rangle$\\
$\langle c,d,b\rangle\leftrightarrow\langle z,x,y\rangle$
\end{tabular}} \psfrag{P}{\begin{tabular}{c}(A)\hspace{.1
in}$\langle a,b,c\rangle\leftrightarrow\langle b,c,d\rangle$
 \end{tabular}}
\psfrag{Q}{\begin{tabular}{c}(B)\hspace{.1 in}$\langle
a,b,c\rangle\leftrightarrow\langle c,d,b\rangle$\\
 \end{tabular}}
\psfrag{R}{\begin{tabular}{c}(C)\hspace{.1 in}$\langle
a,b,c\rangle\leftrightarrow\langle d,b,c\rangle$
 \end{tabular}}
\psfrag{L}{Layer}
        \vspace{0 in}
        \begin{center}
\epsfxsize = 5 in \epsfbox{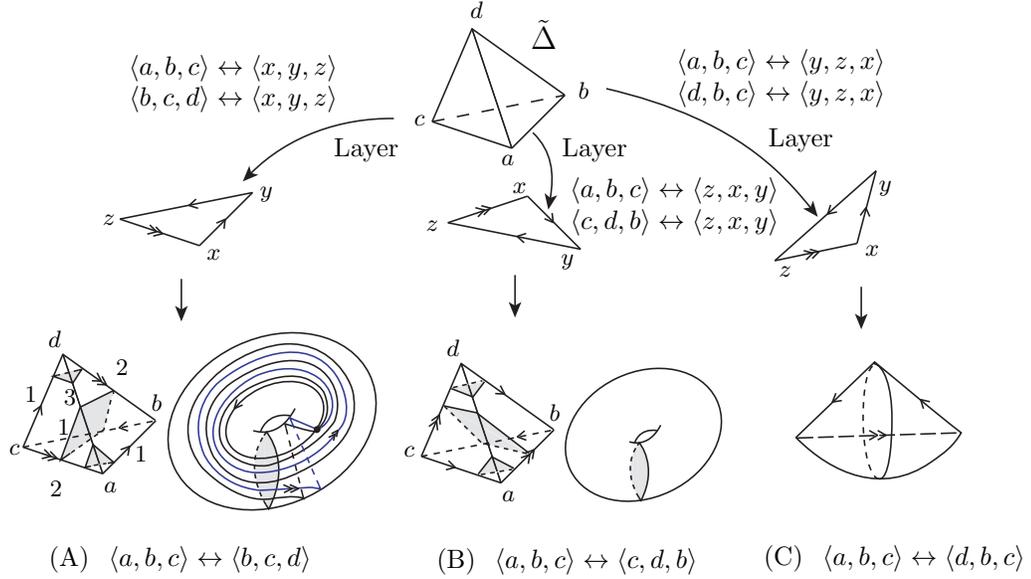}
 \caption{One-tetrahedron solid torus and creased $3$--cell (layering
  of a tetrahedron on a one-triangle M\"obius band).} \label{f-one-tet-layer}
\end{center}
\end{figure}

\begin{figure}[htbp]

            \psfrag{a}{\footnotesize{$1'$}}\psfrag{0}{\footnotesize{$0$}}\psfrag{1}{\footnotesize{$1$}}
            \psfrag{2}{\footnotesize{$2$}}\psfrag{3}{\footnotesize{$3$}}
            \psfrag{x}{\large{(A)} extends $1/1$}
            \psfrag{w}{\large{(B)} extends $0/1$}
            \psfrag{z}{\large{(C)} extends $0/1$}
\psfrag{l}{layer on $3$}\psfrag{m}{layer on $2$}\psfrag{n}{layer on
$1$}

        \vspace{0 in}
        \begin{center}
\epsfxsize = 3.5 in \epsfbox{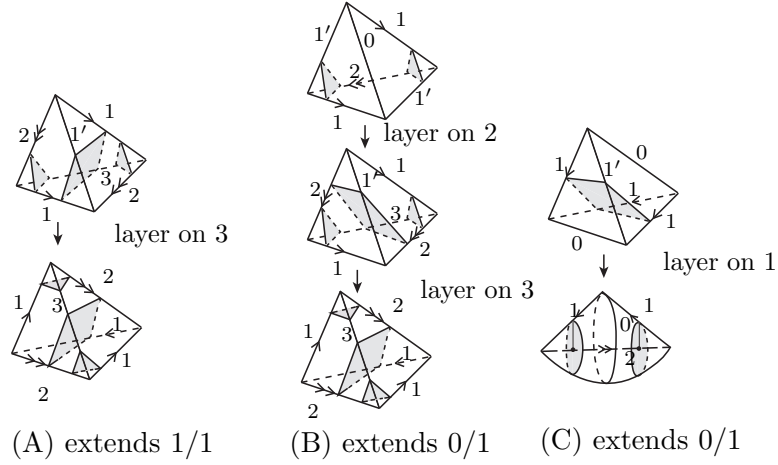} \caption{(A)
Two-tetrahedron (minimal) layered-triangulation of solid torus
extending $1/1$; (B) Three-tetrahedron layered-triangulation of
solid torus extending $0/1$; (C) Two-tetrahedron (minimal)
layered-triangulation of solid torus extending $0/1$}
\label{f-small-layered}
\end{center}
\end{figure}

First, we make some observations about some special small
triangulations of the solid torus. Consider the tetrahedron
$\td{\Delta}$ pictured in Figure \ref{f-one-tet-layer}. There are
three ways to layer the back two faces of the tetrahedron
$\td{\Delta}$ onto the one-triangle M\"obius band. In parts (A) and
(B) the tetrahedron is layered along the interior (orientation
reversing edge) on the one-triangle M\"obius band, the labels and
arrows give the identifications. In these cases, the boundary torus
has a one-vertex triangulation consisting of the two front faces of
the tetrahedron. All four vertices are identified to a single vertex
and the induced edge identifications are indicated in Figure
\ref{f-one-tet-layer}, Parts A and B, respectively. Combinatorially,
these triangulations are the same (exchange the vertices labeled $b$
and $c$). This triangulation of the solid torus will be referred to
as {\it the one-tetrahedron solid torus}. The numbers $1$, $2$ and
$3$ labeling the edges indicate the number of times the edge meets
the boundary of the meridional disk; thus, this is an extension of
the $1/2$--triangulation on the boundary of the solid torus. In the
last case, Figure \ref{f-one-tet-layer}, Part C, we show a {\it
creased $3$--cell} obtained by a single layering of a tetrahedron
along the boundary edge of the one-triangle M\"obius band; again the
labels and arrows give the identification. The M\"obius band and the
creased $3$--cell both have the homotopy type of a solid torus (the
M\"obius band, in particular, plays an important role in our
theory); we think of each as a degenerate layered-triangulation of
the solid torus. In the case of the one-triangle M\"obius band, we
have a (degenerate) extension of a $1/1$--triangulation on the
boundary of the solid torus and in the case of the creased
$3$--cell, we have a (degenerate) extension of a
$0/1$--triangulation on the boundary of the solid torus. Note, if we
layer another tetrahedron along the edge labeled $1$ on the boundary
of the creased $3$--cell, we do get a two-tetrahedron triangulation
of the solid torus extending the $0/1$--triangulation on the
boundary.  This example was shown to us by Eric Sedgwick. See Figure
\ref{f-small-layered}. Layering onto the one-triangle M\"obius band
is useful in using layered triangulations for Dehn fillings and in
our generalization of layered-triangulations of a solid torus to
analogous triangulations of handlebodies.

\begin{figure}[htbp]

            \psfrag{e}{\Large{$e$}}\psfrag{A}{\Large{(A)}}\psfrag{B}{\Large{(B)}}
            \psfrag{D}{\Large{$\td{\Delta}_t$}}
            \psfrag{f}{$\td{e}$}
            \psfrag{e}{$e$}\psfrag{1}{\footnotesize{$1$}}\psfrag{2}{\footnotesize{$2$}}\psfrag{3}{\footnotesize{$3$}}
            \psfrag{4}{\footnotesize{$4$}}\psfrag{5}{\footnotesize{$5$}}\psfrag{6}{\footnotesize{$6$}}
            \psfrag{7}{\footnotesize{$7$}}\psfrag{8}{\footnotesize{$8$}}
            \psfrag{k}{layered}
             \psfrag{l}{\begin{tabular}{c}
            layered\\
        solid\\
        torus\\
            \end{tabular}}
            \psfrag{s}{\begin{tabular}{c}
            $t$\\
       layers\\
            \end{tabular}}
            \psfrag{t}{\begin{tabular}{c}
            $(t-1)$\\
            layers\\
            \end{tabular}}

        \vspace{0 in}
        \begin{center}
\epsfxsize = 3.5 in \epsfbox{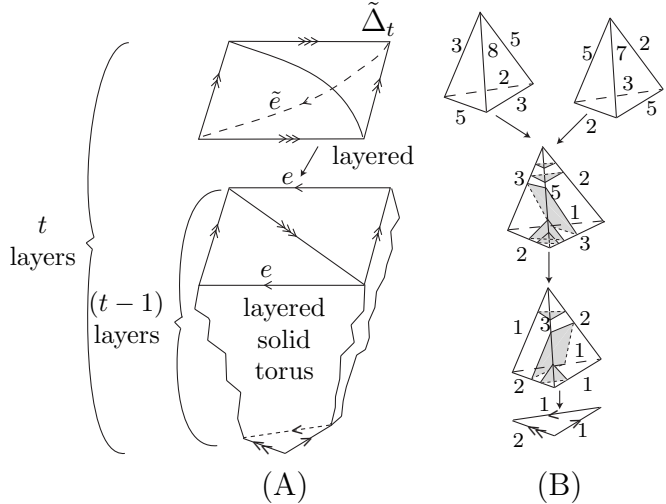} \caption{A
layered-triangulation of a solid torus. In the example, (B), the
numerical values along an edge denote the number of times the
meridional disk meets the edge.} \label{f-layered-torus-def}
\end{center}
\end{figure}

 Inductively, we define a triangulation $\T_t$
to be a {\it layered-triangulation of the solid torus with
$t$--layers} if
\begin{enumerate}
\item $\T_0$ is the one-triangle M\"obius band, \item $\T_1$ is
either the one-tetrahedron solid torus or the creased $3$--cell,
each of which is obtained by a layering of a tetrahedron along an
edge of $T_0$\item $\T_{t} = \T_{t-1} \cup_e \td{\Delta}_t$ is a
layering along the edge $e$ of a layered-triangulation $\T_{t-1}$
having $t-1$ layers, $t \geq 1$. See Figure
\ref{f-layered-torus-def}.
\end{enumerate}

In Figure \ref{f-small-layered}, we give examples of what we might
call``exceptional" layered-triangulations of the solid torus. Part A
is a two-tetrahedron minimal extension of the $1/1$-triangulation on
the boundary; and Part B is a three-tetrahedron
layered-triangulation extending the $0/1$-triangulation on the
boundary. The latter is \emph{not} the minimal extension of the
$0/1$--triangulation on the boundary, which actually is the
two-tetrahedron layered-triangulation, given in Part C.  Notice in
Parts A and B of Figure \ref{f-small-layered}, the edge labeled $3$
is of index two. Generally, in a triangulation with an edge of index
two, we can reduce the number of tetrahedra by crushing the two
tetrahedra containing the common edge to two faces. However, in this
case the crushing does not give a $3$--manifold. In Part B it gives
the M\"obius band and in Part C it gives the creased $3$--cell.
These are supporting observations for thinking of the one-triangle
M\"obius band and the creased $3$--cell as layered-triangulations,
extending the $1/1$-- and the $0/1$--triangulations on the boundary
of a solid torus, respectively.

Except in the degenerate cases, a layered-triangulation of a solid
torus, with $t$ layers, has one vertex, which is in the boundary
torus, $t+2$ edges, three of which are in the boundary torus, and
$2t + 1$ faces, two of which are in the boundary.

Figure \ref{f-layered-torus-def}(B) gives some specific examples of
layered-triangulations of the solid torus. The numbers along the
edges indicate the number of times the edge meets the meridional
disk of the solid torus.

\vspace{.125 in}\noindent{\it Thick and thin edges.} There is a
special edge in a layered-triangulation of a solid torus.
  At the base of a layered-triangulation of a solid torus is the
M\"obius band. Its edges meet the meridional disk $1$ and $2$ times.
The edge that meets the meridional disk $1$ time is called a
\emph{thick edge} in the sense that no small neighborhood of this
edge has boundary a normal torus - the shrinking of the boundary of
a small regular neighborhood of this edge never normalizes and
sweeps across the entire solid torus. On the other hand, every other
edge has arbitrarily small neighborhoods with normal boundaries.
(Even if they meet the meridional disk just once; see, for example,
the edge labeled $1'$ in Part B of Figure \ref{f-small-layered}.)
All edges other than the unique thick edge are called \emph{thin
edges}. If a thin edge is in the interior of the solid torus, then
there are arbitrarily small neighborhoods of the edge having
boundary a normal once-punctured torus; if the thin edge is in the
boundary of the solid torus, then there are arbitrarily small
neighborhoods of the edge having boundary a normal annulus. It is
possible that the thick edge is in the boundary.

There is another edge in a layered-triangulation of a solid torus
that  is distinguished. If $\T$ is a layered-triangulation of a
solid torus, then the unique edge on the boundary having valence one
(meeting only one tetrahedron) is called the \emph{univalent edge}.
There is a unique univalent edge for the one-triangle M\"obius band
(the edge labeled ``$2$") and for the creased $3$--cell (the edge
labeled ``$0$").

\vspace{.125 in}\noindent{\it Slopes and layering.} Suppose $\T$ is
a one-vertex triangulation of a solid torus $\bbb{T}$ and
$\T_{\bdy}$ is the triangulation induced by $\T$ on the boundary of
$\bbb{T}$. Let $e_1,e_2,e_3$ denote the three edges of the
triangulation $\T_{\bdy}$ in the boundary of $\bbb{T}$. Then a
layering along an edge of $\T_{\bdy}$, say for example along the
edge $e_3$, gives a solid torus $\bbb{T}' =
\bbb{T}\cup_{e_3}\Delta$, with a one-vertex triangulation $\T'$. If
$\T_{\bdy}'$ is the triangulation on the boundary of $\bbb{T}'$,
induced by $\T'$, then the layering has the effect of designating
the edge $e_3$ in $\T_{\bdy}$ as the diagonal and ``flipping the
diagonal'' in going from the triangulation $\T_{\bdy}$ to
$\T_{\bdy}'$. We have edges $e_1,e_2,e_3'$ of the triangulation
$\T_{\bdy}'$ in the boundary of $\bbb{T}'$.

Now, suppose $\gamma$ is a slope on $\bdy\bbb{T}$.  The slope
$\gamma$ is determined by its intersection with the edges $e_1$
and $e_2$, which form a basis for the first homology of
$\bdy\bbb{T}$. For a layering along the edge $e_3$, we still have
$e_1$ and $e_2$ edges in the boundary of the resulting torus,
$\bbb{T}' = \bbb{T}\cup_{e_3}\Delta$; and they also form a basis
for the first homology of $\bdy\bbb{T}'$. Furthermore, they meet
$\gamma$ exactly as they did before layering. So, if $\gamma'$ is
the unique slope in the boundary of $\bbb{T}'$ that meets $e_1$
and $e_2$ the same as $\gamma$, we say the slope $\gamma$ {\it
pushes through} the layering to the slope $\gamma'$.

The set of edges $\{e_1,e_2,e_3\}$ in $\T_{\bdy}$ is replaced by
the edges $\{e_1,e_2,e_3'\}$ in the one-vertex triangulation
$\T_{\bdy}'$ on the boundary of the solid torus $\bbb{T}'$. For
any orientation we choose on $e_1$ and $e_2$, we may orient $e_3$
so that either $e_3 = e_1 + e_2$ or $e_3 = e_1 - e_2$, with
respect to homology. Thus, $e_3'$ can be oriented so that
$e_3'=e_1 - e_2$, if $e_3 = e_1 +e_2$; or $e_3' = e_1 + e_2$, if
$e_3 = e_1 - e_2$. Choose an orientation on $\gamma$ and orient
$e_1$ and $e_2$ so that the oriented intersection numbers
$\langle\gamma,e_1\rangle =y_1$ and $\langle\gamma,e_2\rangle =
y_2$. Thus, if $e_3 = e_1 + e_2$, then $y_3 = y_1 + y_2$; and if
$e_3 = e_1 - e_2$, then $y_3 =| y_1 - y_2|$. After layering on
$e_3$, $e_1$ and $e_2$ remain in the boundary and the intersection
number of the new slope, $\gamma'$, with $e_1$ and $e_2$ remains
$y_1$ and $y_2$, respectively. It follows that in layering a
tetrahedron on the boundary of $\bbb{T}$ along the edge $e_3$ of
$\T$, the slope $\gamma$ pushes through the added tetrahedron to a
unique slope $\gamma'$ in the boundary of $\bbb{T}'$ and the
unique triple of numbers associated with $\gamma$ gives the unique
triple of numbers associated with $\gamma'$ as
$$\{y_1,y_2,\abs{y_1-y_2}\} \rightarrow \{y_1,y_2,y_1+y_2\},\hspace{.125 in} y_3
= |y_1 - y_2|$$ and
$$\{y_1,y_2,y_1+y_2\} \rightarrow \{y_1,y_2,|y_1-y_2|\}, \hspace{.125 in}y_3 = y_1 + y_2.$$

\vspace{.125 in}\noindent{\bf Example.} Suppose $\T$ is a one-vertex
triangulation of a solid torus, $\gamma$ is a slope on the boundary
of the solid torus and $\{3,8,11\}$ is the unique triple associated
with $\gamma$. If we layer $\T$ along the edge meeting $\gamma$
three times, then the slope $\gamma$  ``pushes through" to the slope
$\gamma'$ with triple $\{19,8,11\}$; similarly layering along the
edge meeting $\gamma$ eight times, the slope $\{3,8,11\}$ ``pushes
through" to $\{3,14,11\}$; and layering along the edge meeting
$\gamma$ eleven times, the slope $\{3,8,11\}$ ``pushes through" to
$\{3,8,5\}$.

\begin{figure}[htbp]

            \psfrag{a}{$a$}
            \psfrag{b}{$b$}
           \psfrag{c}{$c$}
            \psfrag{e}{$e$}
            \psfrag{1}{\small{$1$}}
            \psfrag{2}{\small{$2$}}
            \psfrag{3}{\small{$3$}}
            \psfrag{4}{\small{$3$}}
            \psfrag{5}{\small{$5$}}
            \psfrag{m}{meridian}
            \psfrag{f}{\begin{tabular}{c}
            slope\\
        of $e$\\
            \end{tabular}}
            \psfrag{k}{\begin{tabular}{c}
            preferred\\
        longitude\\
            \end{tabular}}
            \psfrag{l}{\begin{tabular}{c}
        solid\\
        torus\\
            \end{tabular}}
            \psfrag{A}{\begin{tabular}{c}
            {\Large(A)} slope\\
       of edge $e$\\
            \end{tabular}}
            \psfrag{B}{\begin{tabular}{c}
            {\Large(B)} preferred\\
        longitude\\
            \end{tabular}}
        \vspace{0 in}
        \begin{center}
        \epsfxsize=4 in
        \epsfbox{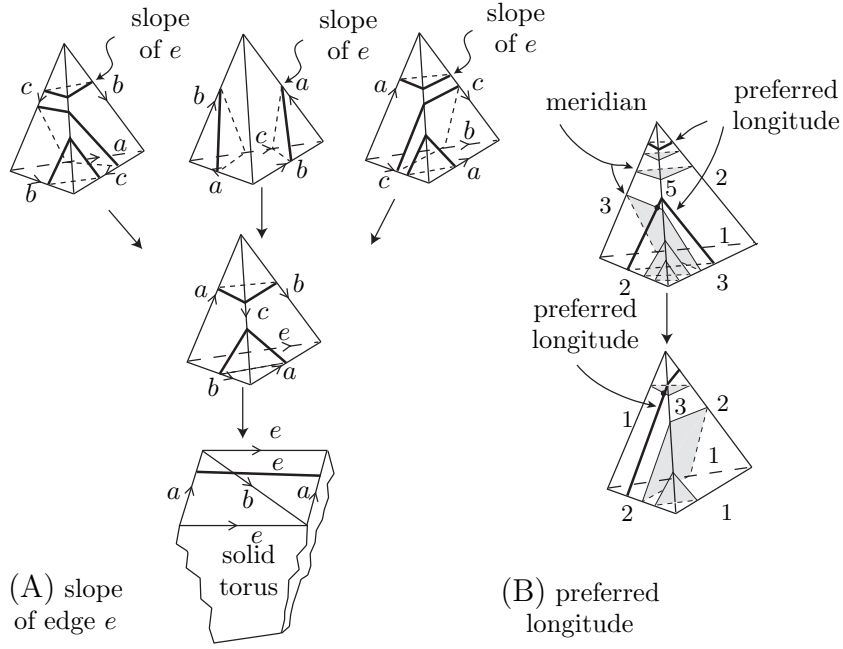}
        \caption{(A) gives slope of edge
        $e$ on boundary; (B) gives slope of thick edge,
        the preferred longitude, on the boundary.}
        \label{f-edge-slope}
        \end{center}

\end{figure}

\vspace{.15 in}\noindent {\it Slope of an edge.} If $\T$ is a
layered-triangulation of a solid torus, then each edge in $\T$ is
``unknotted," in the sense that it is isotopic into the boundary
torus; however, there are infinitely many curves in the boundary to
which it is isotopic through the solid torus -- each differing by a
Dehn twist about a curve representing the meridional slope. On the
other hand, notice that if $e$ is an edge in the triangulation $\T$,
then either $e$ is in the boundary, and so determines a {\it
preferred slope}, or at some unique level of the layering, the edge
$e$ becomes an interior edge. If this is the case, then there is a
unique ``push through" of the slope of the edge $e$, using the
layered-triangulation, which determines a {\it preferred slope} on
the boundary of the layered solid torus for the edge $e$. In Figure
\ref{f-edge-slope} we provide an example of the preferred slope on
the boundary for an interior edge $e$.

We call such a slope on the boundary, the \emph{slope of the edge
$e$} or, more generally, {\it the slope of an edge}. Hence, for a
solid torus with a fixed layered-triangulation there are a finite
number of special slopes on the boundary of the solid torus, the
(preferred) slopes of the edges of the layered-triangulation. In
particular, we have such a preferred slope in the boundary
determined by the unique ``thick edge", and so have a {\it preferred
longitude}. Since there is a unique meridian slope, a
layered-triangulation of a solid torus induces a {\it preferred
framing} (meridian/longitude coordinate system) on the boundary
torus. Below we will see that there is another very special
longitude on the boundary. Also, in Figure \ref{f-edge-slope}, we
exhibit that it is possible that two (or more) edges of a
layered-triangulation give the same slope on the boundary  (in
Figure \ref{f-edge-slope} the interior edge $e$ and the boundary
edge $f$ have the same slope). We use the phrase ``two edges have
the same slope," should this happen. It also may be the case that an
edge has the same slope as the meridian; if this is the case, we say
``an edge has the slope of the meridian". However, for minimal and
nearly-minimal layered-triangulations (defined below), distinct thin
edges give distinct slopes and an edge can not have the slope of the
meridian except for the minimal layered-triangulation extending
$0/1$.

\vspace{.15 in} \noindent{\it Existence: Extending triangulations
from the boundary of a solid torus.}  We wish to show that a
one-vertex triangulation on the boundary of a solid torus can be
extended to a layered-triangulation of the solid torus. In general,
given any triangulation on the boundary of a compact $3$--manifold,
it extends to a triangulation of the $3$--manifold; i.e., there is a
triangulation of the $3$--manifold so that its restriction to the
boundary is the given triangulation. In fact, the triangulation of
the $3$--manifold can be chosen so as not to add any vertices. The
distinction here is that any one-vertex triangulation on the
boundary of a solid torus extends to a {\it layered-triangulation}
of the solid torus.

To see this, suppose $\P$ is a triangulation of a $2$--manifold $S$
and $\gamma$ is a normal curve in $S$ (with respect to $\P$), we
call the cardinality of $\gamma\cap\P^{(1)}$ the {\it length of the
curve $\gamma$}, where $\P^(i)$ denotes the $i$--skeleton of $\P$,
and denote it by $L(\gamma)$. If $S$ is a torus, $\P$ is a
one-vertex triangulation of $S$ with edges $e_1, e_2$ and $e_3$, and
$\gamma$ is a normal simple closed curve in $S$, then choosing the
orientation on $\gamma$ and  $e_1, e_2, e_3$, as above, $L(\gamma)$
is the sum of the intersection numbers of $\gamma$ with the three
edges in $\P$: $L(\gamma) = y_1 + y_2 + y_3$, where $y_i =
\langle\gamma,e_i\rangle, i = 1,2,3$. So, for $S$ a torus and $\P$ a
one-vertex triangulation of $S$, if we perform a diagonal flip along
an edge $e$ in the triangulation $\P$, other than the one with
highest intersection number with $\gamma$, then the curve $\gamma$
has length strictly greater in the resulting triangulation than it
had in the triangulation $\P$. Alternatively, by doing a diagonal
flip along the edge of $\P$ having the highest intersection number
with $\gamma$, then in the resulting triangulation, $\gamma$ has
length strictly less than its length in $\P$, unless, with respect
to some ordering of the edges, we have $y_1 + y_2 =|y_1 - y_2|$
(i.e. one intersection coordinate is zero). The condition $y_1 + y_2
=|y_1 - y_2|$ means $\gamma$ is, in fact, disjoint from some edge
and is therefore the normal representative of that edge; $\gamma$
has intersection numbers $0,1,1$ and $L(\gamma) =2$. A diagonal flip
on either of the edges labeled ``$1$" does not change the
intersection numbers of $\gamma$ in the resulting triangulation from
those in $\P$.

Following our work, these observations were used in
\cite{jac-sedg-dehn}, to prove the following.

\begin{prop}\label{extend-to-layered} Suppose $\P$ is a one-vertex triangulation
of a torus $T$ and suppose $\mu$ is an arbitrary slope in $T$. Then
there is a layered-triangulation $\T$ of the solid torus and an
isomorphism from $\P$ to $\T_{\bdy}$, the triangulation induced by
$\T$ on the boundary of the solid torus, taking $\mu$ to the
meridional slope.\end{prop}

This result was applied in \cite{jac-sedg-dehn} to Dehn fillings of
cusped manifolds (compact, irreducible $3$--manifolds having a
single boundary component an incompressible torus). If a cusped
manifold $M$ is endowed with a one-vertex triangulation $\T$,
$\alpha$ is any slope in $\bdy M$, and $M(\alpha)$ is a Dehn filling
of $M$ along the slope $\alpha$, then by Proposition
\ref{extend-to-layered}, there is an extension of the triangulation
$\T$ to a layered-triangulation of the solid torus (depending on
$\alpha$) resulting in a triangulation $\T(\alpha)$ of $M(\alpha)$
that restricts to $\T$ on $M$. We discuss this more below where we
consider Dehn fillings.

The statement of Proposition \ref{extend-to-layered} is one
interpretation of extending a one-vertex triangulation on the
boundary of a solid torus to a triangulation of the solid torus and
is directly applicable to Dehn fillings. In the following theorem we
provide a different proof from that given in \cite{jac-sedg-dehn}, a
proof that adds some interesting combinatorial information about
layered-triangulations of solid tori.

\begin{thm}\label{th-extend-to-layered}
Suppose $\T_\bdy$ is a one-vertex triangulation on the boundary of a
solid torus. Then $\T_\bdy$ can be extended to a
layered-triangulation of the solid torus.
\end{thm}

\begin{proof} Let  $p,q,p+q$ be the unique triple of intersection
numbers of the meridional slope with the edges of the triangulation
$\T_{\bdy}$. We can assume the triple is not $0,1,1$ or $1,1,2$ and
that $0<p<q$. If we view the Euclidean algorithm for the pair $q,p$
as a subtraction rather than a division algorithm, we have:
$$q = (q-p)+p$$$$q-p = (q-2p) +p$$$$\vdots$$

where the general step is

$$r_{i-2} - (a_i-1)r_{i-1} = r_{i-1} + r_i, \hspace{.125 in}
r_i<r_{i-1}$$
$$\vdots$$

and, finally$$r_{n-1} = (r_{n-1}-1) +1$$$$r_{n-1}-1 = (r_{n-1}-2)
+1$$$$\vdots$$$$2 = 1+1.$$

Reversing the steps and starting with layering on the M\"obius band,
we pass from a triple $\{y_1, y_2,\abs{y_1-y_2}\}$ to the triple
$\{y_1,y_2,y_1+y_2\}$, eventually ending in the triple
$\{p,q,p+q\}$. It follows that we can layer on the M\"obius band to
get a layered-triangulation, $\T'$, of a solid torus having the
$p/q$--triangulation on its boundary. By Theorem \ref{pq-determine},
there is a homeomorphism from our original solid torus to the solid
torus we have built by layering, taking the triangulation
$\T_{\bdy}$ isomorphically onto the triangulation induced by $\T'$
on the boundary. Thus $\T'$ is an extension of $\T_{\bdy}$ to a
layered-triangulation of the solid torus.

In the cases where $\T_{\bdy}$ is a $0/1$-- or $1/1$--triangulation
on the boundary, then we have the two-tetrahedra extensions of these
triangulations given in Figure \ref{f-small-layered} and by applying
Theorem \ref{pq-determine}, we arrive at the desired conclusion in
these cases.\end{proof}

Since $\T_{\bdy}$ can be extended to  a layered-triangulation
without adding vertices, we can take the layered-triangulation
extending $\T_{\bdy}$ with the smallest number of tetrahedra. We
call such a triangulation a {\it minimal layered-triangulation} of
the solid torus extending $\T_{\bdy}$. If a layered triangulation of
the solid torus extends the $p/q$--triangulation on the boundary, we
call it a {\it $p/q$--layered-triangulation} and if it is the
minimal such layered-triangulation, we call it the {\it minimal
$p/q$--layered-triangulation}. We shall try to use ``layered" in
this way to distinguish the triangulation on the boundary from that
of the solid torus. We also find it convenient to refer to an
extension of a one-vertex triangulation on the boundary of a solid
torus to a layered triangulation of the solid torus as an
$\{x,y,z\}$--layered-triangulation, particularly, when we do not
wish to make a distinction as to the relative magnitude of the
intersection numbers.

\vspace{.15 in}\noindent {\it The $L$--graph: Classification of
layered-triangulations of a solid torus.}  From Theorem
\ref{pq-determine},  equivalence classes (up to homeomorphism of the
solid torus) of one-vertex triangulations on the boundary of the
solid torus are in one-one correspondence with reduced rational
numbers $p/q, 0\le p <q$, including  $0/1$ and $1/1$. The rational
number is uniquely determined by intersection numbers of the
meridional slope with the edges of the triangulation in the boundary
of the solid torus. Furthermore, if we add a tetrahedron by layering
along an edge of a one-vertex triangulation on the boundary of a
solid torus, we have a new solid torus and a new one-vertex
triangulation on its boundary. If the triangulation on the boundary
before layering has associated rational $p/q$, then after layering,
the new triangulation on the boundary has associated rational one of
$p/p+q, q/p+q$, or  $p/q-p$, if $p < q-p$ or $q-p/p$, if $q-p < p$.
We call these the {\it adjacency relations}.

\begin{figure}[htbp]

        \begin{center}
\epsfxsize = 4.5 in \epsfbox{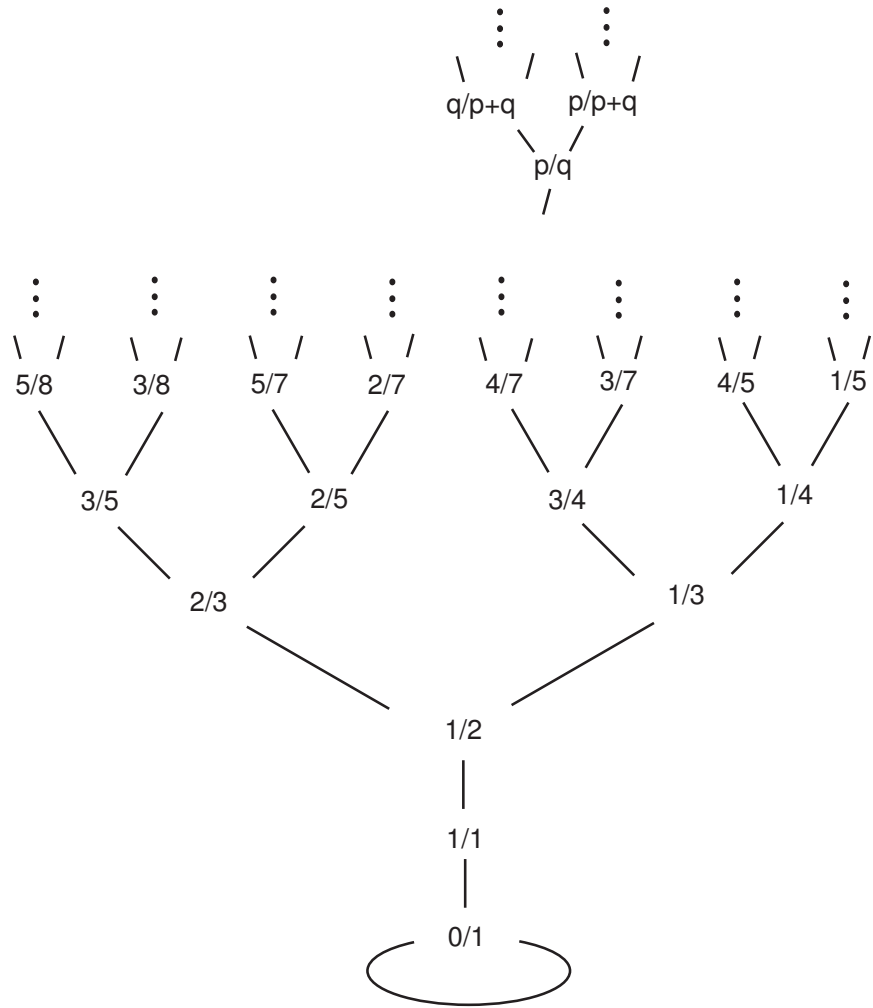}
 \caption{The start of the $L$--graph. The $0$--cells are in
one-one correspondence with one-vertex triangulations on the
boundary of the solid torus.} \label{f-Farey-Complex}
\end{center}
\end{figure}

From these adjacency relations, we construct a graph called the
$L$--graph. The $0$--cells of the $L$--graph are rational numbers
$p/q, 0\le p <q$, with $p,q$ relatively prime and they include the
values $0/1$ and $1/1$. There is a $1$--cell between two $0$--cells
(rationals) if and only if one $0$--cell is $p/q$ and the other
satisfies one of the adjacency relations given above. In Figure
\ref{f-Farey-Complex}, we exhibit some of the $L$--graph. All
$0$--cells except $1/1$ have index three and the $L$--graph minus
the closed edge loop at $0/1$ is a tree. The $0$--cells $1/1$ and
$0/1$ are each a bit of an anomaly.

We can think of a $1$--cell in the $L$--graph going between
$0$--cells labeled, say $p/q$ and $p/p+q$, as a tetrahedron in a
layered-triangulation of the solid torus, where two opposite edges
of the tetrahedron that correspond to the attaching edge and the
univalent edge meet the meridional disk of the solid torus in $q$
and $2p+q$ points; the remaining two edges meet the meridional disk
in $p$ and $p+q$ points. See Figure \ref{f-layered-torus-def}(B). We
have the following connection between layered-triangulations of the
solid torus and paths in the $L$--graph.

\begin{prop} Suppose $p/q$
is a $0$--cell in the $L$--graph (i.e., a one-vertex triangulation
on the boundary of a solid torus), then paths in the $L$--graph
starting at $p/q$ and ending at $1/1$ are in one-one correspondence
with layered-triangulations of the solid torus, extending the
$p/q$--triangulation on the boundary.
\end{prop}

By using the one-triangle M\"obius band as a (degenerate)
triangulation of the solid torus extending the $1/1$--triangulation,
and using the creased $3$--cell as  a (degenerate) triangulation of
the solid torus extending the $0/1$--triangulation, then the number
of $1$--cells in a path from $p/q$ to $1/1$ in the $L$--graph is the
same as the number of tetrahedra in the corresponding
layered-triangulation of the solid torus, extending the
$p/q$--triangulation on the boundary. The layered-triangulation
corresponding to the (unique) shortest path in the $L$--graph from
$p/q$ to $1/1$ is the minimal $p/q$--layered-triangulation. Again,
we point out that the minimal layered-triangulations extending the
$1/1$-- and $0/1$--triangulations are degenerate. If one is actually
wanting a triangulation of the solid torus in these two cases, then
one can use the two-tetrahedron triangulations for both or the
three-tetrahedron layered-triangulation for $0/1$ given in Figure
\ref{f-small-layered}. Notice that because of the loop in the
$L$--graph at the $0$--cell $0/1$, there is an infinite family of
non homotopic paths (keeping end points fixed) corresponding to
cycles about this loop. The homotopy classes of these loops give
rise to infinitely many  interesting $p/q$--layered-triangulations;
however, these will play essentially no role in our study; as we
will see below that these do not give $0$--efficient triangulations.

There is another observation that is useful. Namely, we can extend a
path in the $L$--graph that ends at $1/1$ by a closed loop that goes
to $1/2$ and then back to $1/1$. This is the same as opening the
one-tetrahedron solid torus or the creased $3$--cell along the
M\"obius band and inserting a Dehn filling along the curve with
intersection numbers $1,1,2$ by adding the two-tetrahedron
layered-triangulation of the solid torus extending $1/1$. Observe
that the $1,1,2$ curve is completely determined on the boundary; in
particular, if the edge labeled $2$ is the univalent edge upon
opening, then we have opened the creased $3$--cell and the original
path approached $1/1$ from $0/1$ in the $L$--graph. We refer to this
as {\it opening at $1/1$}. With this new terminology, the
two-tetrahedron layered-triangulation extending the
$1/1$--triangulation is an opening of the M\"obius band at $1/1$ and
the three-tetrahedron layered-triangulation of the solid torus
extending the $0/1$--triangulation is an opening of the creased
$3$--cell at $1/1$. Also, a path in the $L$--graph starting at $0/1$
and looping between $0/1$ and $1/1$ and eventually terminating at
$1/1$ is a creased $3$--cell; i.e., has the homotopy type of a solid
torus but is not a $3$--manifold. If the path transverses the loop
at $0/1$ or goes from $1/1$ to $1/2$ at any time, then we do have a
solid torus. Finally, the number of tetrahedra for the minimal
layered-triangulation of the solid torus extending the triangulation
$p/q$ on the boundary (with the exceptions of $0/1$ and $1/1$) is
$(\sum a_i) - 1$, where $a_i$ is the $i^{th}$ partial quotient in
the continued fraction expansion $p/q = (a_0;a_1; \ldots;a_n).$

Layered-triangulations are very nice ways to extend one-vertex
triangulations on the boundary of the solid torus and the minimal
layered-triangulations exhibit some great properties. We have the
following very compelling conjecture about triangulations of the
solid torus. We have tried but been unable to confirm this
conjecture.

\begin{conj2} The minimal  triangulation of the solid torus that extends
the $p/q$--triangulation on the boundary is the minimal
$p/q$--layered-triangulation of the solid torus.
\end{conj2}

A similar issue also comes up later as to whether
layered-triangulations of lens spaces are the minimal triangulations
of lens spaces. This question appears in the work of S. Matveev
\cite{matveev2}. Of course, we can ask the same question regarding a
minimal extension among all extensions of any triangulation from the
boundary of a $3$--manifold to the $3$--manifold.

In a minimal layered-triangulation of the solid torus, with the
exception of the two-tetrahedron triangulation extending $1/1$, the
univalent edge \emph{always} corresponds to the edge meeting the
meridional disk the greatest number of times. In the two-tetrahedron
layered-triangulation extending $1/1$, the univalent edge meets the
meridional disk just once. In the two-tetrahedron
layered-triangulation extending $0/1$, the univalent edge is labeled
$1$ and satisfies our observation; whereas, in the three-tetrahedron
layered-triangulation extending $0/1$, the univalent edge does not
meet the meridional disk at all (corresponds to the edge labeled
$0$). See Figure \ref{f-small-layered}. In each of these examples,
there are two edges labeled $1$; one is the thick edge, the other is
thin. We also observed above that in these two
layered-triangulations, the first layering is along the edge labeled
$3$, which is the univalent edge. In all other minimal
layered-triangulations, the layering is \emph{never} along the
univalent edge. More, on the implications of layering along
univalent edges, will be discussed below in our classification of
normal and almost normal surfaces in layered-triangulations of the
solid torus and of lens spaces and in the discussions of $0$-- and
$1$--efficient triangulations. These are simple, yet very important,
observations that make the understanding of minimal
layered-triangulations much easier than arbitrary (layered)
triangulations of the solid torus. In particular, we have that for a
minimal  layered-triangulation of a solid torus, the {\it slope of
an edge}, which determines a unique isotopy class on the boundary
and therefore a unique triple $r,s,r+s$, \emph{always} meets the
univalent edge the greatest number of times (of course, except for
the slope of the univalent edge itself).

Later, we define, in an analogous way, layered-triangulations of
handlebodies.

\section{Normal and almost normal surfaces in
layered-triangulations of a solid torus} We can classify all normal
and all almost normal surfaces in a \emph{minimal}
layered-triangulation of the solid torus. In order to apply our
study to efficient triangulations, we also classify all genus $0$
and genus $1$ normal surfaces in a (general) layered triangulation
of the solid torus.

\vspace{.125 in}\noindent{\bf Examples:}

\begin{figure}[htbp]
\psfrag{1}{\small{$1$}}\psfrag{2}{\small{$2$}}\psfrag{3}{\small{$3$}}\psfrag{4}{\small{$4$}}
\psfrag{5}{\small{$5$}}\psfrag{6}{\small{$6$}}
 \psfrag{A}{\begin{tabular}{c}M\"obius Band\\
           slope of\\ edge $2$\end{tabular}}\psfrag{B}{\begin{tabular}{c}Klein Bottle\\
            slope of\\ edge $4$\end{tabular}}\psfrag{C}{\begin{tabular}{c}Genus three\\nonorientable\\
             slope of \\edge $6$
           \end{tabular}}\psfrag{D}{\begin{tabular}{c}Klein Bottle\end{tabular}}
           \psfrag{E}{\begin{tabular}{c}M\"obius Band
           \end{tabular}}
{\epsfxsize = 4 in \centerline{\epsfbox{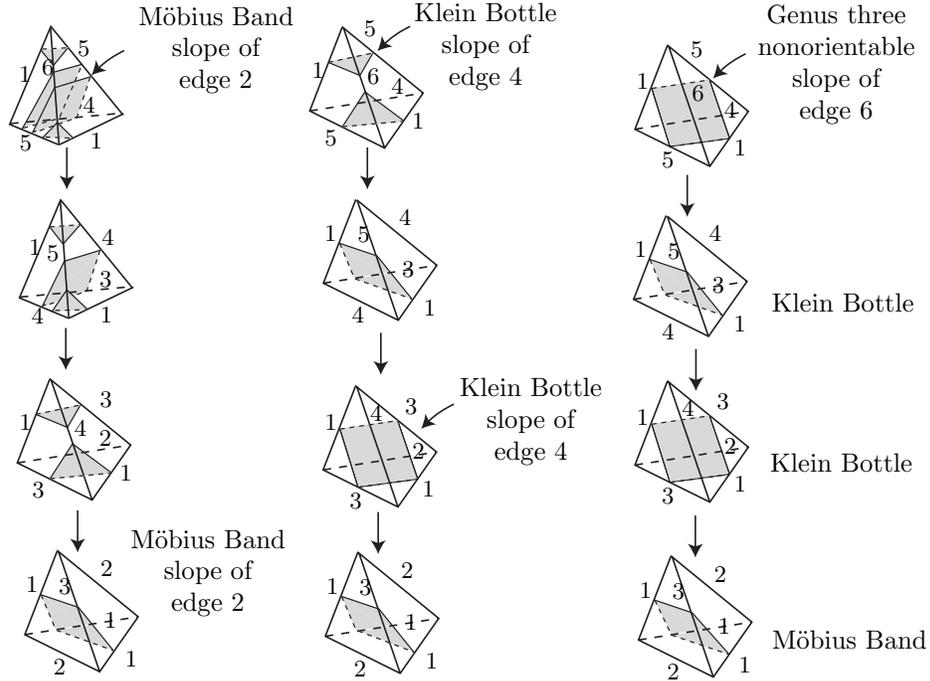}} }
\caption{Examples of nonorientable normal surfaces:  a M\"obius
Band, a Klein Bottle and a nonorientable surface of genus three.
Each has a single boundary curve having the slope of an even edge.}
\label{f-nonorientable}
\end{figure}

We first give some examples of normal surfaces in minimal
layered-triangulations of a solid torus and develop some terminology
for the various families of such normal surfaces. See Figures
\ref{f-merd-vert-disks}-\ref{f-nonorientable}
\begin{enumerate}

\item [A.] {\it Vertex-linking disk.} The unique vertex-linking
surface in a layered-triangulation of a solid torus is a disk.
Notice the general pattern for generating the induced triangulation
of a vertex-linking disk, Figure \ref{f-merd-vert-disks}.

\begin{figure}[htbp]

\psfrag{1}{\small{$1$}}\psfrag{2}{\small{$2$}}\psfrag{3}{\small{$3$}}
\psfrag{4}{\small{$4$}}\psfrag{5}{\small{$5$}}
            \psfrag{A}{\begin{tabular}{c}vertex-linking\\
            disks\end{tabular}}\psfrag{B}{\begin{tabular}{c}meridional\\
            disks\end{tabular}}
{\epsfxsize = 4.5 in \centerline{\epsfbox{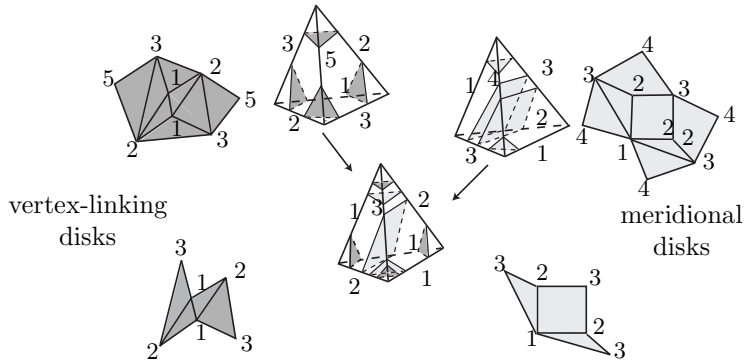}} }
\caption{Examples of normal vertex-linking and meridional disks.}
\label{f-merd-vert-disks}
\end{figure}

\vspace{.125 in} \item [B.] {\it Meridional disk.} There is a {\it
unique} normal meridional disk in a layered-triangulation of a solid
torus. While there may be other normal embedded disks having trivial
boundary in the boundary of the solid torus, there always is only
one normal meridional disk. Again, there is a general pattern for
generating the induced cell-decomposition (must have normal
quadrilaterals) of the meridional disk. See Figure
\ref{f-merd-vert-disks}. The growth (number) of quadrilaterals in
the induced cell decomposition  of the meridional disk in a minimal
layered-triangulation of the solid torus, as a function of the
number of layers, is like that of a Fibonacci sequence and so is
exponential.

\vspace{.125 in}
\begin{figure}[htbp]
\psfrag{1}{\small{$1$}}\psfrag{2}{\small{$2$}}\psfrag{3}{\small{$3$}}
\psfrag{5}{\small{$5$}}\psfrag{7}{\small{$7$}}\psfrag{9}{\small{$9$}}

            \psfrag{A}{\begin{tabular}{c}slope of\\
           edge $7$\end{tabular}}\psfrag{B}{\begin{tabular}{c}thin\\
           edge-linking\\annuli\end{tabular}}\psfrag{C}{\begin{tabular}{c}vertex-linking\\
          disks\end{tabular}}\psfrag{D}{\begin{tabular}{c}edge-linking\\
         annuli\end{tabular}}\psfrag{E}{\begin{tabular}{c}slope of\\
           edge $3$\end{tabular}}
{\epsfxsize = 4 in \centerline{\epsfbox{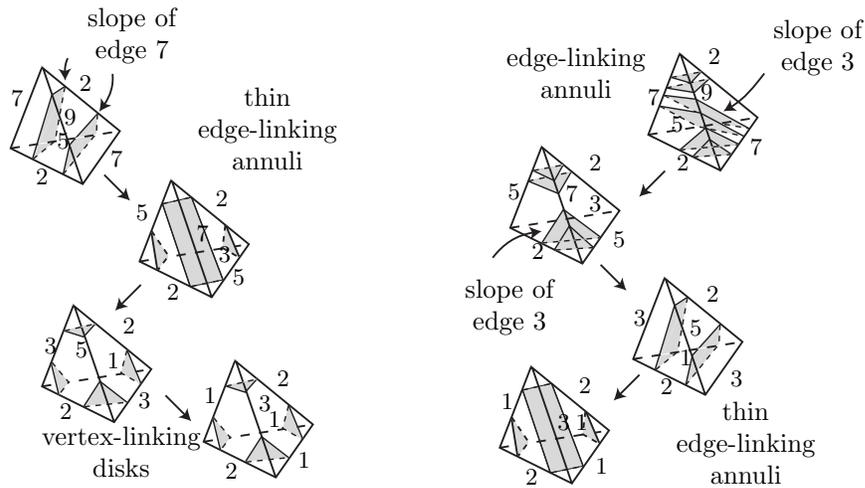}} }
\caption{Examples of normal edge-linking annuli.}
\label{f-edge-link-annuli}
\end{figure}

\vspace{.125 in} \item [C.] {\it Edge-linking annuli.} There are two
models of edge-linking annuli. One is a {\it thin} edge-linking
annulus and is a normal surface determined by taking the boundary of
a small regular neighborhood of a thin edge in the boundary of the
layered solid torus. Figure \ref{f-edge-link-annuli} gives such an
example where the thin edge-linking annulus is about the edge $7$ in
the boundary of the  minimal $2/7$--layered-triangulation of the
solid torus. The other example, in  Figure \ref{f-edge-link-annuli},
gives a normal annulus linking an interior edge; specifically, the
edge $3$ in the interior of $2/7$--layered-triangulation. Note that
the boundaries of these annuli have slopes equal to the slopes of
the edges they link ($7$ and $3$, respectively).

\vspace{.125 in}\item [D.] {\it Vertex-linking disk with thin
edge-linking tubes.} In Figure \ref{f-vert-disk-tube} there is an
example of the vertex-linking disk with a $1$-handle attach. A
$1$--handle is formed by taking a thin edge-linking tube about the
edge $3$ in the minimal $2/7$--layered-triangulation (as well as
already in the minimal $2/5$--layered-triangulation of the solid
torus). There can be numerous thin edge-linking tubes attached to
the vertex-linking disk.

\begin{figure}[htbp]

\psfrag{1}{\small{$1$}}\psfrag{2}{\small{$2$}}\psfrag{3}{\small{$3$}}\psfrag{4}{\small{$4$}}
\psfrag{5}{\small{$5$}}\psfrag{6}{\small{$6$}}\psfrag{7}{\small{$7$}}\psfrag{9}{\small{$9$}}
            \psfrag{A}{\begin{tabular}{c}vertex-linking\\
           disk with thin\\ edge-linking tube\end{tabular}}\psfrag{B}{\begin{tabular}{c}edge-linking annulus\\
            with thin\\ edge-linking tube\end{tabular}}\psfrag{C}{\begin{tabular}{c}thin\\edge-linking\\
             annuli\\
           \end{tabular}}\psfrag{D}{\begin{tabular}{c}thin edge-linking\\
             annuli with thin\\ edge-linking tube\end{tabular}}\psfrag{E}{\begin{tabular}{c}thin\\edge-linking\\
             annulus\\
           \end{tabular}} {\epsfxsize = 4 in
\centerline{\epsfbox{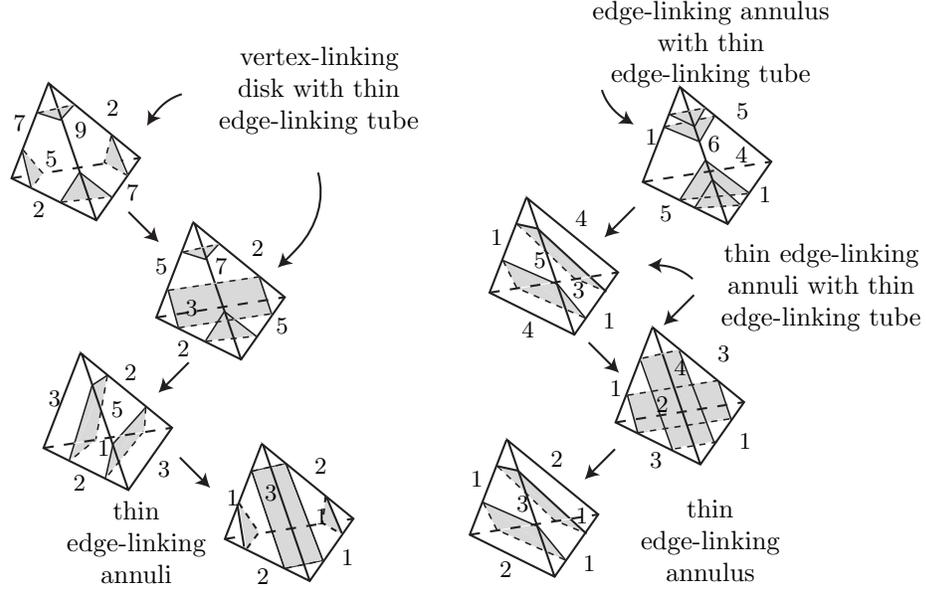}} } \caption{Examples of
a normal vertex-linking disk and a normal edge-linking annulus, each
with a thin edge-linking tube; on the left the tube is about the
edge labeled $3$ and on the right it is about the edge $2$.}
\label{f-vert-disk-tube}
\end{figure}

\vspace{.125 in}\item [E.] {\it Edge-linking annuli with thin
edge-linking tubes.} In Figure \ref{f-vert-disk-tube}  there also is
an example of an edge-linking annulus with a $1$--handle attached.
We already have this at level $1/3$,  where there is a thin
edge-linking annulus about the edge $4$ with a $1$--handle
determined by a thin edge-linking tube about the edge $2$. The
annulus is still thin edge-linking at level $1/4$; however, at level
$1/5$, we have an edge-linking annulus with boundary the slope of
the edge $4$; it is edge-linking but \emph{not} thin edge-linking.

\vspace{.125 in}\item [F.] {\it Nonorientable surfaces.} In Figure
\ref{f-nonorientable} we give three examples of embedded
(incompressible) nonorientable surfaces in the  minimal
$1/5$--layered-triangulation of the solid torus. One is a M\"obius
band and has boundary slope the slope of the edge $2$; one is  a
punctured Klein bottle and has boundary slope the slope of the edge
$4$; and one is a punctured genus three ($3$ crosscaps)
nonorientable surface and has boundary slope the slope of the edge
$6$. Note that the double of such a surface is an edge-linking
annulus with thin edge-linking tubes. Furthermore, in the double
there must be a tube linking every even ordered edge prior to the
edge of the edge-linking annulus. In general, there is a distinct
nonorientable surface for each ``even order" edge and its boundary
has slope equal to the slope of the edge it is associated with.
There are infinitely many nonorientable (incompressible) properly
embedded surfaces in a solid torus; however, for a fixed
layered-triangulation, there are only a finite number that are
isotopic to a normal surface. This displays one of the
discriminating properties of layered-triangulations.
\end{enumerate}

\subsection{Normal surfaces in minimal
layered-triangulations of a solid torus}  We  show that the model
examples given in the previous section represent {\it all} possible
types of embedded normal surfaces in minimal layered-triangulations
of a solid torus.

The classification of normal surfaces in minimal
layered-triangulations of a solid torus is made by induction on the
number of tetrahedra in the layering. We start with the
classification of normal surfaces in the one-tetrahedron solid
torus, which begins our induction step, and is the minimal
layered-triangulation extending the $1/2$--triangulation on the
boundary. We also need to consider, separately, the $2$
two-tetrahedron minimal layered-triangulations extending the $1/1$--
and the $0/1$--triangulations on the boundary.

\begin{lem}\label{one-tet-torus} The connected, embedded normal surfaces
in the one-tetrahedron triangulation of the solid torus are:

\begin{enumerate}
\item Vertex-linking disk, \item Meridional disk, \item Thin
edge-linking annulus with boundary slope the slope of the edge $3$,
\item M\"obius band with boundary slope the slope of the edge
$2$, and \item Thin edge-linking annulus with boundary slope the
slope of the edge $2$, which is also the double of the M\"obius
band.
\end{enumerate}

\end{lem}

\begin{figure}[htbp]
    \psfrag{a}{a}
\psfrag{b}{b} \psfrag{c}{c} \psfrag{d}{$x_1$} \psfrag{e}{$x_2$}
\psfrag{f}{$x_3$} \psfrag{g}{$x_4$} \psfrag{h}{$y_1$}
\psfrag{i}{$y_2$} \psfrag{j}{$y_3$}
            \psfrag{A}{$(a) : x_1 + y_3 = x_2 + y_2$}
            \psfrag{B}{$(b) : x_2 + y_1 = x_4 + y_1  $}
            \psfrag{C}{$(c) : x_4 + y_2 = x_3 + y_3  $}
        \vspace{0 in}
        \begin{center}
        \epsfxsize=4.5 in
        \epsfbox{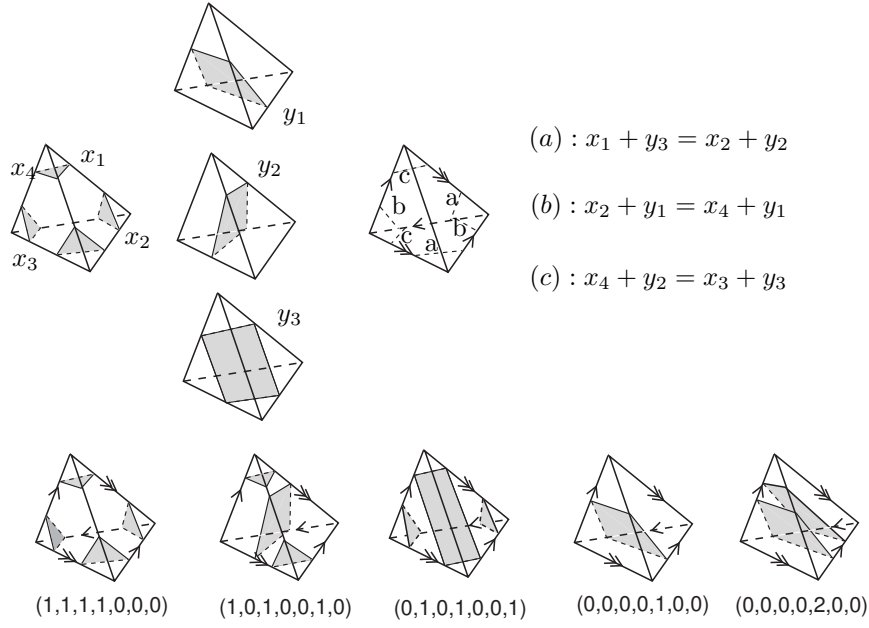}
        \caption{Embedded normal surfaces in the one-tetrahedron solid torus.}
        \label{f-normal-one-tet}
        \end{center}

\end{figure}

\begin{proof}
This is an easy exercise in using normal coordinates and face
matching equations. There are $7$ variables and $3$ matching
equations. We give the equations and their solutions in Figure
\ref{f-normal-one-tet}. Also, one can use quadrilateral coordinates.
In this case, there are three variables and no edge equations (there
are no interior edges). There is precisely one fundamental solution
for each quadrilateral; these solutions and  the vertex-linking disk
give the embedded fundamental solutions. Solution (5) is the double
of solution (4).\end{proof}

 We now use the method of
quadrilateral equations to classify the normal surfaces in the
two-tetrahedron solid torus extending the $0/1$--triangulation on
the boundary of the solid torus.

\begin{lem}\label{0-1-1} The connected, embedded normal surfaces in the
two-tetrahedron triangulation of the solid torus extending the
$0/1$--triangulation are:

\begin{enumerate}
\item Vertex-linking disk, \item Thin edge-linking annulus about
the thin edge labeled $1$ in the boundary, \item M\"obius band with
boundary slope the slope of the interior edge labeled $2$,
\item Double of the M\"obius band, which is an edge-linking annulus
about the interior edge labeled $2$,
\item Meridional disk, \item Vertex-linking disk with a thin
edge-linking tube about the interior edge labeled $2$, \item Thin
edge-linking annulus about the thin edge labeled $1$ in the boundary
with a thin edge-linking tube about the interior edge labeled
$2$,\item Punctured Klein bottle with boundary slope the slope of
the edge labeled $0$, \item Thin edge-linking annulus about the edge
labeled $0$ in the boundary,\item Edge-linking annulus with a tube
about the edge labeled $2$ and with boundary slope that of the edge
labeled $0$, which is the double of the punctured Klein bottle.
\end{enumerate}\end{lem}

\begin{figure}[htbp]
    \psfrag{0}{\small{$0$}}
\psfrag{1}{\small{$1$}} \psfrag{2}{\small{$2$}}
\psfrag{3}{\small{$e$}}\psfrag{4}{\small{$1'$}}\psfrag{a}{\small{$x_1$}}
\psfrag{b}{\small{$x_2$}} \psfrag{c}{\small{$x_3$}}
\psfrag{d}{\small{$x_4$}} \psfrag{e}{\small{$y_1$}}
\psfrag{f}{\small{$y_2$}}
\psfrag{g}{\small{$y_3$}}\psfrag{h}{\small{$x_5$}}
\psfrag{i}{\small{$x_6$}} \psfrag{j}{\small{$x_7$}}
\psfrag{k}{\small{$x_8$}} \psfrag{l}{\small{$y_4$}}
\psfrag{m}{\small{$y_5$}}
\psfrag{n}{\small{$y_6$}}\psfrag{x}{\footnotesize{x}}\psfrag{y}{\footnotesize{y}}
\psfrag{z}{\footnotesize{z}}\psfrag{L}{layered}
\psfrag{A}{\Large{(3)}}\psfrag{B}{\Large{(7)}}
           \psfrag{C}{\begin{tabular}{ll}$(e):$  & $ x_2 + y_3 = x_2 +
            y_2$\\
           &$y_2 = y_3$
            \end{tabular}}
           \psfrag{M}{\begin{tabular}{ll}$y_1 =0$ & $y_4 =0$ \\
          $y_5 = 1$ & $y_6 = 0$
           \end{tabular}}

            \psfrag{N}{\begin{tabular}{ll}$y_1 =1$ &  $y_4 =1$ \\
 $y_5 = 0$ & $y_6 = 0$
            \end{tabular}}

        \vspace{0 in}
        \begin{center}
        \epsfxsize=4.5 in
        \epsfbox{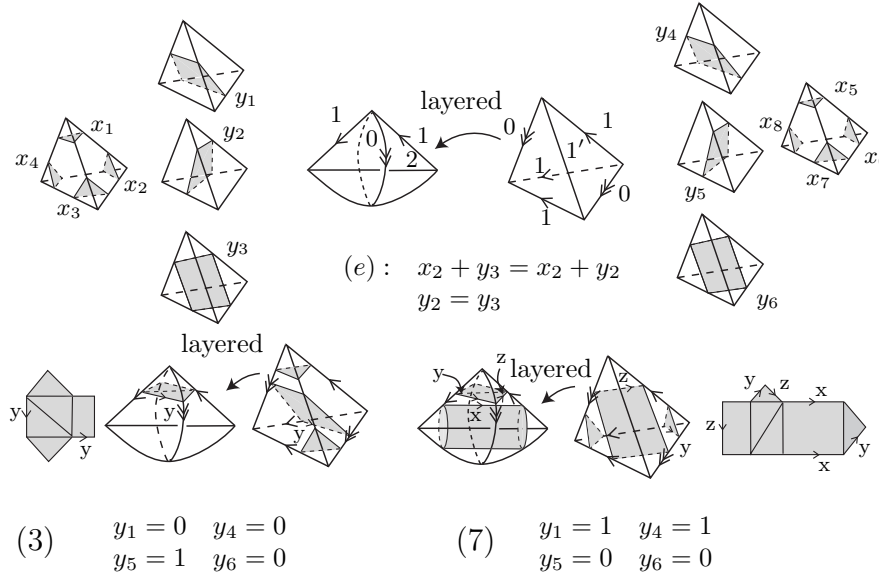}
        \caption{Embedded normal surfaces in the two-tetrahedron
         layered extension of $0/1$.  $(3)$ is the M\"obius
        band with boundary the slope of the edge $2$; and $(7)$
         is a thin edge-linking annulus with a thin edge-linking tube
          about the edge $2$.}
        \label{f-normal-0-1-1}
        \end{center}
        \end{figure}

\begin{proof}  Since there are only two tetrahedra, using quadrilateral
coordinates there are 6 variables (the 6 quad types) and only 1
equation, coming from the one interior edge of the triangulation.
Using the notation from Figure \ref{f-normal-0-1-1}, the variables
can be denoted $y_1,y_2,y_3,y_4,y_5,$ and $y_6$; the one equation is
$y_2 = y_3$. Thus if $y_2 (y_3)\ne 0$, we have a singular surface.
So, we can organize the possibilities as: $y_1 = 0$ or $y_1\ne 0$
and at most one of $y_4\ne 0, y_5\ne 0$ or $y_6\ne 0$. The easiest
way to find the coordinates for the normal triangles is from the
matching equations.\end{proof}

Before we do the general classification, we establish some
terminology generalizing that used above. Suppose $e_1,\ldots,e_n$
are distinct interior edges of a layered-triangulation of a solid
torus and the boundary of an arbitrarily small regular neighborhood
of the wedge of the $e_i's$ is a normal surface. Not every wedge of
distinct interior edges has this property and, in particular, none
of the $e_i's$ are the thick edge. The normal, regular neighborhood
can be thought of combinatorially as removing the collection of
triangles that meet the $e_i's$ from the vertex-linking disk and
replacing them with the quadrilaterals about the $e_i's$. We call
such a surface a {\it vertex-linking disk with thin edge-linking
tubes} or, simply, a {\it vertex-linking disk with tubes}. Notice
that each of the tubes is ``unknotted" in the solid torus in the
sense that its core simple closed curve (an interior edge) is
isotopic in the solid torus to a curve in the boundary of the solid
torus. Examples are given above in the layered-triangulation
extending the $2/7$--triangulation (Example D and Figure
\ref{f-vert-disk-tube}) and in item (6) of Lemma \ref{0-1-1} in the
classification of surfaces in the two-tetrahedra triangulation of
the solid torus extending the $0/1$--triangulation.

The next general form is the edge-linking annulus.  Suppose we have
a layered-triangulation of a solid torus, extending the
$p/q$-triangulation  on the boundary; hence, there is a
corresponding path in the $L$--graph beginning at $p/q$ and ending
at $1/1$. For any $0$--cell $r/s$ in this path, a segment of the
path from $1/1$ to $r/s$ is called an {\it initial segment}. This
initial segment determines a subcomplex of the original
layered-triangulation and is a layered-triangulation of a solid
torus extending the $r/s$-triangulation on the boundary of a solid
torus. (See Figure \ref{f-vert-disk-tube} above for the initial
segment corresponding to the vertex $1/3$ in the minimal
layered-triangulation extending the $1/5$--triangulation.) Now, if
$e$ is a thin edge in the boundary of a layered-triangulation of a
solid torus, the boundary of an arbitrarily small regular
neighborhood of $e$ is a normal surface. It can be thought of
combinatorially as removing the collection of triangles that meet
$e$ from the vertex-linking disk and replacing them with the
quadrilaterals about $e$. We call such a surface a {\it thin
edge-linking annulus}. Next, suppose $\T$ is a layered-triangulation
of a solid torus and $\T'$ is an initial segment of $\T$. If $A'$ is
a thin edge-linking annulus in $\T'$, then we can push the boundary
of $A'$ through the layered-triangulation of $\T$ to get a new
normal, properly embedded annulus $A$ in $\T$ (again see C above and
Figure \ref{f-vert-disk-tube}). We call such a normal surface an
{\it edge-linking annulus}, as opposed to a {\it thin} edge-linking
annulus.

Now, suppose $e$ is an edge in the boundary of a
layered-triangulation of a solid torus and $e_1,\ldots,e_n$ are
edges in the interior of the layered solid torus so that the
boundary of an arbitrarily small regular neighborhood of the wedge
of $e$ and the edges $e_i, 1\le i\le n$ is a normal surface. Such a
normal surface can be combinatorially  described as removing the
collection of triangles that meet $e$ and any $e_i$ from the
vertex-linking disk and replacing them with the quadrilaterals about
$e$ and each $e_i$. We call such a surface a {\it thin edge-linking
annulus with thin edge-linking tubes} or, simply, a {\it thin
edge-linking annulus with tubes} (see Examples  E and F above and in
Figure \ref{f-vert-disk-tube} and Figure \ref{f-nonorientable}, as
well as item (9) of Lemma \ref{0-1-1}).

Finally, suppose $\T$ is a layered-triangulation of a solid torus
and $\T'$ is an initial segment of $\T$ and $A'$ is a thin
edge-linking annulus with tubes in $\T'$. If we push $A'$ through
the layered-triangulation $\T$ to get the normal surface $A$, then
$A$ is an {\it edge-linking annulus with tubes} (as opposed to a
\emph{thin} edge-linking annulus with tubes). Notice that the
boundary slope of the boundary of the edge-linking annulus (with
tubes) is the slope of an edge; it is the slope of the edge of the
thin edge-linking annulus $A'$. Also, there are no tubes possible
between $\bdy A'$ and $\bdy A$.  An example is given in Figure
\ref{f-vert-disk-tube}; in that example we can take $A'$ in the
initial segment corresponding to the $1/3$--layered-triangulation,
where $A'$ is the thin edge-linking annulus about the edge labeled
$4$ with a (thin edge-linking) tube about the edge labeled $2$. The
surface  $A$ is the ``push through" of $A'$ and it is an
edge-linking annulus with a thin edge-linking tube in the
$1/5$--layered-triangulation. We sometimes just use ``edge-linking
annulus with tubes" including the possibility of a thin edge-linking
annulus with tubes; however, they are different in nature. The thin
edge-linking annulus with tubes is the boundary of a regular
neighborhood of a wedge of edges in the layered-triangulation;
whereas, the more general edge-linking annulus with tubes may not
be. A thin edge-linking annulus also has boundary slope the slope of
an edge in the boundary; whereas, an edge-linking annulus that is
not a thin edge-linking annulus has boundary slope the slope of an
edge from the interior of the layered-triangulation. There is
another very important fact about edge-linking annuli. In a
layered-triangulation that is not minimal, there may be distinct
edges determining the same slope on the boundary of the solid torus.
In such a case there can be distinct (not normally isotopic)
edge-linking annuli with the same boundary slope. Even though the
annuli are not normally isotopic, they are isotopic. This phenomenon
can not happen in a minimal layered-triangulation of a solid torus.
Any edge-linking annulus, $A$, is parallel in the solid torus to an
annulus in the boundary containing the vertex. We call this annulus
in the boundary of the solid torus the {\it companion annulus} (to
$A$). The annulus complementary to the companion annulus in the
boundary of the layered solid torus is called the {\it complementary
annulus} (to $A$).

We now have the classification of normal surfaces in minimal
layered-triangula-tions of the solid torus.

\begin{thm}\label{normal} A connected, embedded normal surface in
a minimal layered-trian-gulation of the solid torus is normally
isotopic to one of:
\begin{enumerate}
\item The vertex-linking disk, \item A vertex-linking disk with
tubes, \item The (unique) meridional disk, \item An edge-linking
annulus, having boundary slope the slope of the edge, \item An
edge-linking annulus with tubes, having boundary slope the slope
of an edge, or \item An incompressible, nonorientable surface,
having boundary slope the slope of an even ordered edge.
\end{enumerate}
\end{thm}

\begin{proof} We shall use induction for our proof and assume we
are starting with the one-tetrahedron solid torus. Furthermore, we
will consider the minimal $1/1$--layered-triangulation separately;
hence, we never layer on the univalent edge.

From Lemma \ref{one-tet-torus}, we know exactly the connected,
embedded normal surfaces in the one-tetrahedron triangulation of the
solid torus. These satisfy our conclusions above. Our induction
hypothesis is then that in a minimal layered-triangulation $\T_t$
with $t$ tetrahedra, $(t\geq 1)$, the list (1)--(6) in the
conclusion of Theorem \ref{normal} is the complete list of
connected, embedded normal surfaces; furthermore, each slope of an
edge corresponds to a unique edge. This is a central feature of a
minimal layered-triangulation. Let $\T_{t+1} = \T_t \cup_{e}
\td{\Delta}_{t+1}$ be a minimal layered-triangulation having $t +1$
tetrahedra. Furthermore, suppose the minimal layered-triangulation
$\T_t$ extends the $p/q$--triangulation  on the boundary of the
solid torus. Note also, we may assume that the edge $e$ is {\it not}
a univalent edge.

If $F_{t+1}$ is a normal surface in $\T_{t+1}$, then it meets
$\T_t$ in a normal surface, say $F_t$. Hence, each component of
$F_t$ must fall into one of the categories (1)--(6) in the
conclusion of Theorem \ref{normal}. Furthermore, we have $F_{t+1}$
is obtained from $F_t$ by adding cells that lift to normal
triangles and quads in $\td{\Delta}_{t+1}$. There are two
possibilities.

\vspace{.15 in} \noindent {\bf Case A.} We push $F_t$ through
$\Delta_{t+1}$. Then $F_{t+1}$ is homeomorphic with $F_t$ and either
$F_{t+1}$ is a vertex-linking disk with tubes and has trivial
boundary, $F_{t+1}$ is a push through of an edge-linking annulus
with tubes and has boundary the slope of an edge, or $F_t$ is the
meridional disk and $F_{t+1}$ is also the meridional disk. In the
last two cases, we have that if the boundary of $F_t$ does not have
the slope of the univalent edge, then it meets the univalent edge
maximally and therefore the boundary of $F_{t+1}$ meets the
univalent edge maximally. If the boundary of $F_t$ has the slope of
the univalent edge, then it is still true that the slope of the
boundary of $F_{t+1}$  meets the univalent edge maximally. A push
through in a minimal layered-triangulation can not have boundary
with slope equal to that of the univalent edge.

\begin{figure}[h]
    \psfrag{p}{\footnotesize{$p$}}
\psfrag{q}{\footnotesize{$q$}} \psfrag{a}{\footnotesize{$p+q$}}
\psfrag{B}{{\small layered}}\psfrag{C}{``push
through"}\psfrag{D}{annulus}

           \psfrag{E}{\begin{tabular}{c}
            vertex-disk\\
          with tube
            \end{tabular}}
           \psfrag{F}{\begin{tabular}{c}annulus\\with tube
           \end{tabular}}

            \psfrag{G}{\begin{tabular}{c}not\\orientable
            \end{tabular}}
            \psfrag{H}{\begin{tabular}{c}banding quadrilateral
           \end{tabular}}
{\epsfxsize = 4 in \centerline{\epsfbox{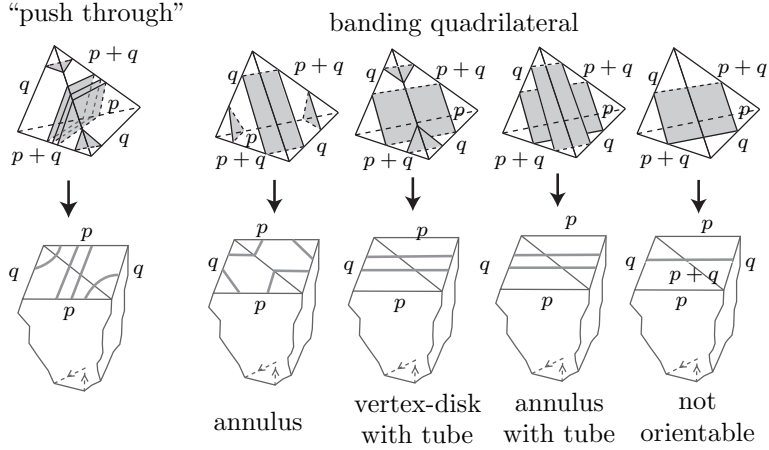}} }
\caption{Examples of ``pushing through" a surface and of adding a
banding quadrilateral.} \label{f-banding-quad}
\end{figure}

\vspace{.15 in} \noindent {\bf Case B.} We add ``banding quads" to
$F_t$ in $\td{\Delta}_{t+1}$. A {\it banding quadrilateral} is a
normal quadrilateral in $\td{\Delta}_{t+1}$ separating the edge that
becomes the univalent edge in $\T_{t+1}$ from its opposite edge,
which is the edge in $\td{\Delta}_{t+1}$ attached to $e$ (see Figure
\ref{f-banding-quad} and \cite{jac-sedg-dehn}).

In this case, we need to consider the various possibilities.
First, in order to add a banding quad, the boundary components of
$F_t$ must be either vertex-linking curves or have the slope of
the edge $e$, along which we are attaching the tetrahedron
$\td{\Delta}_{t+1}$. Hence, by induction a component of $F_t$ is a
vertex-linking disk  (possibly with thin edge-linking tubes), a
thin edge-linking annulus (possibly with thin edge-linking tubes),
or a nonorientable surface. In this last case, the edge $e$ meets
the meridional disk an even number of times. Since we are assuming
$F_{t+1}$ is connected, we have that $F_t$ must be connected and
in the case $F_t$ is an annulus (possibly with thin edge-linking
tubes) a single banding quad must attach to distinct boundary
components of $F_t$.

If the boundary of $F_t$ is the vertex-linking curve in $\T_t$, then
we can add at most one banding quad. Furthermore, in this case by
our induction hypothesis, $F_t$ is either the vertex-linking disk in
$\T_t$ or it is the vertex-linking disk with thin edge-linking
tubes. In the first case, adding the banding quad (along with two
normal triangles, which push through the remaining part of the
boundary of $F_t$), we get that $F_{t+1}$ is a thin edge-linking
annulus with boundary slope that of the univalent edge in
$\T_{t+1}$. If $F_t$ is the vertex-linking disk with tubes, then
adding the banding quad (along with two normal triangles that push
through the rest of the boundary of $F_t$), we get that $F_{t+1}$ is
a thin edge-linking annulus with tubes and with boundary slope that
of the univalent edge in $\T_{t+1}$.

If the slope of the boundary of $F_t$ is the same as the edge $e$
in the boundary of $\T_t$, then again by our induction hypothesis,
we have two possibilities: the boundary of $F_t$ is a single curve
and $F_t$ is the unique nonorientable surface with boundary slope
the edge $e$ in the boundary of $\T_t$ (and this edge has even
order) or the boundary of $F_t$ has two curves.

In the first situation, we have only one possibility for adding a
banding quad and we have that $F_{t+1}$ is the unique nonorientable
surface with boundary slope the slope of the univalent edge in
$\T_{t+1}$ (and the univalent edge has even order for its
intersection number with the meridional disk). The genus of
$F_{t+1}$ is one more than the genus of $F_t$.

If $F_t$ has two curves in its boundary, then it is either a {\it
thin} edge-linking annulus or a {\it thin} edge-linking annulus with
tubes and, in either case, with boundary slope that of the edge $e$.
We can add one band (and push through the remaining part of the
boundary of $F_t$), then $F_{t+1}$ is a vertex-linking disk with
tubes (it must have at least a tube about the edge $e$). We can add
two bands, then $F_{t+1}$ is a thin edge-linking annulus with tubes
(it must have at least a tube about the edge $e$) and with boundary
slope that of the univalent edge of $\T_{t+1}$. This completes the
induction step.

This leaves only determining the normal surfaces in the
two-tetrahedron, minimal, non degenerate layered-triangulation of
the solid torus extending the $1/1$--triangulation. See Figure
\ref{f-small-layered}. We can easily adapt the preceding induction
proof to this simple case of a layering on the one-tetrahedron solid
torus. By pushing surfaces in the one-tetrahedron torus through, we
get the vertex-linking disk, meridional disk, a M\"obius band with
boundary slope that of the edge $2$ in the boundary, its double,
which is a thin edge-linking annulus about the edge $2$, and an
edge-linking annulus with boundary slope that of the interior edge
$3$. We can band the vertex-linking disk in the one-tetrahedron
solid torus , getting a thin edge-linking annulus on the univalent
edge. We can also band on the thin edge-linking annulus about the
edge $3$ in the one-tetrahedron solid torus. Using just one band, we
get a vertex-linking disk with a tube about the interior edge $3$;
using two banding quads, we get a thin edge-linking annulus about
the univalent edge along with a tube about the interior edge $3$.
This completes the proof of Theorem \ref{normal}.\end{proof}

Notice we indicated three methods to find normal surfaces in
layered-triangula-tions: normal coordinates and face-matching
equations, quadrilateral coordinates and edge-matching equations,
and an induction argument, using the layered structure of the
triangulation. The latter is the method of choice for
layered-triangula-tions.

We have from the proof of Theorem \ref{normal} that:
\begin{itemize}
\item There are no closed normal surfaces in a minimal layered-triangulation of a solid torus (this remains true even if the
layered-triangulation is not minimal - see below). \item  There is a
unique normal meridional disk in a minimal  layered-triangulation of
a solid torus and the slope of its boundary meets the univalent edge
a maximum number of times. It remains true that there is a unique
meridional disk even if the layered-triangulation is not minimal. We
think of the nontrivial normal spanning arc in the one-triangle
M\"obius band and the wedge of two triangular cones in the creased
$3$--cell as ``meridional disks" for the (degenerate) $1/1$-- and
$0/1$--triangulations, respectively.
\item There is a unique trivial normal disk, the vertex-linking
disk, in a minimal layered-triangulation (however, this is not
necessarily true if the layered-triangulation is not minimal - see
below). \item The boundary slope of any normal surface other than
the meridional disk and the vertex-linking disk is a slope of an
edge and therefore, meets the univalent edge a maximum number of
times (unless the boundary slope is the slope of the univalent edge
itself). \item There are infinitely many incompressible non
orientable surfaces embedded in a solid torus; however, a
layered-triangulation only has finitely many normal non orientable
surfaces, one for each even order edge in the layered-triangulation.
Furthermore, the slope of the boundary of a non-orientable normal
surface is the slope of the even ordered edge with which it can be
associated. Below we give a recursive function that determines the
genus of a normal non-orientable surface in a minimal
layered-triangulation of the solid torus.
\item We repeat that in the case of an edge-linking annulus with
tubes, all of the tubes are about edges coming in the layering
before the edge about which the annulus links. \item If $\lambda$
and $\gamma$ are two slopes on a torus, let $d(\lambda, \gamma)$
denote the distance between them. For $F$ a vertex-linking disk with
thin edge-linking tubes about edges having slopes $\lambda$ and
$\gamma$, then $d(\lambda,\gamma) \geq 2$; similarly, if $F$ is an
edge-linking annulus about the edge $e$ with slope $\lambda$ and $F$
has thin edge-linking tubes about the edges $e'$ and $e''$ with
slopes $\gamma'$ and $\gamma''$, respectively, then
$d(\lambda,\gamma')$, $d(\lambda,\gamma'')$ and
$d(\gamma',\gamma'')$ are all $\geq 2$.
\end{itemize}

\subsection{Genus $0$ and $1$ normal surfaces in a
(general) layered-triangulation of a solid torus} \vspace{.25 in}

Below it will be necessary for us to understand genus $0$ and genus
$1$ normal surfaces embedded in a (general) layered-triangulation of
a solid torus in order to study $0$--efficient and $1$--efficient
layered-triangulations of lens spaces.

\begin{figure}[h]

\psfrag{a}{\small a}\psfrag{b}{\small b}\psfrag{c}{\small
c}\psfrag{d}{\small d}\psfrag{e}{\small{$\mathbf{e}$}}
\psfrag{f}{\small{$\mathbf{e'}$}} \psfrag{A} {\begin{tabular}{c}
$\mathbf{e}$ and $\mathbf{e'}$\\
 same slope\end{tabular}} \psfrag{B}{\begin{tabular}{c}
$\mathbf{e'}$ slope\\
of meridian\end{tabular}} \psfrag{C}{\begin{tabular}{c}Slope of\\
meridian\end{tabular}}{\epsfxsize = 4 in
\centerline{\epsfbox{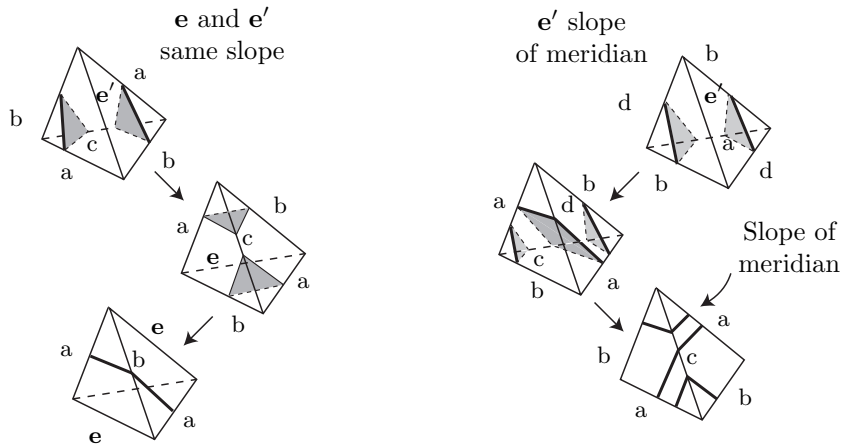}} } \caption{Examples of two
edges ($e$ and $e'$) with the same slope; and an edge ($e'$) with
the slope of a meridional disk.} \label{f-same-slope}
\end{figure}

We are only interested in orientable surfaces of genus $0$ and genus
$1$; as might be expected, some of these occur as the doubles of
embedded normal M\"obius bands or embedded normal punctured Klein
bottles; we will point out these possibilities in our classification
of genus $0$ and genus $1$ embedded normal surfaces.

\vspace{.125 in}\noindent{\bf Examples.} First, we give some
examples  of genus $0$ and genus $1$ normal surfaces in general
layered-triangulations of a solid torus; some are new and are not
seen in minimal layered-triangulations.

\vspace{.125 in} \noindent {\it Genus $0$ normal surfaces}: A genus
$0$ normal surface is obtained in $\T_{t+1}$ by pushing through a
genus $0$ normal surface from $\T_t$ or by adding a banding quad in
$\td{\Delta}_{t+1}$ to a disk or to two copies of a disk.

\begin{enumerate}
\item[(1)--(3)] We have the possibility of any   of the  genus $0$ normal surfaces
that can occur in {\it minimal} layered-triangulations and  are
classified in Theorem \ref{normal}. (1) the {\it vertex-linking
disk}, (2) the {\it meridional disk}, and (3) {\it edge-linking
annuli} (both thin edge-linking  and  non-thin edge-linking annuli).

\item [(4)] In addition, we have the possibility for two new genus $0$ normal surfaces, if
the edge $e$ in $\T_t$, along which we attach $\td{\Delta}_{t+1}$,
has the same slope as the meridian.  To see this, suppose the edge
$e$ has the same slope as the meridian. Let $P_t$ be two copies of
the meridian disk ($P_t$ has two components). First, add a banding
quad in $\td{\Delta}_{t+1}$ that runs between the distinct
components of $P_t$; such a banding quad, along with two triangles
in $\td{\Delta}_{t+1}$, gives an embedded normal disk in $\T_{t+1}$
that is trivial. Its boundary is the vertex-linking curve in the
boundary of $\T_{t+1}$ but it is not the vertex-linking disk. Hence,
it is (4A) a {\it non vertex-linking trivial disk}. See Figure
\ref{f-non-vert-edge-linking}(A).

\begin{figure}[h]
\psfrag{0}{\footnotesize{$0$}}\psfrag{1}{\footnotesize{$1$}}\psfrag{2}{\footnotesize{$2$}}
\psfrag{3}{\footnotesize{$1'$}}\psfrag{A}{\Large (A)}
\psfrag{B}{\Large (B)}\psfrag{J}{\begin{tabular}{c} \small{non
vertex-linking,}\\\small{trivial disk}\end{tabular}}
\psfrag{K}{\begin{tabular}{c}\small {non
edge-linking}\\\small{annulus}\end{tabular}}
\psfrag{C}{\begin{tabular}{c}\footnotesize{ two copies
of}\\\footnotesize{meridional
disk}\end{tabular}}\psfrag{D}{\begin{tabular}{c}\footnotesize{ slope
of}\\\footnotesize{edge
$0$}\end{tabular}}\psfrag{E}{\begin{tabular}{c}\small{
layered}\\\small{solid}\\\small{torus}\end{tabular}}
 {\epsfxsize = 3 in
\centerline{\epsfbox{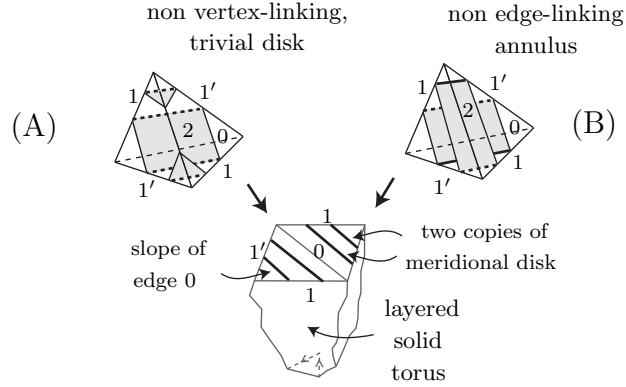}} }
\caption{Examples of a non vertex-linking, trivial, normal disk and
a non edge-linking, normal annulus.} \label{f-non-vert-edge-linking}
\end{figure}

Secondly, and in the previous situation, if we attach two banding
quads in $\td{\Delta}_{t+1}$ to $P_t$, then we get  (4B) a {\it non
edge-linking annulus} in $\T_{t+1}$. Notice that if we take only one
copy of the meridional disk in $\T_t$ and attach one banding quad in
$\T_{t+1}$, then we get a M\"obius band and the non edge-linking
annulus is the double of this M\"obius band. Also, note that the
slope of its boundary corresponds to an edge slope of an edge
meeting the meridional disk $2$ times. See Figure
\ref{f-non-vert-edge-linking} (B). Specific examples of (4A) and
(4B) can be constructed in a four tetrahedron layered-triangulation
of the solid torus extending the $1/1$--triangulation and
corresponding to a path in the $L$--graph starting at $1/1$, going
to $0/1$ back through $1/1$ to $1/2$ and ending at $1/1$. The
boundary of a non edge-linking annulus also has the slope of an
edge; it is analogous to an edge-linking annulus but is a banding of
a non vertex-linking trivial disk and so we do not think of it as an
edge-linking annulus, as we do when we band a vertex-linking disk.

\item[(5)] We have still another possibility for a new genus $0$ normal surface, if
there are two edges having the same slope. Suppose the edges $e$ and
$e'$ in $\T_t$ have the same slope and, say $e$ appears earlier in
the layering than $e'$. Then if we take an edge-linking annulus
about $e$, then at the level that $e'$ first appears, the
edge-linking annulus about $e$ has boundary slope that of the edge
$e'$. We can then continue to shove this annulus through the layered
triangulation, getting an edge-linking annulus unlike those above.
We call this (5) a {\it fat edge-linking annulus}. An example can be
obtained in the three-tetrahedron layered-triangulation of the solid
torus corresponding to an extension of the $1/2$--triangulation for
the path $1/2, 1/1, 1/2, 1/1$ in the $L$--graph. See the first three
layers in Figure \ref{f-fat-tube} and for edges take $e = 3$ and
$e'=3'$.
\end{enumerate}

We have from Theorem \ref{layered-0-1} that these examples give all
the five types of genus $0$ normal surfaces in a general
layered-triangulation of a solid torus. We now consider genus $1$
normal surfaces.

\vspace{.125 in}\noindent{\it Genus $1$ normal surfaces}: A genus
$1$ normal surface is obtained in $\T_{t+1}$ by pushing through a
genus $1$ normal surface from $\T_t$ or by adding a banding quad in
$\td{\Delta}_{t+1}$ to a genus $0$ surface in $\T_t$ and reducing
the number of boundary components (raising the genus) or adding a
banding quad in $\td{\Delta}_{t+1}$ to a genus $1$ surface in $\T_t$
and increasing the number of boundary components (keeping genus
$1$).

\begin{enumerate}
\item [(6)] We have as possible genus $1$ normal surfaces any of
those that appear in minimal layered-triangulations. These can be
(6A) a
 {\it vertex-linking disk with a thin edge-linking tube}, which
gives a once-punctured torus with trivial boundary curve, or  (6B)
an {\it edge-linking annulus with a thin edge-linking tube}. The
last example gives a twice punctured torus.

\item [(7)] Using the above possibility of a non vertex-linking
trivial disk, we obtain new genus $1$ surfaces.  The simplest of
these is a  (7A) {\it non vertex-linking trivial disk with a thin
edge-linking tube}. In this case we have a non vertex-linking
trivial disk with boundary  the trivial vertex-linking curve; hence,
at some later stage, we can add a thin edge-linking tube. On the
other hand, if we have such a non vertex-linking trivial disk with a
thin edge-linking tube at some stage in a layering, then we can band
it and get (7B) a  {\it  non edge-linking annulus with a thin
edge-linking tube}. As with edge-linking annuli, these annuli have
boundary slope the slope of an edge. We do not give a different name
for these, as we did for edge-linking annuli, depending on whether
the boundary slope is that of an edge in the boundary or not.
Examples (7A) and (7B) can happen only if an edge has the slope of
the meridian. Specific examples can be obtained for (7A) and (7B) in
a six-tetrahedron layered-triangulation of the solid torus extending
the $2/5$--triangulation, using the path $2/5, 2/3, 1/2, 1/1, 0/1,
1/1$. In this example, we have a non vertex-linking trivial disk in
the initial three layers ($1/1, 0/1, 1/1$); through the next three
layers, we can attach a thin edge-linking tube about the edge $3$.
Thus in  this six-tetrahedron layered-triangulation of the solid
torus extending the $2/5$--triangulation, we have non vertex-linking
trivial disk with a thin edge-linking tube about the edge $3$ (7A)
and also we have a non edge-linking annulus, with boundary slope
that of the edge $7$, with a thin edge-linking tube about the edge
$3$. See Figure \ref{f-non-link-and-tube}.

\begin{figure}[h]
\psfrag{0}{\small{$0$}}\psfrag{1}{\small{$1$}}\psfrag{2}{\small{$2$}}
\psfrag{3}{\small{$3$}}\psfrag{4}{\small{$4$}}
\psfrag{5}{\small{$5$}}\psfrag{6}{\small{$6$}}\psfrag{7}{\small{$7$}}
\psfrag{8}{\small{$1'$}}\psfrag{A}{\Large{(A)}}
\psfrag{B}{\Large{(B)}}

 {\epsfxsize = 3.5 in
\centerline{\epsfbox{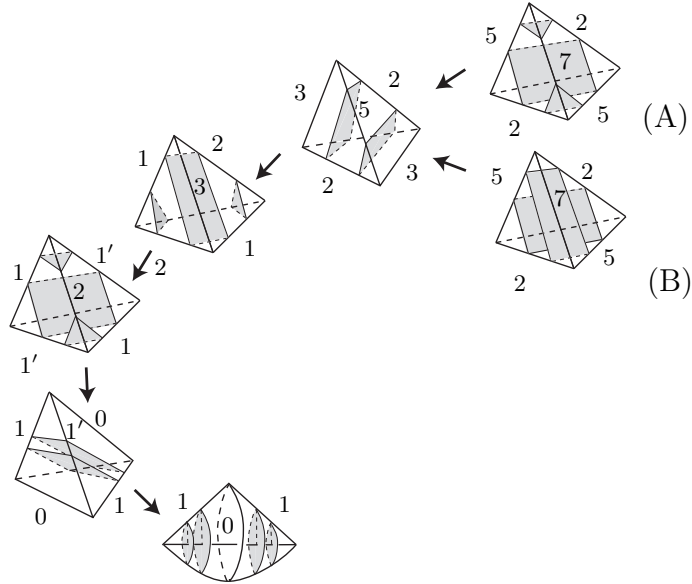}} } \caption{Examples
of a non vertex-linking trivial normal disk with a thin edge-linking
tube (A) and a non edge-linking normal annulus with a thin
edge-linking tube (B).} \label{f-non-link-and-tube}
\end{figure}

\item[(8)] We also have another possibility for new genus $1$
normal surfaces, if the edge $e$ in $\T_t$, along which we attach
$\td{\Delta}_{t+1}$, has the same slope as another edge. To see
this, suppose the edge $e$ has the same slope as the edge $e'$,
which is necessarily an interior edge of $\T_t$. Let $P_t$ be an
edge-linking annulus about $e'$; then $P_t$ is not a thin
edge-linking annulus but does have boundary slope the same as the
boundary edge $e$ in $\T_t$. Thus, we can add a banding quad in
$\td{\Delta}_{t+1}$, which necessarily runs between distinct
components of $\bdy P_t$. Attaching such a banding quad and two
triangles  in $\Delta_{t+1}$ to $P_t$, we get  (8A) a {\it
vertex-linking disk with a non edge-linking tube}. The tube links
two edges with the same slope; we think of these tubes as ``fat"
tubes. This can happen only if another edge has the same slope as
$e$. The resulting surface is a punctured-torus.  If we attach two
banding quads, then we get (8B) an {\it edge-linking annulus with a
non edge-linking tube}; and the surface is a twice-punctured torus.
Also, we can have an edge-linking annulus with a non edge-linking
tube and have still another edge with the same slope as that of the
boundary of the annulus; thus we get a fat edge-linking annulus with
a non edge-linking tube. Specific examples of both (8A) and (8B) can
be obtained in a four-tetrahedron layered-triangulation of the solid
torus corresponding to an extension of the $1/1$--triangulation for
the path $1/1, 1/2, 1/1, 1/2, 1/1$ in the $L$--graph. See Figure
\ref{f-fat-tube}. Having edges with the same slope leads to two
other possibilities. For example, before we see the two edges with
the same slope, we can have a vertex-linking disk with a thin
edge-linking tube; then we have two edges with the same slope. We
then have (8C) a {\it fat annulus with a thin edge-linking tube}.
Also, we could have a vertex-linking disk with a non edge-linking
tube and then see two more edges with the same slope, getting (8D) a
{\it fat annulus with a non edge-linking tube}. Examples can be
obtained from the example in Figure \ref{f-fat-tube} by adding a
couple more layers.

\begin{figure}[h]
\psfrag{0}{\small{$0$}}\psfrag{1}{\small{$1$}}\psfrag{2}{\small{$2$}}
\psfrag{3}{\small{$3$}}\psfrag{4}{\small{$1'$}}
\psfrag{5}{\small{$3'$}}\psfrag{6}{\small{$1''$}}\psfrag{A}{\Large{(A)}}
\psfrag{B}{\Large{(B)}}

 {\epsfxsize = 2.5 in
\centerline{\epsfbox{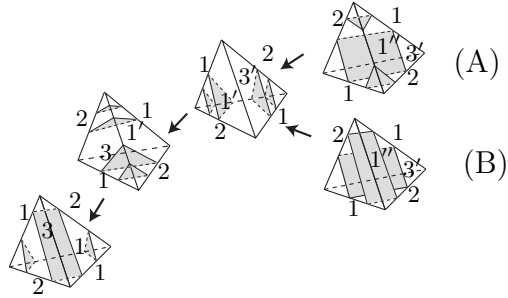}} } \caption{Examples of  (A) a
vertex-linking disk with a non edge-linking (``fat") tube and (B) an
edge-linking annulus with a non edge-linking (``fat") tube .}
\label{f-fat-tube}
\end{figure}

\item[(9)] Here we have both an edge with slope that of a meridian
disk as in Example $(7)$ and two edges with the same slope, as in
Example $(8)$. Hence, we have (9A) a {\it non vertex-linking disk
with a non edge-linking tube} or (9B) a {\it non edge-linking
annulus with a non edge-linking tube}. An example can be easily
constructed by using the initial segment with four tetrahedron
extending the $1/2$--triangulation in Figure
\ref{f-non-link-and-tube}, from this we can get an example of a non
vertex-linking trivial disk; then we layer on these four tetrahedra
the last three tetrahedra given in Figure \ref{f-fat-tube} ending
with a seven tetrahedra layering. Using the two options (A) and (B)
in Figure \ref{f-fat-tube}, we get examples of (9A) and (9B),
respectively.
\end{enumerate}

For the classification of genus $0$ and genus $1$ normal surfaces,
it is possible that we have a creased $3$--cell at the base of our
layered-triangulation and it also is possible that we layer along
the univalent edge in the boundary. In fact, it could be possible to
have a path in the $L$--graph beginning at some point, going through
$1/1$ to $0/1$ and then simply oscillating back and forth between
$0/1$ and $1/1$, finally ending at $1/1$. Such a path gives a
terribly degenerate triangulation; in fact, one in which some
initial segment does not even give a $3$--manifold. Such an initial
segment can be crushed to the M\"obius band. We will assume that the
only initial segment that is not a $3$--manifold is when we have a
one-tetrahedron creased $3$--cell. In Lemma \ref{one-tet-torus}, we
listed all normal surfaces in the one-tetrahedron solid torus; and
in Lemma \ref{0-1-1}, we listed the connected normal surfaces in the
two-tetrahedron minimal triangulation extending the
$0/1$--triangulation. In these cases, there are only genus $0$ and
genus $1$ normal surfaces. We consider the three-tetrahedron
layered-triangulation extending the $0/1$-triangulation on the
boundary as part of the general induction step. With this
information, we can begin our induction argument, building upon the
surfaces in these small examples.

\begin{thm}
\label{layered-0-1}  A connected, embedded genus $0$ or genus $1$
normal surface in a (general) layered-triangulation of a solid torus
is normally isotopic to one of:

\noindent Genus $0$:
\begin{enumerate} \item Vertex-linking disk,\item  meridional
disk (unique),\item an edge-linking annulus with boundary slope the
slope of an edge, \item (A) a trivial, non vertex-linking disk or
(B) a non edge-linking annulus with boundary slope the slope of an
edge (only if an edge has the slope of the meridional disk) or \item
a fat edge-linking annulus (only if two edges have the same slope).
\end{enumerate}

\noindent Genus $1$:
\begin{enumerate}
\item[6] (A) a vertex-linking disk with a thin edge-linking tube
or (B) an edge-linking annulus with a thin edge-linking tube,
\item [7] (A) a non vertex-linking trivial disk with a thin
edge-linking tube or (B) a non edge-linking annulus with a thin
edge-linking tube  (only if an edge has the slope of the meridian),
\item[8] (A) a vertex-linking disk with a non edge-linking tube or
(B) an edge-linking annulus with a non edge-linking tube or (C) a
fat edge-linking annulus with a thin edge-linking tube, or (D) a fat
edge-linking annulus with a non edge-linking tube (in all cases only
two edges have the same slope),
\item[9] (A) a non vertex-linking trivial disk with a non
edge-linking tube or (B) a non edge-linking annulus with a non
edge-linking tube  (only if an edge has the slope of the meridian
and two edges have the same slope).
\end{enumerate}
\end{thm}

The proof uses induction, following the same ideas as in the proof
of Theorem \ref{normal}. We change the topological type of a surface
by attaching a banding quad; however, a banding quad attached to a
surface having boundary the vertex-linking curve gives an annulus or
an annulus with tubes and a banding quad attached to an annulus (or
an annulus with tubes) must run between distinct boundary components
and so raises genus. See \cite{jac-sedg-dehn} for a detailed proof.

If we ever have the meridional slope also being the slope of an {\it
interior} edge, then, except for the two-tetrahedra extension of the
$0/1$--triangulation, we do {\it not} have a minimal
layered-triangulation. On the other hand, we may construct a
layered-triangulation that is not a minimal layered-triangulation,
and does not have the meridional slope as the slope of any edge.

As mentioned above, if we do not have a minimal
layered-triangulation, there is a vertex on the path corresponding
to the layered-triangulation in the $L$--graph where we have layered
along the univalent edge. If we layer along the univalent edge, we
either get two thin edges with the same slope or we pick up the
exceptional longitude (defined below) as a slope of an edge. For
example, in the two-tetrahedron layered-triangulation extending the
$1/1$--triangulation and the three-tetrahedron one extending the
$0/1$--triangulation, we do not have two thin edges with the same
slope; however, we do have the exceptional longitude having the
slope of the edge labeled $1'$ in Figure \ref{f-small-layered}. If
we layer along the univalent edge and do have two edges with the
same slope, then these slopes will necessarily meet the meridional
disk the same number of times. We try to distinguish the edges by
using the same integer but adding primes as superscripts. Just
because edges meet the meridional disk the same number of times does
not mean they determine the same slope on the boundary. For example,
at the vertex labeled $0/1$ in the $L$--graph, we have an edge
(cycle) with both its end points at $0/1$. If we take a
layered-triangulation corresponding to a path in the $L$--graph that
transverses this cycle a number of times, we get a number of edges
meeting the meridional disk just once but they determine different
slopes (longitudes) on the boundary. We summarize the conclusions we
need later to classify $0$-- and $1$--efficient
layered-triangulations of lens spaces in the following corollaries.
Most notably is that the number of boundary components of a genus
$0$ or a genus $1$ normal surface is either one or two and we have a
normal disk or normal annulus possibly with a single tube.

\begin{cor} A connected, embedded, normal, genus $0$ surface in a
layered-triangulation of a solid torus is either the unique meridian
disk, a trivial disk with boundary the vertex-linking curve, or an
annulus with boundary slope the slope of an edge.

\end{cor}

\begin{cor} A connected, embedded, normal, genus $1$ surface in a
layered-triangulation of a solid torus is either a disk with a tube
(a punctured-torus) with boundary the vertex-linking curve or an
annulus with a tube (a twice-punctured torus) with boundary slope
the slope of an edge.
\end{cor}

\subsection{Almost normal surfaces in minimal layered-triangulations
of a solid torus}

\vspace{.1 in} In this section we classify the almost normal
surfaces in a minimal layered-triangulation of a solid torus that
extends a given one-vertex triangulation on its boundary.

 We begin by giving some examples of almost
normal octagonal surfaces in minimal layered-triangulations of a
solid torus. There are, of course, various almost normal tubed
surfaces formed from normal surfaces in these triangulations.

\vspace{.125 in}\noindent{\bf Examples.}
\begin{enumerate}
\item Three different almost normal octagonal annuli, two with
boundary slope a ``longitude" in the minimal
$2/3$--layered-triangulation of the solid torus and one with
boundary slope the slope of the thin edge $2$ in the minimal
$1/4$--layered-triangulation. See Figure \ref{f-longitude-an}.
Notice that (A) has slope the preferred longitude, which is a slope
of an edge, while (B) has slope also a longitude but not what we are
calling the slope of an edge. The third (C) has slope the same as a
thin edge and is isotopic to an edge-linking normal annulus about
the edge $2$. Examples (A) and (B) are not isotopic to edge-linking
normal annuli. Example (B) is a unique exception; all other
octagonal annuli have boundary slope the slope of an edge. We will
refer to such an octagonal annulus in a layered solid torus as an
\emph{exceptional} octagonal annulus and its boundary is called an
\emph{exceptional longitude}.

\begin{figure}[h]
\psfrag{1}{\small{$1$}}\psfrag{2}{\small{$2$}}\psfrag{3}{\small{$3$}}\psfrag{4}{\small{$4$}}
\psfrag{5}{\small{$5$}}\psfrag{A}{\Large{(A)}}\psfrag{B}{\Large{(B)}}
\psfrag{C}{\Large{(C)}}\psfrag{a}{\begin{tabular}{c}preferred\\longitude\end{tabular}}
\psfrag{b}{\begin{tabular}{c}exceptional\\longitude\end{tabular}}\psfrag{c}{\begin{tabular}{c}slope
of\\edge
$2$\end{tabular}}\psfrag{y}{\begin{tabular}{c}\LARGE{$\approx$}\\isotopy\end{tabular}}
 {\epsfxsize = 3.5 in
\centerline{\epsfbox{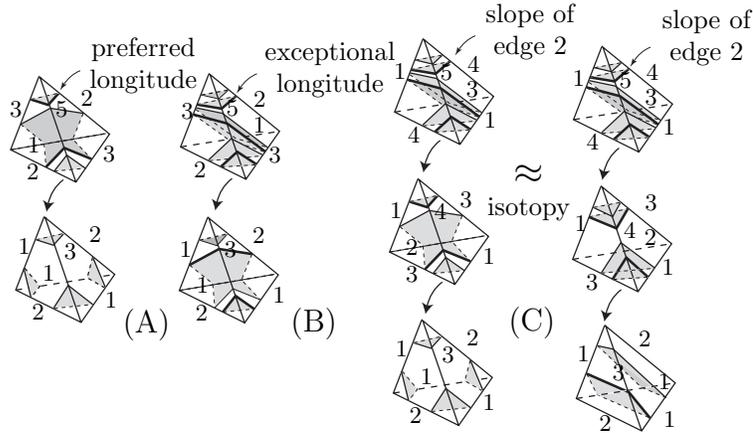}} } \caption{Examples of
(A) an almost normal octagonal annulus with slope the preferred
longitude; (B) an almost normal octagonal annulus with slope the
exceptional longitude; and an almost normal octagonal annulus with
slope that of the edge $2$.} \label{f-longitude-an}
\end{figure}

\item An
almost normal, octagonal annulus with a thin edge-linking tube in
the minimal $1/4$--layered-triangulation. See Figure
\ref{f-an-annulus-tube}. Often it will be the case that an almost
normal annulus with tubes is isotopic to a normal annulus with
tubes; however, this example is not. The slope of the boundary of
the almost normal annulus is that of the edge $3$. The tubes must
come before the octagon in the layering.

\begin{figure}[h]
\psfrag{1}{\small{$1$}}\psfrag{2}{\small{$2$}}\psfrag{3}{\small{$3$}}\psfrag{4}{\small{$4$}}
\psfrag{5}{\small{$5$}}\psfrag{a}{\begin{tabular}{c}slope of\\edge
$3$\end{tabular}}
 {\epsfxsize = 2.5 in
\centerline{\epsfbox{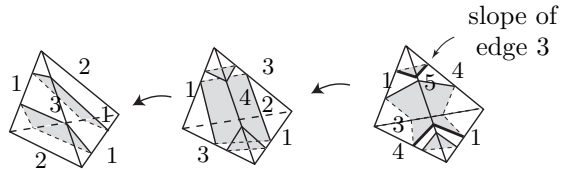}} } \caption{Example of
an almost normal octagonal annulus with a tube and boundary slope
 that of the edge $3$.} \label{f-an-annulus-tube}
\end{figure}

\item Two (interesting) almost normal octagonal disks with trivial
boundary slope: one in the three-tetrahedra layered-triangulation
extending the $0/1$--layered-triangulation (Figure
\ref{f-an-disk}(A)); and the other in the minimal two-tetrahedra
triangulation of the solid torus extending the $0/1$--triangulation
(Figure \ref{f-an-disk}(B)). The first is formed by using an almost
normal octagon to band two copies of the meridional disk together.
The second is a ``push through" of an almost normal octagonal disk
embedded in the creased $3$--cell. Notice that the latter almost
normal octagon may be viewed as the solution $y_2 = y_3$ in Lemma
\ref{0-1-1}. In general, almost normal octagons can be found as
solutions to the normal equations (or quad equations) having two
distinct quads in the same tetrahedron. Both normal octagonal disks
isotopically shrink in one direction to the vertex-linking disk.
These octagonal disks lead to almost normal octagonal $2$--spheres
in examples of $0$--efficient triangulations of $S^3$.

\begin{figure}[h]
\psfrag{1}{\small{$1$}}\psfrag{2}{\small{$2$}}\psfrag{3}{\small{$3$}}
\psfrag{0}{\small{$0$}}
\psfrag{4}{\small{$1'$}}\psfrag{A}{\Large{(A)}}\psfrag{B}{\Large{(B)}}
 {\epsfxsize = 2.5 in
\centerline{\epsfbox{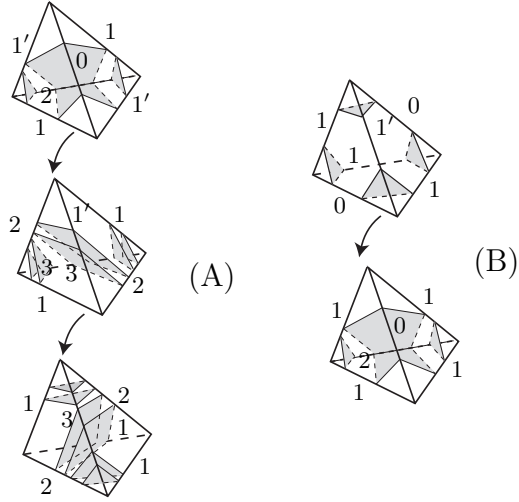}} } \caption{Examples of almost
normal octagonal disks with trivial boundary slope.}
\label{f-an-disk}
\end{figure}

\end{enumerate}

  We now
give the classification of almost normal surfaces in a minimal
layered-triangulation of a solid torus. In many of these cases the
almost normal surface is isotopic to a normal surface; we note this
feature in such examples. Neither of the  longitudinal, octagonal
annuli are isotopic to a normal annulus; whereas,  other octagonal
annuli (typically) are isotopic to a normal edge-linking one. We
also observe that an almost normal tube along a thin edge is
(typically) isotopic over the edge to a thin edge-linking tube. The
minimal (two-tetrahedron) triangulation and the three-tetrahedron
layered-triangulation of the solid torus extending the
$0/1$--triangulation both admit an almost normal octagonal disk,
with boundary the vertex-linking curve, which is isotopic to the
vertex-linking disk. These do not fit into the general pattern;
however, layered-triangulations extending $0/1$ lead to ``layered"
triangulations of the $3$--sphere and these octagonal disks, which
appear, become part of an almost normal octagonal $2$-sphere, which
is necessary in certain triangulations of the $3$--sphere. Finally,
there are a number of interesting non orientable, almost normal
octagonal surfaces in layered-triangulations of the solid torus; we
avoid these.

\begin{thm}\label{almostnormal} A connected, embedded, orientable, almost normal
surface in a minimal layered-triangulation of the solid torus is
normally isotopic to one of the following surfaces.
\begin{enumerate}
\item If octagonal:
\begin{itemize}
\item An almost normal, octagonal annulus with boundary slope the
exceptional longitude, \item An almost normal, octagonal annulus
with boundary slope a slope of an edge, possibly the preferred
longitude, \item An almost normal, octagonal annulus with thin
edge-linking tubes and having boundary slope the slope of an edge
(possibly the longitude), or \item An almost normal, octagonal disk
with boundary the vertex-linking curve, (only in the minimal
layered-triangulation extending $0/1$).

\end{itemize}
\item If tubed:
\begin{itemize}
\item The vertex-linking disk (possibly) with thin edge-linking
tubes and an almost normal tube, \item An edge-linking annulus
(possibly) with thin-edge-linking tubes and an almost normal tube,
\item Two copies of the meridional disk with an almost normal tube
from one to the other (not in the product region between them),
 \item A vertex-linking disk (possibly) with thin edge-linking
tubes and an almost normal tube between it and a distinct and
disjoint normal surface, or \item An edge-linking annulus
(possibly) with thin edge-linking tubes and an almost normal tube
between it and a distinct and disjoint normal surface.

\end{itemize}
\end{enumerate}
\end{thm}

We note, in general,  both normal annuli and almost normal octagonal
annuli have boundary slope the slope of an edge, with the single
exception the unique almost normal octagonal annulus with boundary
slope the exceptional longitude; however, all such annuli, including
the exceptional longitudinal one, are parallel into an annulus in
the boundary containing the vertex. In particular, we have {\it
companion} and {\it complementary} annuli defined for almost normal
octagonal annuli.

\begin{proof} We defer the consideration of the three-tetrahedron layered-triangulation
of the solid torus extending $0/1$ and the minimal triangulations of
the solid torus extending $1/1$ and $0/1$ until after we have
handled the general case.

We first consider the case where the almost normal surface is
octagonal. There are three possible octagons in a tetrahedron. The
octagons are distinguished by the opposite edges in the tetrahedron
that they meet twice (there are three such pair of edges). See
Figure \ref{f-layer-octagon}. This distinction is particularly
significant in a minimal layered-triangulation.

\begin{figure}[h]
\psfrag{a}{\small{$p+q$}}\psfrag{b}{\small{$q$}}\psfrag{c}{\small{$p$}}
\psfrag{1}{\Large{$(B_1)$}}\psfrag{2}{\Large{$(B_2)$}}\psfrag{3}{\Large{$(B_3)$}}
\psfrag{d}{\Large{$\td{\Delta}_{t+1}$}} {\epsfxsize = 3.5 in
\centerline{\epsfbox{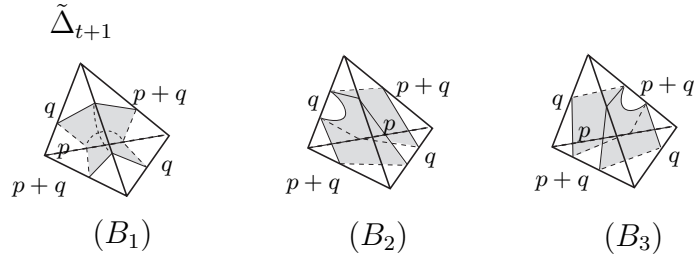}} } \caption{Examples of
octagons in the $(t+1)$-layer.} \label{f-layer-octagon}
\end{figure}

In the one-tetrahedron solid torus, we have an octagonal annulus
when the octagon meets the univalent edge (and its opposite edge)
twice. This is the exceptional octagonal annulus having boundary
slope a longitude. The other two possibilities give octagonal,
once-punctured Klein bottles.

So, let us assume our statement is true for a minimal
layered-triangulation of a solid torus having $t$ layers and let
$\T_{t+1} = \T_t \cup_{e} \td{\Delta}_{t+1}$ be a minimal
layered-triangulation having $t +1$ tetrahedra. Furthermore, suppose
the minimal layered-triangulation $\T_t$ extends the
$p/q$--triangulation  on the boundary of the solid torus.

As in the proof of Theorem \ref{normal}, we are only considering
minimal layered-triangula-tions; so, we have that the attaching edge
$e$ is either the edge labeled $p$ or the edge labeled $q$ and
$\T_{t+1}$ extends either the $(q/p+q)$--triangulation  or
$(p/p+q)$--triangulation on the boundary of the solid torus.

If $F_{t+1}$ is an  almost normal octagonal surface in $\T_{t+1}$,
then it meets $\T_t$ either in a normal surface $F_t$, each
component of which satisfies the conclusion of Theorem \ref{normal},
or in an almost normal surface $F_t$ each component of which
satisfies the conclusion of this theorem, except that none of the
components of $F_t$ can be a surface with an almost normal tube.
Thus there are two possibilities for $F_{t+1}$:

\vspace{.15 in} \noindent {\bf A.} $F_{t+1}$ is obtained from $F_t$,
which is an almost normal octagonal surface, by adding normal
triangles and quads in $\td{\Delta}_{t+1}$.

 \vspace{.15 in}
\noindent {\bf B.} $F_{t+1}$ is obtained from $F_t$, which is
normal, by adding an octagon along with possibly some triangles in
$\td{\Delta}_{t+1}$.

In case A,  we may assume that we have a component of $F_t$ that is
an  almost normal octagonal surface and therefore by induction is an
 almost normal octagonal annulus, possibly with thin edge-linking
tubes, and with boundary the slope of an edge; or $F_t$ is the
exceptional, longitudinal  almost normal octagonal annulus. Recall
the slope of the boundary of the exceptional, longitudinal almost
normal octagonal annulus behaves as a slope of an edge but is never
the slope of an edge (unless it is the exceptional octagonal annulus
in the two-tetrahedron $1/1$--layered-triangulation, whose boundary
has slope that of the univalent edge, as we shall see below). So, it
follows that we can only ``push through" the boundary of $F_t$.
Hence, $F_t$ is connected and $F_{t+1}$ is homeomorphic to $F_t$ and
its boundary is either the slope of an edge (possibly the preferred
longitude) or $F_{t+1}$ is the exceptional octagonal annulus.

In case B, we must consider adding an octagon, along with possibly
some triangles, in $\td{\Delta}_{t+1}$ to $F_t$. The situation is
symmetric by attaching along the edge $p$ or $q$ in the boundary of
$\T_t$; so, say we are attaching along the edge that is labeled $p$.

There are the three possibilities for an octagon in
$\td{\Delta}_{t+1}$, shown in Figure \ref{f-layer-octagon}; we
designate the cases as  $B_1, B_2$ and $B_3$.

In case $B_1$, the octagon meets the univalent edge in the boundary
of $\T_{t+1}$ twice, and therefore its opposite edge, which is the
attaching edge $p$, twice. We analyze the possibilities by
considering how we might add triangles, along with the octagon, in
$\td{\Delta}_{t+1}$. If there are no triangles, then the boundary of
$F_t$ would meet the edge $p$ twice, which would be a maximal number
of times; however, this is possible only if $T_t$ is the
two-tetrahedron minimal layered-triangulation of the solid torus
extending the $1/1$--triangulation, which it is not. If we add
triangles meeting the edge $p$, then the edge $p$ would have to be
the univalent edge in $\T_t$, which it is not since we have a
minimal layered-triangulation distinct from those extending $0/1$
and $1/1$. If we add triangles not meeting the edge labeled $p$,
then we have a vertex-linking curve as a curve in the boundary of
$F_t$. Since we have $F_{t+1}$ connected, then the only possibility
is adding two triangles, $F_t$ is the vertex-linking disk with
possibly some thin edge-linking tubes, making $F_{t+1}$ an almost
normal octagonal annulus with possibly some thin edge-linking tubes
and with boundary slope the slope of the edge $p$.

In case $B_2$, the octagon meets the edge $q$, twice. If we do not
add triangles in $\td{\Delta}_{t+1}$, then, as above, $F_t$ would
meet the edge $q$ a maximal number of times, which is not possible
in our situation. If we add triangles, then any added triangle meets
the edge $q$; again this would force $F_t$ to meet the edge $q$ a
maximal number of times.  Hence, this situation is impossible.

In case $B_3$, the octagon meets the univalent edge of $\T_t$, the
edge $p+q$, twice. In this case, if we do not add any triangles,
then the single octagon component in $\td{\Delta}_{t+1}$ meets the
boundary of $\T_t$ in a slope that is possibly the boundary of a
normal surface. However, the normal surface must have {\it
connected} and nontrivial boundary and meet the univalent edge
precisely two times. It follows that the only possibility for $F_t$
would be a nonorientable surface. We have excluded such a
possibility.  If we add triangles at every corner, then again we
would have neither $F_t$ nor $F_{t+1}$ connected. If we add
triangles that meet the edge labeled $p$, then the surface $F_t$
would have a vertex-linking curve and so could not be connected; but
then the octagonal surface would only band on the vertex-linking
disk of $F_t$ and $F_{t+1}$ would not be connected. Finally, the
only possibility is to add triangles that do not meet the edge
labeled $p$. In this case we can get a permissible slope for a
surface $F_t$; however, $F_t$ must have connected boundary and so
must be nonorientable. Notice that in this case, adding the octagon
in $\td{\Delta}_{t+1}$ increases the number of cross caps in
$F_{t+1}$ to two more than in $F_t$.

Our conclusions are that if $F_{t+1}$ is orientable and octagonal,
then in Case (A) the surface $F_t$ is itself an  almost normal
octagonal surface and we ``push through" $F_t$ to get $F_{t+1}$; in
Case ($B_1$) the surface $F_t$ is normal and we add an octagon and
two triangles in $\td{\Delta}_{t+1}$ to $F_t$, which makes $F_{t+1}$
an octagonal annulus (possibly) with thin edge-linking tubes. The
almost normal octagonal annulus $F_{t+1}$ can be moved by an isotopy
to a normal annulus, unless the attaching edge for
$\td{\Delta}_{t+1}$ is the preferred longitude or there is a thin
edge-linking tube in $F_t$ that meets the attaching edge.

Now, we consider the situation where the almost normal surface is a
tubed surface and is obtained either from a connected normal surface
by adding a tube along an edge (Case ($C_1$)) or is obtained from
two disjoint connected normal surfaces (not necessarily distinct) by
adding a tube from one surface to the other (Case ($C_2$)). Since we
have classified all normal surfaces in a minimal
layered-triangulation, we only need consider the possibilities for
adding an almost normal tube to these normal surfaces.

\vspace{.125 in} $\mathbf{C_1.}$ Adding an almost normal tube to a
connected normal surface.

We make some observations. If the normal surface is the
vertex-linking disk and we add an almost normal tube along any edge
except the unique ``thick edge" or an edge in the boundary, then the
surface is isotopic to a normal surface with a thin edge-linking
tube; the edge-linking tube is about the edge along which we added
the almost normal tube.  If the normal surface is a vertex-linking
disk with thin edge-linking tubes, we can not add an almost normal
tube along an edge that has a thin edge-linking tube. If we add an
almost normal tube along an edge other than the unique ``thick edge"
or an edge in the boundary, then we have the possibility of an
isotopy of the almost normal tube over the edge to get a normal
surface with a thin edge-linking tube about that edge. However,
there is an obstruction; if there is a thin edge-linking tube in the
surface having a quadrilateral that meets the edge adjacent to the
almost normal tube. If this is the case, then we can not just shove
the almost normal tube over the edge to get a thin edge-linking tube
without destroying the existing thin edge-linking tube. If the
normal surface is the meridional disk and we add an almost normal
tube, we get a non orientable surface. If the normal surface is an
edge-linking annulus, having boundary slope the slope of an edge. We
can not add a tube along the edge the annulus is linking. However,
we may add the tube along any edge in the initial segment before we
come to the edge the annulus is linking. We can move such an almost
normal tube over the edge by an isotopy to get a thin edge-linking
tube or the almost normal tube is adjacent to the unique ``thick
edge" or an edge that meets a quad in the edge-linking annulus. So,
the possibilities can be summarized as an edge-linking annulus
having boundary slope the slope of an edge and an almost normal tube
along an edge. We also can have the almost normal tube added after
the edge that the annulus is linking. In this case, the tube can be
moved by an isotopy to be along an edge in the boundary of the solid
torus (the argument is inductive and is the same as that below). If
the normal surface is an edge-linking annulus with tubes, having
boundary slope the slope of an edge. As in the preceding, we can not
add a tube along the edge the annulus is linking or along any of the
edges having thin edge-linking tubes. The almost normal tube may be
added before the edge of the edge-linking annulus. If the edge along
which it is added does not meet a quad of a thin edge-linking tube
or a quad in the edge-linking annulus, then we can move the almost
normal tube over the edge by an isotopy to get a thin edge-linking
tube. Again, it is possible the almost normal tube is added after
the edge that the annulus links. We then have an edge-linking
annulus with thin edge-linking tubes and having boundary slope the
slope of an edge and an almost normal tube along an edge in the
boundary of the solid torus.

\vspace{.25 in} $\mathbf{C_2.}$ Adding an almost normal tube between
two distinct normal surfaces.

If two distinct normal surfaces are connected via an almost normal
tube, then up to isotopy, we show that we may assume the almost
normal tube is along an edge in the boundary of the solid torus.
To do this, we use induction.

The induction begins easily enough for the one-tetrahedron solid
torus - all edges are in the boundary.

So, we assume the statement is true for a layered-triangulation of
length $t$ and consider the minimal layered-triangulation,
$\T_{t+1}$, having $t+1$ tetrahedra. If we have two connected normal
surfaces, $F_{t+1}$ and $F_{t+1}'$ in $\T_{t+1}$, then each must
meet $\T_t$ in a connected normal surface, say $F_t$ and $F_{t}'$,
respectively. If an almost normal tube is added to connect $F_t$ and
$F_{t}'$, then, by the induction hypothesis, it is along the
boundary. This tube will be along the boundary of $\T_{t+1}$ unless
it is along the edge of the layering. However, a tube along this
edge must be either between parallel normal triangles or parallel
normal quadrilaterals or between a triangle and a quadrilateral. In
any case, it is isotopic to an almost normal tube along an edge in
the boundary of the solid torus with the triangulation $\T_{t+1}$.
Thus we see that by an isotopy we can have the almost normal tube
along an edge in the boundary.

 Finally, we consider the almost normal
surfaces in the special cases for the minimal extensions of the
triangulations $1/1$ and $0/1$.

\vspace{.125 in} $\mathbf{D}.$ Almost normal surfaces in the
two-tetrahedra $1/1$--layered-triangulation.

In the case of the two-tetrahedron extension of $1/1$, we can use
quadrilateral equations to see two octagonal solutions: one with an
octagon in the one-tetrahedron solid torus, which meets the
univalent edge twice, and another with an octagon in the second
tetrahedron in the layering, also having two points on the univalent
edge. Thus we have just the two possibilities shown in the {\it
first two layers} in Figure \ref{f-0-1-1}($D_1)$ and $(D_2)$:
$(D_1)$ is an octagonal annulus with boundary slope the exceptional
longitude, which happens in this case to be the slope of the thin
edge labeled $1'$,  and $(D_2)$ is an octagonal annulus with
boundary slope the slope of the edge labeled $3$. The latter can be
moved via an isotopy to an edge-linking annulus about the edge
labeled $3$.


The almost normal tubed surfaces in the two-tetrahedron extension of
$1/1$ can be easily listed from the normal surfaces. There are no
surprises.

\vspace{.125 in} {\bf E.} Almost normal surfaces in the
 small layered-triangulations extending the $0/1$--triangulation.

We consider both the three-tetrahedron layered-triangulation  and
the minimal two-tetrahedron triangulation of the solid torus
extending $0/1$.

For the three-tetrahedron triangulation, we have the layering of a
tetrahedron onto the two-tetrahedron layered-triangulation extending
$1/1$; the almost normal surfaces of which were analyzed in the
preceding paragraphs.

First, there are the almost normal octagonal surfaces in minimal
$1/1$--layered-triangulation that can be pushed through. See Figure
\ref{f-0-1-1}$D_1$ and $D_2$. These give an almost normal octagonal
annulus with boundary slope the exceptional longitude, labeled $1'$,
and an almost normal annulus with boundary slope the slope of the
edge $3$. So, we need only consider the possibilities if we add an
octagonal disk in the third layer. In this situation, we have the
same possibilities as in the induction argument above. The first we
consider is when the octagonal disk in the third layer meets the
univalent edge, labeled $0$, and the attaching edge, labeled $2$,
twice. Just as above, in this situation, we can attach the octagonal
disk to the vertex-linking curve. The possibilities are that we can
attach the octagonal disk to either the vertex-linking disk, getting
an almost normal octagonal annulus with boundary slope the slope of
the edge $2$ (see Figure \ref{f-0-1-1}$E_1$), or to the
vertex-linking disk with a thin edge-linking tube about the edge
$3$, getting an almost normal octagonal annulus with a thin
edge-linking tube and, again, having boundary slope the slope of the
edge $2$ (see Figure \ref{f-0-1-1}$E_2$).

\begin{figure}[h]
\psfrag{1}{\footnotesize{$1$}}\psfrag{2}{\footnotesize{$2$}}
\psfrag{3}{\footnotesize{$3$}}
\psfrag{0}{\footnotesize{$0$}}\psfrag{4}{\footnotesize{$1'$}}
\psfrag{c}{\Large{$(D_1)$}}\psfrag{C}{\Large{$(D_2)$}}\psfrag{e}{\Large{$(E_1)$}}
\psfrag{E}{\Large{$(E_2)$}}
\psfrag{F}{\Large{$(F)$}}{\epsfxsize = 3.5 in
\centerline{\epsfbox{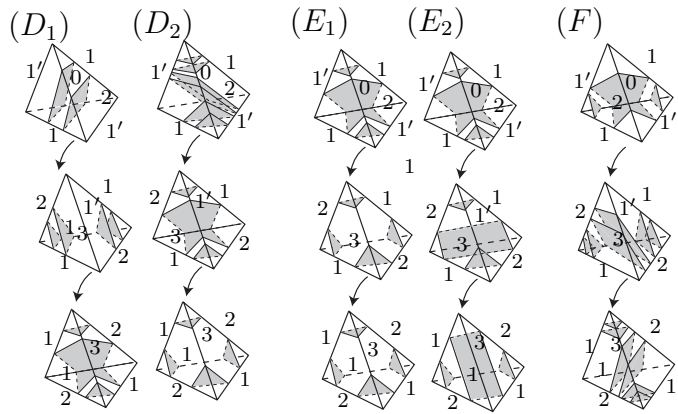}} } \caption{Almost normal
octagonal surfaces in the three-tetrahedron, layered extension of
the $0/1$--triangulation. In $D_1$ and $D_2$ are almost normal
octagonal surfaces in the two-tetrahedron extension of the
$1/1$--triangulation, which can be seen in the initial segments. }
\label{f-0-1-1}
\end{figure}

Continuing, we have a possibility here that we did not have in the
general induction argument above. The univalent edge in the
$1/1$--layered-triangulation is {\it not} the edge meeting the
meridional disk a maximal number of times; thus we get a {\it new}
almost normal surface by attaching the octagonal disk and two
triangles meeting the edge $2$ to two copies of the meridional disk
in the $1/1$--layered-triangulation. See Figure \ref{f-0-1-1}$F$.
This gives us an almost normal, trivial, non vertex-linking disk
having boundary the vertex-linking curve. It is a banding (by the
octagon) of two copies of the meridional disk. We can also band just
one copy of the meridional disk by not adding any triangles in the
third layer; however, this gives an almost normal M\"obius band. The
other two possibilities for adding an octagonal disk in the third
layer can be analyzed as above in the induction step. There are no
orientable surfaces possible in these two cases.

For the minimal two-tetrahedron triangulation extending $0/1$, we go
back to the quadrilateral solutions in the proof of Theorem
\ref{0-1-1} and Figure \ref{f-normal-0-1-1}. There is only one
quadrilateral equation (from the interior edge $2$). It is $y_2 =
y_3$. Otherwise, we have that $y_4, y_5$ and $y_6$ can take on
independent values.

For $y_2 = y_3 =1$ and $y_i = 0$ otherwise, we have a solution with
triangular coordinates, $x_2 = x_4 = x_6 = x_5 = x_7 = 1$, which is
an  almost normal, trivial, octagonal disk with boundary the
vertex-linking curve. For $y_2 = y_3 =1, y_4 = 1$ and $y_1 = y_5 =
y_6 = 0$, we have a solution with triangular coordinates, $x_2 = x_4
= x_6 = x_8 = 1$, which is a  longitudinal, almost normal octagonal
annulus having boundary slope that of the univalent edge labeled
$1'$ in the boundary, a longitude. Note that here this is not an
exceptional longitude and is the slope of an edge. For $y_2 = y_3
=1, y_5 = 1$ and $y_1 = y_4 = y_6 = 0$, we have a non orientable
solution (an almost normal M\"obius band with boundary slope that of
the edge $2$). For $y_2 = y_3 =1, y_6 = 1$ and $y_1 = y_4 = y_5 =
0$, we have a solution with triangular coordinates $x_2 = x_4 = 2$
and $x_5 = x_7 = x_6 = x_8$, which gives us disjoint copies of the
meridional disk and the  almost normal, trivial, octagonal disk with
boundary the vertex-linking curve. Notice that two such surfaces
also are in the three-tetrahedron minimal layered-triangulation of
the solid torus extending $0/1$ and are similarly disjoint there.
This takes care of all possibilities for an octagonal disk in the
first tetrahedron.

For $y_4 = y_5 = 1, y_i = 0$, otherwise, and $y_4 = y_6 = 1, y_i =
0$, otherwise, we have nonorientable solutions. For $y_5 = y_6 = 1,
y_i = 0$, otherwise, we have a solution with triangular coordinates
$x_1 = x_3 = x_2 = x_4 = x_5 = x_7$, which is an  almost normal,
exceptional, longitudinal octagonal annulus.

Thus in the two-tetrahedron, minimal triangulation of the solid
torus extending $0/1$ we have an almost normal  trivial octagonal
disk, an almost normal octagonal annulus with boundary slope the
exceptional longitude, and an almost normal octagonal annulus with
boundary slope that of the new longitude that is an edge in the
boundary (not the preferred longitude, also the slope of an edge in
the boundary). The almost normal tubed surfaces can be listed off by
using the possibilities for normal surfaces from Theorem
\ref{0-1-1}.\end{proof}

\subsection{Nearly-minimal, layered-triangulations of a solid
torus} As is evident, one of the key features of
layered-triangulations is the ability to classify the normal and
almost normal surfaces in such a triangulation. If we have two thin
edges with the same slope or have an edge with the slope of the
meridian, we get more complex normal surfaces than in the case of a
minimal layered-triangulation. It follows that
layered-triangulations of the solid torus, in which the meridional
slope is not a slope of an edge and no two {\it thin} edges have the
same slope, have the same structure for normal surfaces as do
minimal layered-triangulations. We say a layered-triangulation of a
solid torus is {\it nearly-minimal} if the meridional slope is not a
slope of an edge and no two thin edges have the same slope. For
layered-triangulations of the solid torus a minimal layered
triangulation is nearly-minimal, except for the two-tetrahedra
layered-triangulation extending the $0/1$--triangulation, where the
edge labeled $0$ has slope that of the meridian.

\vspace{.125in}\noindent{\bf Examples.}
\begin{enumerate}
\item The opening along $1,1,2$ of a minimal layered-triangulation is nearly-minimal.

\item A nearly-minimal but non minimal layered-triangulation
extending $1/6$. Start with the minimal layered-triangulation
extending $1/6$. At the point $1/3$ we add the path to $3/4$ and
back to $1/3$.  This increases the layered-triangulation by two
tetrahedra, it has a layering onto the univalent edge of the
$3/4$--layered-triangulation, the edge labeled $7$, and has two
edges having the same slope. One is the thick edge and the other is
the new edge meeting the meridional disk just once. This is a
nearly-minimal layered-triangulation. Notice that we get a different
one for each point $1/n$ between, and including $1/6$ and $1/2$;
that is, in this case five different nearly-minimal, non minimal
layered-triangulations. However, one can not repeat these at more
than one point or extend beyond the points $n/n+1$. Since the edge
$1$ is the thick edge, there are not two thin edges having the same
slope.

\end{enumerate}

\section{Layered-triangulations of lens spaces}

In this section we define layered-triangulations of a lens spaces,
as well as minimal and nearly-minimal layered-triangulations of lens
spaces. We show that every lens space has a layered-triangulation
(in fact, infinitely many); each lens space, except $\mathbb{R}P^3$
and $S^2\times S^1$, has a $0$--efficient layered-triangulation with
an arbitrarily large number of tetrahedra; and each lens space,
having a $0$--efficient triangulation, except $S^3$, has only a
finite number of $1$--efficient layered-triangulations. In
particular, the (unique) minimal layered-triangulation of these lens
spaces is $1$--efficient.

\begin{figure}[h]
\psfrag{e}{$e$}\psfrag{i}{\footnotesize
{simplicial}}\psfrag{j}{\footnotesize {isomorphism}}
\psfrag{p}{\footnotesize{ pinching}}\psfrag{q}{\footnotesize{ about
$e$}} \psfrag{f}{\footnotesize{ folding}}\psfrag{g}{\footnotesize{
along $e$}}{\epsfxsize = 3.5 in
\centerline{\epsfbox{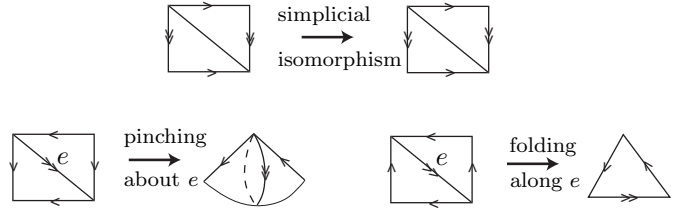}} } \caption{Simplicial
attaching maps from the boundary of a non degenerate
layered-triangulation of the solid torus to the boundary of a
layered-triangulation of the solid torus, including the degenerate
creased $3$--cell and the M\"obius band. } \label{f-fold-pinch}
\end{figure}

\vspace{.15 in}\noindent {\it Definition.}  Having layered
triangulations of the solid torus, then the obvious method for
getting similar triangulations of a lens space is to exploit the
genus one Heegaard decompositions of lens spaces. In doing this we
include the degenerate solid tori, the creased $3$--cell and the
one-triangle M\"obius band, as factors in a Heegaard decomposition
of a lens space. For this we need to specify what we mean by a
simplicial attachment of two layered-triangulations of the solid
torus along their boundaries when one of the factors is a degenerate
solid torus. If both factors are non degenerate, the attaching map
is a simplicial isomorphism. If one factor is non degenerate and the
other is a creased $3$--cell, then the simplicial map is determined
by the edge $e$ on the boundary of the non degenerate
layered-triangulation of the solid torus identified with the
univalent edge on the boundary of the creased $3$--cell; we call
this a {\it pinching} about the edge $e$. If one factor is non
degenerate and the other is the one-triangle M\"obius band, then the
simplicial map is determined by the edge $e$ on the boundary of the
non degenerate layered-triangulation of the solid torus identified
with the boundary edge of the M\"obius band; we call this a {\it
folding} along the edge $e$. Finally, we have the possibilities
where both solid tori are degenerate. If both are creased
$3$--cells, the double does not give a $3$--manifold and the one
remaining simplicial attachment gives one of the $2$ two-tetrahedra,
one-vertex triangulations of $S^3$. If one of the factors is a
creased $3$--cell and the other a one-triangle M\"obius band, then
we do not get a manifold if we take the univalent edge to the
boundary of the M\"obius band, leaving only one possibility which
gives us the one-tetrahedron, one-vertex triangulation of the
$3$--sphere. In Figure \ref{f-fold-pinch}, we give the options, just
described, for a simplicial attachment of two layered-triangulations
of the solid torus. Typically, a layered-triangulation of a lens
space can be viewed in a number of different ways as a union of two
layered solid tori.

\vspace{.125in}\noindent {\bf Examples.}
\begin{enumerate}
\item  Suppose we have a one-tetrahedron solid
torus. There are three ways to identify the boundary torus to the
one-triangle M\"obius band. For any one of these, we get a lens
space with a one-tetrahedron triangulation. The results are
one-tetrahedron triangulations for the lens spaces $L(1,1) = S^3$,
$L(4,1)$, and $L(5,2)$. These are layered (and minimal)
triangulations. The example that gives $S^3$ can also be viewed as
an identification of the two faces of the creased $3$--cell. All of
these examples are $1$--efficient triangulations.\item  Up to
isomorphism, there is only one triangulation of a $3$--manifold
obtained by attaching two creased $3$--cells along their boundaries.
We get a two-tetrahedron, one-vertex triangulation of $S^3$, which
is $1$--efficient. \item  There are $4$ two-tetrahedra solid tori:
the two-tetrahedra triangulations of the solid torus extending the
$0/1, 1/1, 1/3$ and $2/3$ triangulations on the boundary of the
solid torus.

First, we consider the identifications of the two faces in the
boundary of the two-tetrahedron extension of $0/1$. We have: $\rppp
= L(2,1)\equiv\{0,1,1'\}\lra \{2,1,1\}$ (we fold over the edge
labeled $0$); $S^3 = L(1,0)\equiv \{0,1,1'\}\lra \{1,2,1\}$; and
$S^3= L(1,1)\equiv \{0,1,1'\}\lra \{1,1,2\}$, where we use the label
$1'$ to denote the thin (univalent) edge rather than the thick edge.
By reversing our view in the first example giving $S^3$, we see that
it is the same as attaching a creased $3$--cell to the boundary of
the one-tetrahedron solid torus with the edge labeled $0$ in the
creased $3$--cell attached to the thick edge in the one-tetrahedron
solid torus. The last example is the same as two creased $3$--cells
identified along their boundaries. These last two examples are non
isomorphic two-tetrahedron triangulations of $S^3$, which has $5$
non isomorphic two-tetrahedron triangulations. These two are
$1$--efficient.

For the two-tetrahedron triangulation of the solid torus extending
$1/1$, we have: $L(3,1)\equiv\{1',1,2\}\lra \{2,1,1\}$;
$L(3,1)\equiv\{1',1,2\}\lra \{1,2,1\}$; and $S^2\times S^1 =
L(0,1)\equiv\{1',1,2\}\lra \{1,1,2\}$, where, again, we use $1'$ to
denote the thin (univalent) edge. These two triangulations of
$L(3,1)$ are not isomorphic and give two of the three non isomorphic
triangulations of $L(3,1)$. The first one of these is not
$0$--efficient; we are folding over the univalent edge and the
triangulation is the same as attaching the creased $3$--cell to the
one-tetrahedron solid torus with the edge labeled $0$ in the creased
$3$--cell identified to the edge labeled $3$ in the one-tetrahedron
solid torus (we have an edge bounding a disk). The second is
$1$--efficient.

For the two-tetrahedron triangulation of the solid torus extending
$1/3$, we have: $L(7,3)\equiv\{1,3,4\}\lra \{2,1,1\}$;
$L(5,1)\equiv\{1,3,4\}\lra \{1,2,1\}$; and $\rppp =
L(2,1)\equiv\{1,3,4\}\lra \{1,1,2\}$. This is one of the $2$
two-tetrahedron triangulations of $\rppp$ and is isomorphic to the
one above.

For the two-tetrahedron triangulation of the solid torus extending
$2/3$, we have: $L(8,3)\equiv\{2,3,5\}\lra \{2,1,1\}$;
$L(7,2)\equiv\{2,3,5\}\lra \{1,2,1\}$; and $S^3 =
L(1,2)\equiv\{2,3,5\}\lra \{1,1,2\}$. This triangulation of $L(7,2)
= L(7,3)$ is equivalent to that above and the last triangulation of
$S^3$ is isomorphic to the one above with $S^3 = L(1,0)\equiv
\{0,1,1\}\lra \{1,2,1\}$.\item  If we have a layered solid torus,
then we can attach a creased $3$--cell to its boundary via a
simplicial map. This is the same as adding a $2$-- and $3$--handle
to the solid torus, where the $2$--handle is attached along the
slope of an edge in the boundary of the layered solid torus. For
example, if we have the one-tetrahedron solid torus, an extension of
the $1/2$--triangulation on the boundary of the solid torus, then
attaching the creased $3$--cell gives: $L(1,0) = S^3$ for the edge
labeled $0$ in the boundary of the creased $3$--cell attached to the
edge labeled $1$; $L(2,1) = \rppp$ for the edge $0$ attached to the
edge $2$; and $L(3,1)$ for the edge $0$ attached to the edge $3$.
These give examples of two-tetrahedra, layered-triangulations for
$S^3$, $\rppp$ and $L(3,1)$. These two-tetrahedra
layered-triangulations for $\rppp$ and $L(3,1)$ are minimal but
neither is $0$--efficient. The one for $S^3$ is $1$--efficient (this
distinction will be discussed below). Recall that except for these
minimal triangulations of $\rppp$ and $L(3,1)$, a minimal
triangulation of an irreducible $3$--manifold must be
$0$--efficient.

\item A more typical example of a layered-triangulation of a lens space
is given as $\{9,7,2\}\lra\{5,3,8\} = L(62,27)$, meaning any
layered-triangulation of a solid torus extending $2/7$ identified
with any layered-triangulation of a solid torus extending $3/5$ so
that $9\lra 5, 7\lra 3$ and $2\lra 8$ gives the lens space
$L(62,27)$.
\end{enumerate}

\vspace{.15 in}\noindent {\it Existence.} The next lemma follows
from the fact that we can extend any one-vertex triangulation on the
boundary of a solid torus to a layered-triangulation of the solid
torus, along with the fact that lens spaces have a genus one
Heegaard splitting.

\begin{thm} Every lens space admits a layered-triangulation.
\end{thm}

\begin{proof} Start with a genus one Heegaard splitting; fix a one-vertex
triangulation of the splitting torus. By Theorem
\ref{extend-to-layered}, we can extend the one-vertex triangulation
of the splitting torus to a layered-triangulation of the solid torus
on each side.\end{proof}

\vspace{.15 in}\noindent {\it Combinatorics of
layered-triangulations of lens spaces.} As mentioned above, a fixed
layered-triangulation of a lens space, typically,  can be viewed as
a simplicial attachment along the boundaries of two
layered-triangulations of the solid torus in a number of ways. On
the other hand, a layered-triangulation of a lens space can be
viewed as being formed from a layered-triangulation of a single
solid torus along with one of the three possible ways of
simplicially attaching the boundary torus to the one-triangle
M\"obius band; i.e., folding along one of the three edges in the
boundary.  This gives us a very nice way to address the
combinatorics of what lens spaces we get from a given
layered-triangulation of the solid torus and conversely, given the
relatively prime integers $X$, $Y$, we provide an algorithm to
construct a layered-triangulation of the lens space $L(X,Y)$.

The combinatorics we develop work equally well for determining the
lens space obtained by folding along one of the three edges in the
boundary of any one-vertex triangulation of a solid torus. So, to
this end, suppose $\T$ is a one-vertex-triangulation of the solid
torus extending the $p/q$--triangulation  on the boundary torus.
There are three edges in the boundary and we can choose a labeling
and an orientation of the edges $e_1,e_2,e_3$ so that the meridional
disk meets these edges in $p, q$ and $p+q$ points, respectively. The
one-triangle M\"obius band can be considered as a degenerate solid
torus having ``meridian slope" the nontrivial spanning arc and
``longitude" the interior edge in the one-triangle triangulation. It
follows that the meridian of the M\"obius band lifts to a normal
curve in the boundary of the solid torus meeting one edge twice and
each of the other two edges in the boundary once each. We can think
of the edge being met twice as the edge in the boundary of the solid
torus along which we are folding, see Figure \ref{f-fold-pinch}. It
is the edge identified with the boundary edge in the one-triangle
M\"obius band. The other two edges in the boundary of the solid
torus are identified to the single interior edge in the M\"obius
band, which functions as a longitude for the M\"obius band.

We have that $e_1$ and $e_2$ form a basis for the homology of the
boundary torus and $e_3 = e_1 - e_2$. Thus the meridional slope,
$\mu'$, for the one-vertex triangulation of the solid torus can be
expressed in terms of the basis $\{e_1, e_2\}$ as $\mu' = qe_1 +
pe_2$. There are three possibilities for identifying the two faces
in the boundary torus, depending on which edge, $e_1, e_2$, or
$e_3$, is the edge along which we are folding.

\vspace{.125 in}\noindent {\bf Type I.} Folding along $e_1$. Then
the meridional slope for the M\"obius band lifts to the curve $\mu =
e_1 - 2e_2$; $\lambda = e_2$ and $\hat{\lambda} = e_3 = e_1 - e_2$
are lifts of the longitude in the M\"obius band. It follows that
$d(\mu,\mu')$, the distance (or number of intersections) between
$\mu$ and $\mu'$, is $p + 2q$; whereas, $d(\lambda, \mu') = q$ and
$d(\hat{\lambda},\mu') = p+q$. Thus using the longitude $\lambda$,
we conclude that we have the lens space $L(p + 2q, q)$; and using
the longitude $\hat{\lambda}$, we conclude that we have the lens
space $L(p + 2q, p+q)$. Of course, these are the same lens space;
$p+q$ is the additive inverse of $q$,$ \mod(p+2q)$.

\vspace{.125 in}\noindent {\bf Type II.} Folding along $e_2$. This
is exactly the same as Type I with a change in the role  of $p$ and
$q$. For this case, the meridional slope for the M\"obius band lifts
to the curve $\mu = 2e_1 - e_2$; $\lambda = - e_1$ and
$\hat{\lambda} = e_3 =
 e_1- e_2$ are lifts of the longitude in the M\"obius band. It follows
that $d(\mu,\mu') = 2p + q$; whereas, $d(\lambda, \mu') = p$ and
$d(\hat{\lambda},\mu') = p+q$. Thus using the longitude $\lambda$,
we conclude that we have the lens space $L(2p + q, p)$; and using
the longitude $\hat{\lambda}$, we conclude that we have the lens
space $L(2p + q, p+q)$. As above, these are the same lens space.

\vspace{.125 in}\noindent {\bf Type III.} Folding along $e_3$. Now,
the meridional slope for the M\"obius band lifts to the curve $\mu =
e_1 + e_2$; and $\lambda = e_1$ and $\hat{\lambda} = -e_2$ are the
lifts of the longitude in the M\"obius band. We have in this case
that $d(\mu,\mu')= \abs{q-p}$; whereas, $d(\lambda, \mu') = p$ and
$d(\hat{\lambda},\mu') = q$. Thus using the longitude $\lambda$, we
conclude that we have the lens space $L(\abs{q-p},p)$; and using the
longitude $\hat{\lambda}$, we conclude that we have the lens space
$L(\abs{q-p}, q)$. These are the same lens space as $p$ is
equivalent to $q, \mod(\abs{q-p})$.

\vspace{.15 in}Conversely, suppose we wish to construct a
layered-triangulation for the lens space $L(X,Y)$.

If $0\le Y\le X/2$, then choose a layered-triangulation of the solid
torus extending the one-vertex triangulation on the boundary whose
edges, $e_1,e_2,e_3$, meet the meridional slope in $X-2Y,Y,X-Y$
points, respectively, and fold over the edge $e_1$ (with
intersection number $X-2Y$).

If $X/2\le Y\le X$, then choose a layered-triangulation of the solid
torus extending the one-vertex triangulation on the boundary whose
edges, $e_1,e_2,e_3$, meet the meridional slope in $X-Y,2Y-X,Y$
points, respectively, and fold over the edge $e_2$ (with
intersection number $2Y-X$).

If $0< X< Y$, then choose a layered-triangulation of the solid torus
extending the one-vertex triangulation on the boundary whose edges,
$e_1,e_2,e_3$, meet the meridional slope in $Y,Y-X,2Y-Y$ points,
respectively, and fold over the edge $e_3$ (with intersection number
$2Y-X$). More generally, for any $X>0, Y\ge 0$, choose a
layered-triangulation of the solid torus extending the one-vertex
triangulation on the boundary whose edges, $e_1,e_2,e_3$, meet the
meridional slope in $Y,X+Y,X+2Y$ points, respectively, and fold over
the edge $e_3$ (with intersection number $X+2Y$). Notice that in
each of these last two cases, we are folding over the edge with
highest intersection number with the meridional slope; hence, these
triangulations admit a reduction in the number of tetrahedra.

Note that if we have a minimal layered-triangulation of a solid
torus, then the univalent edge is the furthest from the meridional
slope (meets the meridional slope the greatest number of times) and
if we were to fold over the univalent edge, then, typically, we
could reduce the number of tetrahedra in the triangulation. Hence,
the preferred situation is that we have a
$p/q$--layered-triangulation and we fold over one of the edges
labeled $p$ or $q$. The following is an easy mnemonic. If we have a
$p/q$-triangulation on the boundary of the solid torus and we fold
over the edge labeled $p$; i.e., the lift of the meridian for the
one-triangle M\"obius band meets the edge labeled $p$ two times,
then set down the array

\centerline{\(\begin{array}{cc} $2$ &$1$\\
$p$ &$q$\\
\end{array}\)}

\noindent cross multiply and add, getting $2q + p$. This is $X$ in
the resulting lens space $L(X,Y)$ and $Y = q$ is the value in the
column with $1$. If we use as the longitude the edge labeled $p+q$
rather than the edge labeled $q$, then we have the array

\centerline{\(\begin{array}{cc} $2$ &$1$\\
$p$ &$p+q$\\
\end{array}\)}

\noindent and in this case we take the difference $2p+2q - p$,
after cross multiplying, and have $X = 2q + p$; but now we have $Y
= p+q$, which is the value in the column with $1$.

\vspace{.125 in}\noindent{\bf Example.} Suppose we have a
$4/5$--triangulation on the boundary of the solid torus. If we
obtain a lens space by folding over the edge labeled $5$, we have
the arrays

\centerline{\(\begin{array}{cc} $2$ &$1$\\
$5$ &$4$\\
\end{array}\) or \(\begin{array}{cc} $2$ &$1$\\
$5$ &$9$\\
\end{array}\)}

\noindent depending on a choice of longitude. In the first array,
by cross multiplying and adding, we see that we get the lens space
$L(13,4)$; whereas, in the second array, we cross multiply and
take the difference and have the lens space $L(13,9)$.

If we fold over the edge labeled $4$ and use as longitude the edge
labeled $5$, we have an array

\centerline{\(\begin{array}{cc} $2$ &$1$\\
$4$ &$5$\\
\end{array}\)}

\noindent and get the lens space $L(14,5)$.

If we have a layered-triangulation of a lens space obtained by
folding along an edge in the boundary of a
$p/q$--layered-triangulation, then there is a symmetry; namely, once
we have folded over, then we have a layered-triangulation of a solid
torus in the reverse direction and have obtained the
layered-triangulation of the lens space by folding along an edge in
its boundary. We describe the combinatorial relationship of looking
at a layered-triangulation of a lens space by considering it as
obtained by folding at either end. The layered-triangulation of the
solid torus in the reverse direction is a
$p'/q'$--layered-triangulation for some $p', q'$. We call these two
views of the same layered-triangulation of a lens space as {\it
reversed layerings}.

Since a layered-triangulation extending $p/q$ on the boundary can be
viewed as a path in the $L$--graph, it has inside it the path of the
minimal layered-triangulation extending $p/q$. This provides us with
an algorithm to determine $p'/q'$ from $p/q$. We show this in a
couple of examples.

First we establish some notation. If we have a layered-triangulation
extending a one-vertex triangulation on the boundary where the
intersection numbers of the meridional slope with the three edges is
the triple $\{x,y,z\}$, we use $\{x,\underline{y},z\}$ to denote
that the edge in the triangulation of the boundary meeting the
meridional slope $y$ times is the univalent edge for our layering.
If we have obtained a lens space by folding over the edge labeled
$x$, then we can think of the simplicial identification
$\{x,y,z\}\lra\{2,1,1\}$ of the boundary of the solid torus with the
one-triangle M\"obius band, where we have the edge labeled $x$
identified with the boundary of the M\"obius band and the edges $y$
and $z$ identified with the interior edge. Now, using the above
notation with $y$ the univalent edge, we have
$\{x,\un{y},z\}\lra\{2,\un{1},1\}$ or $\{x,y',z\}\lra\{2,3,1\}$,
where $y' = x+z$ or $y' = \abs{x-z}$.

\vspace{.125 in}\noindent{\bf Examples:} \vspace{-.05 in}
 \begin{tabbing}
 \= \bf{(a)} \=$\{\un{1},2,1\} \lra \{\un{14},3,11\} = L(25,11)$
\hspace{.5 in} \=\bf{(b)} \=$\{\un{1},2,1\} \lra
\{\un{5},3,8\} = L(13,8)$\\
\>  \>$\{3,2,\un{1}\} \lra \{8,3,\un{11}\} $ \hspace{.75 in} \>
\>$\{\un{3},2,1\} \lra
\{\un{11},3,8\} $\\
\> \>$\{\un{3},2,5\}\lra\{\un{8},3,5\}$  \> \>$\{1,2,\un{1}\}\lra\{5,3,\un{8}\}$\\
\> \>$\{7,2,\un{5}\}\lra\{2,3,\un{5}\}$  \> \>$\{\un{1},2,3\}\lra\{\un{5},3,2\}$\\

\> \>$\{7,2,9\}\lra\{2,\un{3},1\} $   \> \>$\{5,2,3\}\lra\{1,\un{3},2\} $\\
\> \>$\{7,16,9\}\lra\{2,1,1\} = L(25,9)$   \> \>$\{5,8,3\}\lra\{1,1,2\} = L(13,5)$\\
 \end{tabbing}

More generally, suppose we have two solid tori and
layered-triangulations of each. If we identify the boundaries of the
two solid tori via a simplicial isomorphism, we have a
layered-triangulation of a lens space. The above method provides us
with an algorithm to determine the lens space obtained by knowing
the original layered-triangulations. We exhibit this in the
following example.

\vspace{.25 in}\noindent{\bf (c)} Suppose we have
$\{9,7,2\}\lra\{3,7,4\}$, where the notation means that we have
identified the boundaries of two solid tori, one has a
layered-triangulation extending $2/7$, the other has a
layered-triangulation extending $3/4$ and the identification takes
the edges $9\lra 3, 7\lra 7$ and $2\lra 4$.
 \begin{tabbing}
 \=\=$\cdots\cdots\cdots\cdots\cdots\cdots$\=\hspace{.375 in}\=\hspace{.375 in}\=\hspace{.25 in}
\=$\{31,11,20\} \lra \{1,1,2\} = L(42,11)$\\
 \>\>$\cdots\cdots\cdots\cdots\cdots\cdots$\>\>\>
\>$\{9,11,20\} \lra \{3,1,2\} $\\
 \>\>$\cdots\cdots\cdots\cdots\cdots\cdots$\>\>\>\>$\{9,11,\un{2}\}\lra\{3,1,\un{4}\}$\\
\> \>$\{\un{9},7,2\} \lra \{\un{3},7,4\}$ \>\>$\equiv$ \>\>
$\{9,\un{7},2\} \lra
\{3,\un{7},4\}$\\
\> \>$\{\un{5},3,2\}\lra\{\un{11},15,4\}$  \> \>\>\>$\cdots\cdots\cdots\cdots\cdots\cdots$\\
\> \>$\{1,3,2\}\lra\{19,\un{15},4\} $ \> \> \>\>$\cdots\cdots\cdots\cdots\cdots\cdots$\\
\> \>$\{1,1,2\}\lra\{19,23,4\} = L(42,19)$ \> \> \>\>$\cdots\cdots\cdots\cdots\cdots\cdots$\\
 \end{tabbing}

Of course, we have that our combinatorial description of the lens
space obtained by attaching two layered solid tori gives in one
reduction the lens space description as L(X,Y) and by {\it
reversing} the description we get $L(X,Y')$, where $Y\cdot Y'\equiv
\pm 1 mod(X)$.

We summarize our observations in the following two propositions.

\begin{prop} Any lens space can be formed from identifications of Type I (or Type
II).
\end{prop}

\begin{prop}If the lens space $L(X,Y) = \{x,y,z\}\lra\{1,1,2\}$,
then the reverse identification $\{1,1,2\}\lra\{x',y',z'\} =
L(X,Y')$ and $Y\cdot Y' \equiv \pm 1 mod(X)$.\end{prop}

A layered-triangulation of a lens space is said to be a {\it
minimal} layered-triangulation if it has the minimal number of
tetrahedra among all layered-triangulations of the lens space.
Clearly, if we have a minimal layered-triangulation of a lens space,
then it is obtained by folding along an edge in the boundary of a
minimal layered-triangulation of a solid torus and, in general, the
edge along which the boundary is folded is not the univalent edge
(the only exception is the one-tetrahedron layered-triangulation of
$S^3$ obtained by folding the one-tetrahedron solid torus along the
univalent edge.  Later, by using the uniqueness of Heegaard
splittings of lens spaces, we will see that if a
layered-triangulation of a lens space is obtained by folding the
boundary of a minimal $p/q$--layered-triangulation of the solid
torus along an edge other than the univalent edge, then the
triangulation is a minimal layered-triangulation of the lens space.
At this point, we could employ the classification of lens spaces and
the last proposition of the previous section for this conclusion;
however, uniqueness of Heegaard splittings gives us the standard
classification of lens spaces.  It seems to us that these results
should follow from the combinatorics above but this is not clear.
Finally, with the exception of $L(3,1)$ there is a unique minimal
layered-triangulation of a lens space.

\begin{prop} Every lens space has a minimal layered-triangulation.
\end{prop}

There is another curious question which we have not been able to
answer. Namely, in general, is a minimal triangulation of a lens
space a layered-triangulation? There are three exceptions, $S^3$,
$\rppp$ and $L(3,1)$; each has a minimal triangulation that is not
layered. $S^3$ has a one-tetrahedron layered-triangulation and also
has a one-tetrahedron, two-vertex triangulation. Both $\rppp$ and
$L(3,1)$ have minimal, two-tetrahedra triangulations that are not
layered. They also have two-tetrahedra layered-triangulations. $S^3$
has 5 distinct two-tetrahedra triangulations, 3 of which are not
layered. Triangulations of these three manifolds also provide the
exceptions in the reduction of an arbitrary triangulation to a
$0$--efficient one \cite{jac-rub0}.

\begin{conj2} With the exceptions of $S^3, \rppp$ and $L(3,1)$, a minimal
triangulation of a lens space is a minimal layered-triangulation.
\end{conj2}

This has also been conjectured by S. Matveev \cite{matveev2}.

As a final observation, every layered triangulation of a lens space
L(p,q) is invariant under the canonical involution. The quotient
space of this involution is the $3$-sphere and the fixed set of the
involution projects to a $2$-bridge knot or link. Therefore we can
view the classification of minimal layered-triangulations as a
version of the classification of 2-bridge knots and links.

\section{Normal and almost normal surfaces in layered-triangulations of lens spaces}

We first classify the orientable normal and almost normal surfaces
in a minimal layered-triangulation of a lens space; then we classify
the nonorientable normal surfaces. The last uses our
layered-triangulations for a new proof of earlier work by G. Bredon
and J. Wood, \cite{bredon-wood}. From now on we shall refer to the
two-tetrahedra layered-triangulation of $L(3,1)$ that is obtained by
folding over along the univalent edge in the two-tetrahedra
$1/1$--layered-triangulation as the {\it bad} minimal
layered-triangulation of $L(3,1)$. Recall that we also get $L(3,1)$
if we fold over along the thick edge; however, this latter
triangulation of $L(3,1)$ does not admit the unnecessary embeddings
of normal and almost normal surfaces as does the former.

\begin{thm2}\label{surface-in-lens} The embedded, orientable, normal
and almost normal surfaces in a minimal layered-triangulation of a
lens space are:

If normal:
\begin{itemize}

\item  the vertex-linking $2$--sphere, possibly, with thin
edge-linking tubes,  except \item a non vertex-linking normal
$2$--sphere for the two-tetrahedra layered-triangulations of
$\rppp$, $S^2\times S^1$, and the bad minimal layered-triangulation
of $L(3,1)$.
\end{itemize}

If almost normal and not normal:
\begin{itemize}

\item  an  almost normal tubed surface obtained either by adding an almost
normal tube to the vertex-linking $2$--sphere (possibly) with thin
edge-linking tubes, or by adding an almost normal tube between two
such surfaces, except \item  an  almost normal octagonal
$2$--sphere, in the minimal layered-triangulations of $S^3$,
$\rppp$, and the bad minimal layered-triangulation of $L(3,1)$, and
\item an almost normal octagonal  torus in the bad two-tetrahedra
layered-triangulation of $L(3,1)$.
\end{itemize}

\end{thm2}
\begin{proof} A layered-triangulation of a lens space is obtained from
a layered-triangulation of a solid torus by folding along its
boundary. Hence, a normal or almost normal surface in the
layered-triangulation of the lens space,  must come from a normal or
almost normal surface embedded in the layered-triangulation of the
solid torus after some identifications in its boundary. In
particular, the boundary of the surface in the solid torus must be
invariant under the folding. It follows that if the edges in the
boundary of a layered solid torus are denoted $e_1,e_2,e_3$ and if
we fold over along the edge $e_3$, say, then the number of
intersections of the normal or almost normal surface in the layered
solid torus with the edges labeled $e_1$ and $e_2$ must be the same.
Hence, following the classification of normal and almost normal
surfaces in layered solid tori, the only possible curves, given by
their intersection number with $e_1, e_2$ and $e_3$, respectively,
are: $1,1,0$; $1,1,2$; and $2,2,2$. These very restricted
possibilities are the reason for the limited types of normal and
almost normal surfaces in layered-triangulations of lens spaces.

We shall defer consideration of the normal and almost normal
surfaces in the small, exceptional triangulations obtained by
folding along an edge in the boundaries of the minimal $0/1$--,
$1/1$-- and $1/2$--layered-triangulations until we have done the
general case. Otherwise, for a  minimal layered-triangulation of a
solid torus:
\begin{itemize}\item[-] An orientable normal
surface whose boundary curves have intersection numbers $1,1,0$ is a
thin edge-linking annulus about the edge $e_3$ (possibly) with thin
edge-linking tubes.  Hence, upon folding the boundary of the solid
torus, such a surface gives a vertex-linking $2$-sphere with thin
edge-linking tubes.  \item[-] An orientable normal surface whose
boundary curves have intersection numbers $1,1,2$ can only happen if
$e_3$ is  the univalent edge; however, in a minimal
layered-triangulation of a lens space, we do not fold over along the
univalent edge.\item[-] An orientable normal surface whose boundary
curves have intersection numbers $2,2,2$ is  a vertex-linking disk
(possibly) with some thin edge-linking tubes. Hence, upon folding
the boundary of the solid torus, such a surface gives the
vertex-linking $2$--sphere, possibly, with some thin edge-linking
tubes. For intersection $2,2,2$ we have the same surfaces under
consideration no matter which edge we fold along.\end{itemize}

Thus, to complete our proof, we only need to consider lens spaces
obtained by folding the boundary of the one-tetrahedron solid torus
extending $1/2$, or the two two-tetrahedra solid tori extending
$1/1$ or $0/1$. The normal surfaces in these layered-triangulations
are given in Lemmas \ref{one-tet-torus}, \ref{0-1-1}, and in the
case of $1/1$ in the final paragraph of the proof of Theorem
\ref{normal}.

\vspace{.1 in}\noindent For the minimal
$1/2$--layered-triangulation:

This is the one-tetrahedron solid torus. Folding along the edge $3$
(the univalent edge), we get a one-tetrahedron triangulation of
$S^3$. There is a thin edge-linking annulus about the edge $3$
having boundary with intersection numbers $1,1,0$ that gives the
vertex-linking 2-sphere with a thin edge-linking tube about the edge
$3$.  There is  an almost normal octagonal annulus having boundary
with intersection numbers $1,1,2$ (the slope of the exceptional
longitude) that gives an  almost normal octagonal $2$--sphere. The
vertex-linking curve with intersection numbers $2,2,2$ bounds only
the vertex-linking disk and gives the vertex-linking $2$--sphere.
Hence, the exception in this case is an almost normal octagonal
$2$--sphere.

Folding along the edge $2$, we get the one-tetrahedron triangulation
of $L(4,1)$. The boundary of a normal M\"obius band has intersection
numbers $1,0,1$ and gives us a Klein bottle, which we consider
below; its double gives the vertex-linking $2$-sphere with a thin
edge-linking tube about the edge $2$. There are no surfaces having
boundary with intersection numbers $1,2,1$. The intersection $2,2,2$
gives the vertex-linking $2$--sphere.

Folding along the edge $1$ (the thick edge), we get the
one-tetrahedron triangulation of $L(5,2)$. There are no normal
surfaces with intersection numbers $0,1,1$ or $2,1,1$. The
intersection $2,2,2$ gives the vertex-linking $2$--sphere.

\vspace{.1 in} \noindent For the minimal
$0/1$--layered-triangulation:

Folding along the univalent (thin) edge $1$ or the thick edge $1$,
gives two distinct two-tetrahedra triangulations of $S^3$. These are
not minimal; however, the examples are interesting so we briefly
list their orientable normal and almost normal surfaces.

First, we consider the case of folding along  the univalent (thin)
edge.  For boundary with intersection $0,1,1$, we have a thin
edge-linking annulus about the univalent edge  and this same annulus
with a thin edge-linking tube about the edge $2$; these give the
vertex-linking $2$--sphere with a thin edge-linking tube about the
univalent edge and this same surface with a second tube about the
edge $2$. There also is an  almost normal octagonal annulus having
boundary slope the same as the univalent edge. This gives an almost
normal octagonal torus; however, we do not list this among the
conclusions of the theorem because this is not a minimal layered
triangulation for $S^3$.  For intersection numbers $2,2,2$ we have
the vertex-linking disk and the vertex-linking disk with a thin
edge-linking tube about the edge $2$; these give the vertex-linking
$2$--sphere and the vertex-linking $2$--sphere with a thin
edge-linking tube, respectively. We also have a  almost normal
octagonal disk with trivial boundary curve; this gives an almost
normal octagonal $2$--sphere.

If we fold along the thick edge, we only get surfaces with the slope
$2,2,2$ and we have the vertex-linking $2$-sphere, the
vertex-linking $2$--sphere with a thin edge-linking tube on the edge
$2$ and an  almost normal octagonal $2$--sphere.

Folding along the edge labeled $0$, we get $\rppp$. Here we are
concerned with surfaces having boundary slope one of $2,2,2$ or
$1,1,0$; there are no surfaces with boundary slope $1,1,2$. For
$2,2,2$ we get the vertex-linking $2$--sphere, the vertex-linking
$2$--sphere with a thin edge-linking tube along the edge $2$ and an
 almost normal octagonal $2$--sphere. For $1,1,0$, we have the
meridional disk in the solid torus, which leads to the projective
plane in $\rppp$ and its double is a normal $2$--sphere. The twice
punctured Klein bottle, also with boundary intersection $1,1,0$,
leads to a non orientable surface having genus $3$; its double is a
genus two orientable surface, which is the vertex-linking
$2$--sphere with thin edge-linking tubes about edges $0$ and $2$.

\vspace{.1 in}\noindent For the minimal
$1/1$--layered-triangulation::

Folding along the edge $2$ is the same as doubling the
one-tetrahedron solid torus and we get $S^2\times S^1$. Here we have
all three possible slopes. For $2,2,2$, we get the vertex-linking
$2$--sphere and the vertex-linking $2$--sphere with a thin
edge-linking tube on the edge $3$. For $1,1,2$, we have the
meridional disk which gives the essential $2$--sphere; and for
$1,1,0$, the M\"obius band doubles to a Klein bottle and its double
is the vertex-linking $2$--sphere with a thin edge-linking tube
about the edge $2$. There are no  almost normal octagonal surfaces.

Folding along either the thick edge  or the univalent edge, we get
the two distinct minimal layered-triangulations of $L(3,1)$. In the
first case,  folding along the thick edge $1$, we only have the one
slope $2,2,2$ and thus we have the vertex-linking $2$--sphere and
the vertex-linking $2$--sphere with a thin edge linking tube about
the edge $3$. In the last case, folding along the univalent edge, we
get the bad minimal layered-triangulation of $L(3,1)$ and have all
three sets of intersection numbers possible. As above, for $2,2,2$,
we get the vertex-linking $2$--sphere and the vertex-linking
$2$--sphere with a thin edge-linking tube along the edge $3$. For,
intersection numbers $2,1,1$, we have an edge-linking annulus about
the interior edge $3$, which gives an embedded normal $2$--sphere;
and an  almost normal octagonal annulus, which gives an almost
normal octagonal $2$--sphere. Finally, we have intersection numbers
$0,1,1$ and get a vertex-linking $2$--sphere with thin edge-linking
tube along the thin edge labeled $1$, the vertex-linking $2$--sphere
with thin edge-linking tubes along the thin edge labeled $1$ and
along the edge labeled $3$ and an  almost normal octagonal torus.
The last example is \emph{not} $0$--efficient.\end{proof}

Notice that we have almost normal octagonal $2$--spheres in the
minimal layered-triangulations of $S^3$, $\rppp$ and  the bad
layered-triangulations of $L(3,1)$. Every triangulation of $S^3$
must have an almost normal $2$--sphere and a one-vertex
triangulation  must have such an octagonal  $2$--sphere. The
triangulations of $\rpp$ and the bad minimal layered-triangulation
for $L(3,1)$ are not $0$--efficient; thus they have non
vertex-linking normal $2$--spheres. The  almost normal octagonal
$2$--spheres are formed as the minimax in a sweep-out between these
non vertex-linking normal $2$--sphere and the vertex-linking
$2$--sphere. See more on this below.

\begin{cor} For a lens space distinct from $S^3$, $\rppp$, $L(3,1)$, and
$S^2\times S^1$, the only embedded, orientable normal surfaces are
the vertex-linking $2$--spheres with thin edge-linking tubes. There
are no almost normal octagonal surfaces.\end{cor}

\section{Applications of layered-triangulations of the solid torus and lens spaces}
In this section we use layered-triangulations to capture some
classical results about the topology of lens spaces. We obtain
earlier results of Bredon and Wood \cite{bredon-wood} characterizing
those lens spaces that admit embedded nonorientable surfaces and
classifying the embedded nonorientable surfaces in each such lens
space. We apply layered-triangulations to obtain new proofs of the
results of Waldhausen \cite{wald-HS3} and those of Bonahon and Otal
\cite{bon-otal} to classify Heegaard splitings of  $S^3$ and
$S^2\times S^1$ and all (other) lens spaces. We also use
layered-triangulations of the solid torus to construct what we
consider canonical triangulations for Dehn-fillings of cusped
manifolds (called triangulated Dehn fillings); in particular, we use
these triangulations as a tool to study Heegaard splittings of
Dehn-fillings of cusped manifolds. Finally, we give sufficient
conditions for $0$--efficient triangulations of cusped manifolds to
extend to $0$--efficient triangulated Dehn fillings of these
manifolds.

\subsection{Embedded, non orientable surfaces in lens
spaces.}

A determination of those lens spaces that admit embeddings of non
orientable surfaces and the classification of non orientable
surfaces embedded in such lens spaces is given in
\cite{bredon-wood}. Here we apply layered-triangulations of lens
spaces to obtain these results, including the cases of $S^3$ and
$S^2\times S^1$.

We shall use $U_h$ to denote the non orientable surface of genus $h$
(connected sum of $h$ copies of $\rpp$). In an orientable
$3$--manifold a non orientable surface is obtained by possibly
adding handles to an embedded, incompressible non orientable surface
or to an incompressible, non separating, orientable surface. For a
lens space, except  $S^2\times S^1$,  we only have  the first
possibility; for $S^2\times S^1$, we also have the second
possibility, interpreting the essential non separating $2$--sphere
as incompressible.

If the lens space $L(X,Y)$ is distinct from $S^2\times S^1$ and
contains a non orientable surface, then such a surface is isotopic
to an embedded, incompressible non orientable surface with possibly
some trivial handles attached; it follows that if $L(X,Y)$ contains
an embedded, non orientable surface, then for any triangulation,
there must be a \emph{normal} (incompressible) non orientable
surface and our given surface is isotopic to this normal surface
with possibly some trivial handles attached. For $S^2\times S^1$ a
similar statement applies where we have in place of the
incompressible nonorientable surface the ``incompressible", non
separating $2$--sphere. Of course, once the handles are attached,
there is no reason the surface need be normal or even isotopic to a
normal surface. In any case, to understand all embedded non
orientable surfaces, it is sufficient to understand the embedded,
\emph{normal}, non orientable surfaces in a triangulation of
$L(X,Y)$. To do this, we use a minimal layered-triangulation for
$L(X,Y)$.  As above for orientable surfaces, such a normal surface
in a minimal layered-triangulation of a lens space must come from an
embedded normal surface in a minimal layered-triangulation of the
solid torus with identifications in its boundary determined by the
folding of the boundary of the solid torus to get the lens space.
Furthermore, and again just as above, such a normal surface in the
layered-triangulation of the solid torus must have intersection
numbers with the edges $e_1, e_2$ and $e_3$ in the boundary of the
solid torus one of $1,1,0$; $1,1,2$; or $2,2,2$, respectively, where
we are, again, denoting the edge on which we fold by $e_3$.

A normal surface with intersection number $2,2,2$ can not give a non
orientable surface upon folding. A normal surface with intersection
numbers $1,1,2$ indicates that the normal surface meets the folding
edge $e_3$ maximally; but this can only happen if we fold on the
univalent edge or the layered-triangulation of the solid torus is
the minimal $1/1$--layered-triangulation and we fold on the edge
$2$. For a minimal layered-triangulation, we only fold on the
univalent edge when we have the one-tetrahedron
$1/2$--layered-triangulation and then we get $S^3$; if we fold on
the edge $2$ in the $1/1$--layered-triangulation, we get $S^2\times
S^1$. There are no normal surfaces in the one-tetrahedron
$1/2$-layered-triangulation with boundary slope having intersection
numbers $1,1,2$. Hence, except for $S^2\times S^1$, a non orientable
normal surface in a minimal layered-triangulation of a lens space
must come from a normal surface in a minimal layered-triangulation
of a solid torus having boundary parallel to the folding edge. From
our classification of normal surfaces in Theorem \ref{normal} such a
surface is either a thin edge-linking annulus, possibly with thin
edge-linking tubes, or a non orientable surface, or we have the
two-tetrahedra $0/1$ triangulation of the solid torus, the normal
surface is the meridional disk, and folding over along the edge $0$
gives $\rppp$. In the first case, the edge-linking annulus
identifies to a vertex-linking $2$--sphere with thin edge-linking
tubes, which is orientable. Thus, except for possibly $\rppp$ and
$S^2\times S^1$, we only get a non orientable normal surface in the
lens space, if we have a non orientable surface in the layered solid
torus with boundary slope the slope of the edge upon which we are
folding; in particular, the folding edge is even.

Now, we wish to determine the genus of the unique non orientable
normal surface in a minimal $p/q$--layered-triangulation of a solid
torus that has its boundary slope the slope of an edge of the
triangulation that is in the boundary; furthermore, we want the
genus in terms of the layered-triangulation of the solid torus. Let
$e(p,q)$ denoted the number of even order edges in the minimal
$p/q$--layered-triangulation of the solid torus. From Theorem
\ref{normal} we have that the genus of the surface of interest is
just $e(p,q)$. We use  the M\"obius band as the minimal
$1/1$--layered-triangulation; hence, we set $e(1,1)=1$; since we
have one edge of even order.

\begin{lem}\label{e-p-q}
 For, $p,q$ relatively prime, $q>p\ge 1$, we have \[ e(p,q) = \left\{\begin{array}{ll}
e(q-p,p) + 1& \mbox{if p and q are both odd,}\\
e(q-p,p)    & \mbox{if p or q is even.}
\end{array}
\right. \]
\end{lem}
\begin{proof} If we consider Euclid's original algorithm (or the
difference version of the Euclidean algorithm) for $p/q
=(a_0,\ldots,a_n)$ the continued fraction expansion of $p/q$, we get
a sequence of numbers that index the order of the edges in the
corresponding minimal $p/q$--layered-triangulation of the solid
torus. Hence, for $\{p,q\}$ we have the sequence:
$$p+q, q, q-p,\ldots,q-a_0p, p, p-r_1,\ldots,
p-a_1r_1,r_1,\ldots,2,1,$$ which contains the sequence for
$\{q-p,q,p\}$. So, if $p+q$ is even we have one more even edge in
the minimal $p/q$--layered-triangulation than we have for the
minimal layered-triangulation associated with the triple
$\{q-p,q,p\}$.\end{proof}

\vspace{.125 in}\noindent{\bf EXAMPLES:}\begin{enumerate}\item
$e(2,q) = 1$; $q$ is odd,\item $e(2k,1) = k$,\item $e(16,7) = e(9,7)
= e(2,7) + 1 = 2$.\end{enumerate}

In the statement of the following theorem, we exclude the lens
spaces $\rppp$ and $\S^2\times S^1$, both of which admit embedded
non orientable surfaces; the non orientable surfaces in these two
lens spaces arise in a different way from those in other lens
spaces. We provide the answers for $L(2,1) = \rppp$ and $L(0,1) =
S^2\times S^1$ in the examples following the proof of the theorem.

\begin{thm}\cite{bredon-wood} Suppose $L(X,Y)$ is a lens space distinct
 from $\rppp$ and $\S^2\times S^1$. Then $L(X,Y)$ admits an
embedded non orientable surface iff $X $ is even. Furthermore, there
is a unique embedded, \underline{incompressible} non orientable
surface $U_h$ in $L(X,Y)$, where $ h = e(X,Y)$. \end{thm}
\begin{proof} Assume $L(X,Y)$ has a minimal layered-triangulation.
If $L(X,Y)$ admits an embedded non orientable surface, then there is
an embedded, incompressible, non orientable surface in $L(X,Y)$;
hence, a normal such surface. From the above discussion, we have
that the minimal layered-triangulation for $L(X,Y)$ is derived by
folding over a minimal layered-triangulation of a solid torus where
the edge along which we fold has even order and is not univalent. We
denote the triangulation of the layered solid torus by
$\{2k,q,2k+q\}$. We fold over on the edge labeled $2k$; therefore,
$X = 2k + 2q$ and for $0<Y\le X/2$, $Y = q$, whereas, for $X/2\le
Y<X$,  $Y = 2k+q$. We have $X$ is even and  \[ \left\{
\begin{array}{lll}
k = X/2-Y & q = Y, &\mbox{if $0 < Y \le X/2,$}\\
k = Y-X/2 & q = X-Y,    & \mbox {if $X/2 \le Y <
X.$}\end{array}\right.
\]

The minimal layered-triangulation of the solid torus therefore
extends the triangulation on the boundary of a solid torus where the
edges meet the meridional disk in the triples \[ \left\{
\begin{array}{ll}
$\{X-2Y,Y,X-Y\}$, &\mbox{if $0 < Y < X/2,$ or}\\
$\{2Y-X,X-Y,Y\}$, & \mbox {if $X/2 < Y < X$.}\end{array}\right.
\]

 In
 general, we have the unique non orientable surface having boundary
 slope that of the edge $X-2Y$ for $0 < Y < X/2,$ and that of the edge
$2Y-X$ for $X/2 < Y <
 X$. Hence, the genera of the  embedded, incompressible, non orientable
surfaces in $L(X,Y)$ for these two possibilities are $e(X-2Y, Y)
+1$, for $0 < Y < X/2$, and $e(2Y-X,X-Y)+1$, for $X/2 < Y <
 X$.

But since $X$ is even, we have from Lemma \ref{e-p-q}  that $e(X,Y)
= e(X-Y,Y) = e(X-2Y,Y)+1$ ($X-Y,Y$ are both odd), $0<Y<X/2$; and
$e(X,Y) = e(X-Y,Y) = e(X-Y,(2Y-X)+(X-Y))= e(X-Y, 2Y-X)+1$ ($X-Y,Y$
are both odd), $X/2<Y<X$. Thus in both cases, we have $h =
e(X,Y)$.\end{proof}

\vspace{.125 in}\noindent{\bf Examples:}\begin{enumerate}\item The
unique incompressible, non orientable surface in $L(30,7)$ is $U_3$.

In this case we have $X = 30, Y = 7$; hence, $Y<X/2$ and e(X-2Y,Y) =
e(16,7) = 2  and h = 3. \item For $L(30,23)$, we have $X=30, $Y = 23
and now $X/2<Y$. We use $e(2Y-X,X-Y)= e(16,7)$, and, as expected, $h
= 3$.

\item $L(2,1) = \rppp$ is the only lens space in which we can embed $\rpp$.

The minimal layered-triangulation of the lens space $\rppp$ is
obtained from the two-tetrahedra extension of the
$0/1$-triangulation by ``folding over" along the edge $0$. The
$0/1$--layered-triangulation of the solid torus admits a M\"obius
band with boundary having the slope of the interior edge $2$ and a
once-punctured Klein bottle with boundary having slope of the edge
$0$ in the boundary. However, we also have the meridional disk with
boundary having the slope of the edge $0$. Folding along the edge
$0$, the boundary of the M\"obius band has intersection numbers
$3,1,2$,  the boundary of the once-punctured Klein bottle has
intersection numbers $1,1,0$, and the boundary of the meridional
disk also has intersection numbers $1,1,0$. From the meridional disk
we get an embedded $\rpp$ and from the the once-punctured Klein
bottle we get a genus $3$ nonorientable surface, which is the
embedded $\rpp$ with a trivial handle.

It follows that $\rppp$ admits an embedded nonorientable surface of
every \emph{odd} genus.
\item Those lens spaces $L(X,Y)$ that admit an embedding of the
Klein bottle have the form
 $L(4n,2n-1),
n\ge 1$.

In this case, the minimal layered-triangulation of the solid torus
that gives the minimal layered-triangulation of $L(X,Y)$ has $e(p,q)
= 1$; hence, the even edge in the boundary, say $p$, must be $2$. It
follows that $X = 2+2q$, $Y = q\ge 1$ and odd. Set $q = 2n-1$. In
particular, for $n = 1$, we have $L(4,1)$; for $n = 2$, we have
$L(8,3);\ldots$; for $n= 11$, we have $L(44,21); \ldots$.
\item There are no \emph{incompressible} non orientable surfaces
embedded in $S^2\times S^1$; however, there is an embedding of a
nonorientable surface of every even genus in $S^2\times S^1$.

 If we add a handle from one side of the non separating
$2$--sphere in $S^2\times S^1$ to its other side, then we have an
embedded Klein bottle; so, any non orientable surface with even
Euler number can be embedded in $S^2\times S^1$. On the other hand,
the minimal two-tetrahedra, layered-triangulation of $S^2\times S^1$
is obtained by folding the boundary of the two-tetrahedra,
layered-triangulation of the solid torus extending $\{2,1,1\}$ along
the edge labeled $2$.

The normal surfaces in this two-tetrahedra triangulation of the
solid torus are given in the proof of Theorem \ref{normal}. The only
possibilities leading to a nonorientable normal surface in
$S^2\times S^1$ is the M\"obius band with boundary having slope that
of the boundary edge $2$. Upon folding over the edge $2$, this
surface becomes a Klein bottle. There are no incompressible non
orientable surfaces. It is curious that this minimal genus non
orientable surface is captured by the minimal layered-triangulation.

\end{enumerate}

\subsection{$0$-- and $1$--efficient  layered-triangulations of lens spaces}

Layered-triang-ulations of lens spaces are quite simple and do not
exhibit the full richness of $0$-- and $1$--efficient
triangulations; however, they do provide examples of triangulations
that are not efficient as well as allow a complete characterization
of those layered-triangulations that are efficient.

A triangulation of a closed $3$--manifold  is said to be {\it
$0$--efficient} if and only if the only normal $2$--spheres are
vertex-linking \cite{jac-rub0}. A layered-triangulation of a lens
space is said to be {\it $1$--efficient} if and only if it is
$0$--efficient and the only normal tori are thin edge-linking. There
is a definition of $1$--efficient that applies to general
triangulations of closed $3$--manifolds; however, for
layered-triangulations of lens spaces, what we have given is an
equivalent definition.

\vspace{.15 in} \noindent{\it ${0}$--efficient
layered-triangulations of lens spaces.} A layered-triangulation of a
lens space can be viewed in two ways as a layered-triangulation of a
solid torus with its boundary folded over; one view is the reverse
of the other. In other words, once we fold over, we can start from
either end and we have a layered solid torus with its boundary
folded over. If in either of these views we have the meridian slope
of the solid torus the slope of an edge in the triangulation, we say
``a meridional slope is the slope of an edge". If a meridional slope
is the slope of an edge, there is a level in the triangulation where
the meridional slope of a solid torus  is the slope of an edge in
the boundary of that solid torus. It follows that the
layered-triangulation of the lens space can be decomposed as the
union of two subcomplexes, each a layered-triangulation of a solid
torus where one is an extension of the $0/1$--triangulation on the
boundary of the solid torus. It is possible in such a case that we
have the creased $3$--cell. We have the following characterization
of $0$--efficient layered-triangulations of lens spaces.

\begin{thm} \label{0-eff-lens} A layered triangulation of a lens space,
distinct from $ S^3$ and $S^2\times S^1$,  is $0$--efficient if and
only if no edge has the slope of a meridian.
\end{thm}
\begin{proof}
Suppose an  edge has the slope of a meridian. Then split the
layered-triangulation of the lens space at a level where the
meridional slope is the slope of an edge in the boundary of the
solid torus; i.e., extends the $0/1$--triangulation. One possibility
is that we are  folding the boundary of a layered-triangulation that
extends the $0/1$--triangulation; hence, we get $S^3$ or we are
folding along the edge $0$ and we get $\rppp$. Since we are assuming
we do not have $S^3$, we must get $\rppp$ and the meridional disk
identifies to a projective plane and its double is a non
vertex-linking $2$--sphere. So, the triangulation is not
$0$--efficient. The other possibility is that we have two solid tori
and in one solid torus, the meridional curve has the slope of an
edge in the boundary. In the other solid torus there are two
possibilities: the edge in the boundary of this solid torus, which
is identified to the edge having the slope of the meridian,  is
either thick or it is thin. If it is thick, then we must have $S^3$.
Again this is ruled out by hypothesis, the edge in this solid torus,
identified to the edge having slope that of a meridian, must be
thin. Hence, there is a thin edge-linking annulus about the edge.
This annulus and two copies of the meridional disk form a non
vertex-linking normal $2$--sphere. Hence, if we do not have a
layered-triangulation of $S^3$, then the existence of an edge having
the slope of a meridian gives a non vertex-linking normal
$2$--sphere and the triangulation is not $0$--efficient.

Now, suppose a layered-triangulation of a lens space is not
$0$--efficient. Hence, there must be a non vertex-linking, normal
$2$--sphere, say $\Sigma$. The layered-triangula-tion of a lens
space is obtained by folding over a layered-triangulation of a solid
torus and $\Sigma$ meets this solid torus in a planar surface (the
$2$-sphere $\Sigma$ is obtained from this  planar surface after
identifications along its boundary). Again, the possible
intersection numbers for such a planar normal surface with the edges
$e_1, e_2, e_3$ in the boundary of the solid torus, folding over
edge $e_3$, are: $2,2,2; 1,1,2;$ or $1,1,0$, respectively.

Now, by Theorem \ref{layered-0-1}, the possible  planar (genus $0$)
normal surfaces embedded in a layered-triangulation of a solid torus
are (using the labeling from Theorem \ref{layered-0-1})
\begin{itemize}
\item [(1)] the vertex-linking disk with boundary intersection $2,2,2$, which
gives a vertex-linking $2$--sphere in the lens space; \item [(2)]
the meridional disk either with boundary intersection $1,1,2$, which
then gives $ S^2\times S^1$ and contradicts our hypothesis, or with
boundary intersection $1,1,0$ in which case we have a meridional
slope the slope of an edge and we have $\rppp$;
\item[(3)] an edge-linking annulus with boundary slope the slope of an
edge. If the intersection numbers are $1,1,2$, then we are folding
along the univalent edge and thus by reversing our layered solid
torus, we see that we have a creased $3$--cell and an edge with
meridional slope; this includes the possibility of a fat
edge-linking annulus, where we have two edges with the same slope.
If the intersection numbers are $1,1,0$, then we get a torus upon
identification of the boundary and not a $2$--sphere;\item[(4)] a
trivial, non vertex-linking disk having boundary with intersection
numbers $2,2,2$, in which case an edge has the slope of the
meridional disk;
\item [(5)] a non edge-linking annulus with boundary slope the slope
of an edge; so, again, an edge has the slope of a meridional
disk.\end{itemize}

Hence, if the lens space is distinct from $ S^3$ or $S^2\times S^1$
and no edge in the layered-triangulation has slope that of a
meridian, then the layered-triangulation is
$0$--efficient.\end{proof}

 \vspace{.125 in}\noindent{\bf Examples.}\begin{enumerate}
\item There are infinitely many distinct $0$--efficient layered
triangulations of $S^3$, each with an edge having the slope of a
meridian.

In the previous proof we noted that if a lens space is formed by
identifying the boundaries of two layered-triangulations of the
solid torus via a simplicial identification with an edge having the
slope of the meridian in one layered solid torus being identified to
a thick edge in the other layered solid torus, then we get a
layered-triangulation of $S^3$. This is used to construct infinitely
many distinct $0$--efficient, layered-triangulations of $S^3$, each
with an edge having the slope of the meridian.

Consider the layered-triangulations of $S^3$ given by folding over
along the edge $2n+1$ in the minimal
$(n/n+1)$--layered-triangulation of the solid torus. We write this
as $\{n,n+1,2n+1\}\lra \{1,1,2\}$,
 $n\ge 0$. Since we have folded over along the univalent edge,
each of these triangulations is the same as attaching the creased
$3$--cell to the minimal layered-triangulation of the solid torus
extending the $1/n$--triangulation; namely,
\begin{tabbing}\=$\{n,n+1,\underline{2n+1}\}\lra
\{1,1,\underline{2}\}$\\\>$\{n, n+1, 1\}\lra
\{1,1,0\}$.\end{tabbing} So, in each triangulation, we have an edge
having the slope of a meridian being attached to the thick edge.

Now, to see that these are $0$--efficient triangulations, we analyze
the possible normal surfaces in these layered-triangulations.  A
normal surface in such a triangulation of $S^3$ is determined by a
normal surface in the minimal $(n/n+1)$--layered-triangulation of
the solid torus, along with some identifications in its boundary.
So, we can proceed, as above, and consider those normal surfaces in
the minimal $(n/n+1)$--layered-triangulation of the solid torus,
where the surface has boundary intersection numbers with the edges
labeled $n$, $n+1$ and $2n+1$ one of: $2,2,2$; $1,1,2$; or $1,1,0$,
respectively. The only interesting set of intersection numbers here
is $1,1,2$. Typically there would be an edge-linking annulus in the
solid torus with boundary slope having the intersection numbers
$1,1,2$; however, in these examples, the intersection numbers
$1,1,2$ correspond to the slope of the thick edge and as such there
is no associated edge-linking annulus. We do, however, have in each
an almost normal annulus with boundary slope having the intersection
numbers $1,1,2$ and from these we get almost normal $2$--spheres in
each of these triangulations of $S^3$.

Notice for $n = 0$, we have one of the two, two-tetrahedra
$0$--efficient layered triangulations of $S^3$, which is two creased
$3$--cells identified along their boundaries ({\it not} doubled);
and for $n = 1$, we have the one-tetrahedron, one-vertex minimal
triangulation of $S^3$. For these examples, we  fold along the
univalent edge. We also can  get a two-tetrahedra, layered
$0$--efficient triangulation of $S^3$ by folding along the thick
edge in the minimal $0/1$--layered-triangulation of the solid torus.

\item All lens spaces (except
$\rppp$ and $S^2\times S^1$) admit an infinite family of distinct
$0$--efficient layered-triangulations; furthermore, by Theorem
\ref{0-eff-lens}, except for $S^3$, none of these
layered-triangulations have an edge with the slope of a meridian.

Recall that a path in the $L$--graph starting at $p/q$ and ending at
$1/1$ determines a layered-triangulation of the solid torus and if
the path  runs through $0/1$, then we have an edge with the slope of
a meridian. To obtain a layered-triangulation of a lens space we
fold along an edge in the boundary of the
$p/q$--layered-triangulation of the solid torus. Now, the reverse
layering gives us another path in the $L$--graph. Does it run
through $0/1$? If we have $\{p,q,p+q\}\lra \{1,1,2\}$ ( we fold over
along the univalent edge, except in the case
$\{1,1,2\}\lra\{1,1,2\}$ where $2$ is not univalent and we get
$S^2\times S^1$), then the reverse path has
\begin{tabbing}\=$\{p,q,\underline{p+q}\}\lra
\{1,1,\underline{2}\}$\\\>$\{p,q, \abs{p-q}\}\lra
\{1,1,0\}$,\end{tabbing} which gives an edge having the slope of a
meridian. So, we can not fold along the univalent edge. If we fold
along the edge $p$ or $q$, say $q$, then the reverse path passes
through $1/1$ every time our original path  passes through $p/q$.
Hence, the reverse path has
\begin{tabbing}\=$\{p,\underline{q},p+q\}\lra
\{1,\underline{2},1\}$\\\>$\{p,2p+q, p+q\}\lra \{1,0,1\}$,
\end{tabbing} whenever we go from $p/q$ in the $L$--graph to
$p/p+q$ in the $L$--graph in the original path.

We conclude,  with the exception of folding along an edge in the
boundary of a $0/1$-layered-triangulation or a
$1/1$--layered-triangulation,  for any point $p/q$ in the $L$--graph
and any path from $p/q$ to $1/1$ that does not pass through $0/1$ or
$p/p+q$ we have a layered-triangulation of the solid torus so that
folding along the edge $q$ gives a layered-triangulation of
$L(2p+q,p)$ with no edge having slope a meridian; hence, we have a
$0$--efficient layered-triangulation by Theorem \ref{0-eff-lens}.

Recall for paths beginning at $0/1$, if we fold, we get $S^3$ or
$\rppp$; for paths beginning at $1/1$, we get $S^2\times S^1$ if we
fold over the edge $2$ and get $L(3,1)$, otherwise. If we fold over
the univalent edge  (thin edge) labeled $1$ we do not get a
$0$--efficient triangulation; whereas, if we fold over the thick
edge labeled $1$ we do get a $0$--efficient triangulation. Hence,
paths in the $L$--graph beginning at $1/1$ and never layering along
the edge paired with $2$ when passing through $1/1$ provide an
infinite family of $0$--efficient triangulations of
$L(3,1)$.\end{enumerate}

\begin{cor} For any lens space other than $S^2\times S^1$ and $\rppp$, there exists a
$0$--efficient layered-triangulation having an arbitrarily large
number of tetrahedra. Hence, for each such lens space there are
infinitely many distinct $0$--efficient layered-triangulations.
\end{cor}

\vspace{.15 in}\noindent{\it $1$--efficient layered-triangulations
of lens spaces.}  We have the following characterization of
$1$--efficient layered-triangulations of lens spaces. Just as in the
case of a layered-triangulation of the solid torus, we say a
layered-triangulation of a lens space is {\it nearly-minimal} if no
edge has the slope of a meridian and no two \emph{thin} edges have
the same slope. Recall that in a layered-triangulation of a lens
space there are two thick edges; whereas, only one in a
layered-triangulation of a solid torus. For lens spaces the
two-tetrahedron layered-triangulation of $\rppp$ is minimal but is
not nearly-minimal and the two-tetrahedron layered triangulation of
$L(3,1)$ obtained by folding over along the univalent edge in the
two-tetrahedron layered-triangulation of the solid torus extending
the $1/1$--triangulation is minimal but is not nearly-minimal.
Otherwise, a minimal triangulation is nearly-minimal.

\begin{thm}\label{1-eff-lens} A layered-triangulation of a lens space, distinct from
$S^3$ and $S^2\times S^1$, is $1$--efficient if and only if it is
nearly-minimal.
\end{thm}
\begin{proof} Suppose the layered-triangulation is obtained by
folding along an edge in the boundary of a layered-triangulation of
the solid torus $\bbb{T}$ and it is not nearly-minimal. Then either
an edge has the slope of a meridian or two thin edges have the same
slope.

If an edge has the slope of a meridian, then by Theorem
\ref{0-eff-lens} and our assumption that the lens space is not
$S^3$, the triangulation is not $0$--efficient and so can not be
$1$--efficient.

Thus we may assume, no edge has the slope of a meridian; hence,
there are two thin edges (in the lens space) having the same slope.
It follows that two thin edges, say $e$ and $e'$, in $\bbb{T}$ have
the same slope. Now, by Theorem \ref{layered-0-1}, there is either a
vertex-linking disk with a non thin edge-linking tube in $\bbb{T}$
or an edge-linking annulus about both $e$ and $e'$ and having
boundary slope the slope of these two edges. In the first case, we
get a non thin edge-linking torus in the lens space and the
triangulation of the lens space is not $1$--efficient. In the second
case, we have two possibilities for the intersection numbers for the
slope of $e$ and $e'$: either  $1,1,2$  or $1,1,0$. Now, if the
slope of $e$ and $e'$ has intersection number $1,1,2$, then in the
reverse layering there is an edge having the same slope as $e$ and
$e'$ (we could actually choose one of these to be such an edge) and
having the slope of a meridian. An easy way to see this is that the
slope with intersection numbers $1,1,2$ bounds the meridional disk
in the degenerate layered-triangulation of the solid torus, the
one-triangle M\"obius band. But we have assumed no edge has the
slope of a meridian. If the intersection numbers for the boundary
slope are $1,1,0$ and since we have that both $e$ and $e'$ are both
thin edges, then we have a non thin edge-linking torus in the lens
space.

Hence, under the assumption of an edge having the slope of a
meridian or two thin edges having the same slope, the
layered-triangulation of the lens space is not $1$--efficient. It
follows that our condition is necessary.

Now, suppose the layered-triangulation of the lens space is
\emph{not} $1$--efficient. By Theorem \ref{0-eff-lens}, we may
assume it is $0$--efficient. Thus we must have a normal torus, say
$T$, that is not thin edge-linking in the lens space. We again
consider the layered-triangulation of the lens space as coming from
a layered-triangulation of a solid torus $\bbb{T}$ with its boundary
folded over. The torus $T$ meets $\bbb{T}$ in a genus $0$ or a genus
$1$ normal surface. By Theorem \ref{layered-0-1}, the possibilities
that can lead to a normal torus are (we use the numbering in the
conclusions of Theorem \ref{layered-0-1}):
\begin{enumerate}\item[(3)] An edge-linking annulus with boundary slope
the slope of an edge. The  possibilities are: the thin edge-linking
annulus whose boundary slope has  intersection numbers $1,1,0$ and
is about the edge along which we are folding; or a thin edge-linking
annulus whose boundary slope has  intersection numbers $1,1,2$ and
meets the edge along which we are folding in $2$ points. The former
gives a thin edge-linking torus in the lens space, which is not $T$;
the latter, as above, can only happen if the reverse layering has an
edge having the slope a meridian, which contradicts that the
layered-triangulation of the lens space is $0$--efficient and the
lens space is not $S^3$.\item[(4)] We  have (4B), a non edge-linking
annulus with boundary slope the slope of an edge. But we only have
such an annulus if there is an edge having the slope of a
meridian.\item[(5)] A fat edge-linking annulus. We only have such an
annulus if there are two edges with the same slope.
\item[(6)] In the case of (6A), a vertex-linking disk with a thin
edge-linking tube, we have a thin edge-linking torus and not $T$. In
the case of (6B), an edge-linking annulus with a thin edge-linking
tube, we have the possibilities of intersection numbers $1,1,0$ and
$1,1,2$. In the former case, after folding we get a genus $2$
surface and not $T$; in the latter, we would have to have an edge
with the slope of a meridian, which we have ruled out.
\end{enumerate}

In the cases of (7), (8) and (9) of Theorem \ref{layered-0-1}, we
have an edge with the slope of a meridian, two edges with the same
slope or both possibilities.

Hence, if the lens space is distinct from $ S^3$ or $S^2\times S^1$
and  the layered-triangulation is nearly-minimal, then the
layered-triangulation is $1$--efficient.\end{proof}

Note there are triangulations of $S^3$ that are $1$--efficient and
not nearly-minimal; also, the minimal two-tetrahedra triangulation
of $S^2\times S^1$ is nearly-minimal but is not $1$--efficient.

\begin{figure}[h]
\psfrag{e}{$e$}\psfrag{4}{\small{$1'$}}\psfrag{1}{\small{$1$}}\psfrag{2}{\small{$2$}}
\psfrag{3}{\small{$3$}}\psfrag{A}{\Large (A)} \psfrag{B}{\Large (B)}
\psfrag{a}{\begin{tabular}{c}fold\\
           over\\ \end{tabular}}\psfrag{b}{\begin{tabular}{c}identify\\
           \end{tabular}}\psfrag{c}{\begin{tabular}{c}two-tetrahedra\\extension
           of\\$1/1$--triangulation\\
           \end{tabular}}
 {\epsfxsize = 3 in
\centerline{\epsfbox{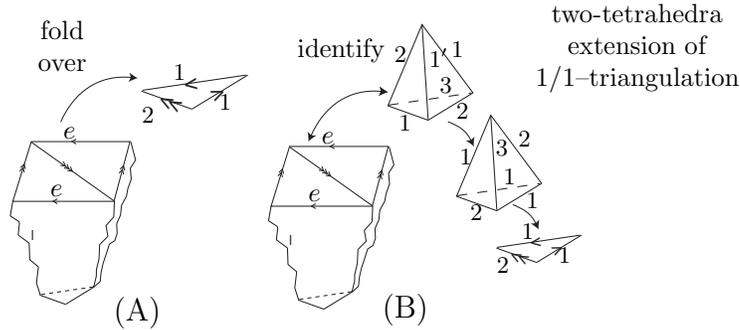}} } \caption{Examples of
opening a layered-triangulation (A)  at a thick edge  and (B) at
$1/1$.} \label{f-open-edge}
\end{figure}

\vspace{.125 in}\noindent{\bf Examples.}\begin{enumerate}
\item There are infinitely many distinct $1$--efficient
triangulations of $S^3$. In fact, the triangulations,
$\{n,n+1,2n+1\}\leftrightarrow \{1,1,2\}$,  given in the above
example of infinitely many $0$--efficient triangulations for $S^3$
also are  $1$--efficient.
\item An opening at $1/1$ of a nearly-minimal layered-triangulation
of a lens space gives a nearly-minimal layered-triangulation; hence,
we can generate from a minimal layered-triangulation of a lens
space, that is distinct from the two-tetrahedra minimal
layered-triangulation for $\rppp$, three additional nearly-minimal
layered-triangulations (two by opening at one end, or the opposite
end, and a third by opening at both ends).
\item An opening at $1/1$ is a special case of {\it opening at a
thick edge}; that is, if we have the minimal $1/k$--layered
triangulation, $k\ge 2$, and layer on $1$, we get the minimal
$(k/k+1)$--layered-triangulation. Then we can undo this layering by
layering on the univalent edge $2k+1$ to return to a
$1/k$--layered-triangulation (not minimal). We still have an edge
that meets the meridional disk once but now it is not a thick edge.
The resulting $1/k$--layered-triangulation is nearly-minimal. It
follows that a minimal $1/n$--layered-triangulation, $n\ge 2$, of
the solid torus has $n-1$ such openings, one for each $k, 2\le k
<n$,  giving $n-1$ distinct nearly-minimal
$1/n$--layered-triangulations. See Figure \ref{f-open-edge}.
\end{enumerate}

 \begin{cor} A lens space distinct from $S^3$ has only finitely many nearly-minimal
 (hence, $1$--efficient) layered-triangulations.\end{cor}

\begin{proof} Suppose $\T$ is a layered-triangulation of a solid torus
and $e$ is an edge of $\T$ in the boundary. Let $\T'=\T\cup_e\Delta$
be a layering of $\T$ along the edge $e$. Now, let $e'$ be the
univalent edge in $\T'$ and let $\T'' = \T'\cup_{e'}\Delta'$ be a
layering of $\T'$ along the univalent edge $e'$ of $\T'$. Finally
let $e''$ be the univalent edge of $\T''$. Then $e''$ has the same
slope as $e$. Hence, any time we layer on the univalent edge we get
two edges with the same slope.

Now, using this notation, if the edge $e$ is a thick edge, then we
do not have two thin edges with the same slope; however, this, which
we call opening at a thick edge, can only happen once. See Example 3
in the previous paragraph; this is precisely what is happening in
that example.

When $e$ is not a thick edge, this leads to a layered-triangulation
of the solid torus that is not a nearly-minimal
layered-triangulation. If such a triangulation of the solid torus is
folded over to form a lens space, then the layered-triangulation of
the lens space is typically not nearly-minimal. However, it is
possible that we fold over the boundary of $\T''$ and $e''$ becomes
a thick edge; i.e., we have \emph{not} folded over along $e''$.

If this is the case,  the two tetrahedra $\Delta\cup_{e'}\Delta'$,
 become a two-tetrahedra extension of
the $1/1$-triangulation after folding along an edge in the boundary
of $\Delta'$. Thus the layered-triangulation of the lens space is
just an opening at $1/1$ of the layered-triangulation of the {\it
same} lens space had we folded over the boundary of the
layered-triangulation $\T$. In this case, we have the situation
given in Example 2 in the paragraphs preceding this Corollary.

Finally, if we fold over along $e''$, which is the univalent edge,
then we have an edge with slope the meridian and the triangulation
can not be nearly-minimal. It follows there are only finitely many
ways to make the triangulation of a fixed lens space a
nearly-minimally layered-triangulation.
\end{proof}

Note that the typical situation would be at most five $1$--efficient
triangulations, the minimal one, opening along $1/1$ at either end,
which typically gives two more,  opening at both ends and, finally,
opening at the thick edge. There also are those paths in the
$L$--graph that run through a number of vertices of the form $1/k$.
For example, layered-triangulations of the lens space $L(n,1), n\ge
4$ are formed from a $(1/n-2)$--layered-triangulation by folding
along the edge labeled $n-2$ in the boundary. We saw above in
Example 3 that each of these have $n-3$ places where the
triangulation can be opened at a thick edge. If we fold the
two-tetrahedra $1/1$--layered-triangulation along the thick edge, we
get $L(3,1)$, which we can then open along the thick edge to get a
nearly minimal four-tetrahedron triangulation of $L(3,1)$.

\subsection{Heegaard splittings of lens spaces}

We are able to use layered-triangula-tions for a new proof of the
uniqueness of Heegaard splittings for $S^3$ and $S^2\times S^1$,
first done by F. Waldhausen \cite{wald-HS3}, and for other lens
spaces, first done by F. Bonohan and J.P. Otal \cite{bon-otal}.
While the combinatorics of this proof make it seem quite elementary,
we have lurking in the background, the observation, first made by
Rubinstein \cite{rubin-jac}, that for any triangulation of a closed
$3$--manifold, a strongly irreducible Heegaard surface is isotopic
to an almost normal surface; and the observation by A. Casson and C.
McA. Gordon \cite{cass-gord} that for non Haken manifolds, an
irreducible Heegaard surface is strongly irreducible. We comment
here that the methods initiated by Rubinstein and presented in
\cite{rubin-jac} do \emph{not} use the classification of Heegaard
splittings of $S^3$, as do those methods coming from thin position.
One of our reasons for considering this alternate proof of
uniqueness of Heegaard splittings for lens spaces is a belief that
one can use nice triangulations to gain important
topological/geometric information about three-manifolds as exhibited
in this proof. In particular, these methods may extend to
layered-triangulations of higher genera three-manifolds and offer
analogous results. We use simple triangulations in
\cite{jac-rub-sSFS} to study the Heegaard splittings of small
Seifert fibered spaces.

We assume familiarity with basic notions and results regarding
Heegaard splittings. We refer the reader to the exposition by M.
Scharlemann \cite{shar-HS}.

We remind the reader that, except for the lens spaces $S^3, \rppp,
S^2\times S^1$, and $L(3,1)$, the only normal surface in a minimal
layered-triangulation of a lens space is a vertex-linking
$2$--sphere (possibly) with some thin edge-linking tubes. Hence, an
almost normal tubed surface must be a vertex-linking $2$--sphere
(possibly) with some thin edge-linking tubes and an almost normal
tube; or two such surfaces with an almost normal tube between them.
There are no almost normal octagonal surfaces. We are able to
characterize which of these surfaces are Heegaard surfaces.

\begin{figure}[h]
\psfrag{a}{\begin{tabular}{c}layered\\
          solid torus\\
          \end{tabular}}\psfrag{b}{\begin{tabular}{c}from\\above\\
           \end{tabular}}\psfrag{c}{\begin{tabular}{c}from\\below\\
           \end{tabular}}
 {\epsfxsize = 2.5 in
\centerline{\epsfbox{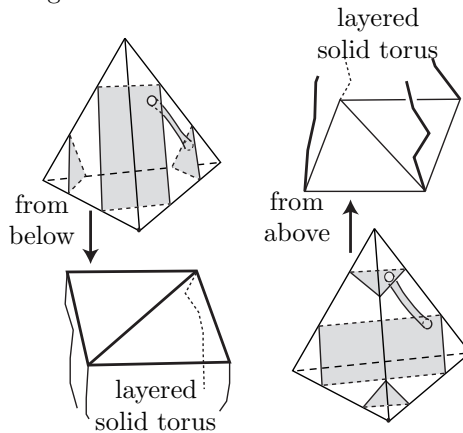}} } \caption{An almost normal
tube at the same level as a thin edge-linking tube.}
\label{f-samelevel}
\end{figure}

If we have an  almost normal tubed surface in a
layered-triangulation of a lens space, we say the {\it almost normal
tube is at the same level as a thin edge-linking tube} if the almost
normal tube has at least one end in a quadrilateral of the thin
edge-linking tube. See Figure \ref{f-samelevel}.

The following lemma and its corollary are used in our understanding
of Heegaard splittings of lens spaces and also used later when we
consider Heegaard splittings of Dehn-fillings. Note the following
lemma is not necessarily true if we have an almost normal tubed
surface or an almost normal octagonal surface without thin
edge-linking tubes; also, in the case of the two-tetrahedron
$0/1$--layered-triangulation, we have an edge-linking annulus that
is compressible.

\begin{lem}\label{no-comp} Suppose $\T$ is a non-degenerate minimal
layered-triangulation of the solid torus $\bbb{T}$ distinct from the
$0/1$--layered-triangulation. Furthermore, suppose $A$ is an
edge-linking annulus possibly with thin edge-linking tubes or an
almost normal octagonal annulus with thin edge-linking tubes. Let
$B'$ denote the complementary annulus to $A$. Then $A$ does not
$\bdy$--compress into $B'$.\end{lem}

\begin{proof} The embedded tori which correspond to a boundary
torus at each level of the layering are called {\it level tori}.  If
we have a thin edge-linking tube of $A$, then in each level torus
that meets the tube, we have a companion and complementary annulus.
If we denote the component of the complement of $A$ that meets $B'$
by $H'$, then the complementary annuli in the various level-tori
that meet a thin edge-linking tube of $A$ determine a family of
parallel, properly embedded, annuli in $H'$.

We order the thin edge-linking tubes in $A$ from top to bottom,
$tube_1, tube_2,\ldots,\\tube_k$. If we split the layered
triangulation along a level torus meeting the $i^{th}$ tube of $A$,
then we split off an initial segment of the layered solid torus
$\bbb{T}$, say $\bbb{T}_i$, and split $A$ into an edge-linking
annulus $A_i$ with tubes embedded as a normal surface in
$\bbb{T}_i$. For simplicity, we shall always assume we have done
this at the lowest level possible for $tube_i$.

If $A$ has no tubes, then a $\bdy$-compression of $A$ into $B'$
gives that $A$ must be longitudinal. But in a minimal
layered-triangulation of a solid torus, there are no longitudinal
normal annuli (only almost normal octagonal annuli). But we have
assumed that if $A$ is longitudinal, then it does have tubes. Thus,
we may suppose $A$ has tubes.

If $A$ $\bdy$--compresses into $B'$, denote the compressing disk by
$D_0'$. Two copies of $D_0'$ along with a disk in $B'$ determined a
compressing disk $D_0$ for $A$ and we may assume that $D_0$ does not
meet any of the complementary annuli in the various level-tori that
are parallel in $H'$ to $B'$, including, of course, $B'$.

We have that $A$ has tubes. Consider the way $D_0$ meets the
collection of complementary annuli to the tubes of $A$. We may
assume, up to isotopy, $D_0$ has minimal intersection with these
complementary annuli.

Suppose first that it may be the case that $D_0$ does not meet any
of the complementary annuli. Then $D_0$ is a compressing disk in a
region between two level tori in the layering and $D_0$ meets each
of these level tori in annuli that correspond to complementary
annuli for some edge slope. Thus we must have the edge slope of $A$
the same as the slope of the edge that $tube_1$ links. But two edges
can not have the same slope in a minimal layered-triangulation of
the torus. So, if $D_0$ exists, it must meet some of the
complementary annuli associated with the thin edge-linking tubes of
$A$.

If we have a simple closed curve component in the intersection, then
we have an innermost, on $D_0$, simple closed curve, say in the
intersection with $B_i'$. If this curve were essential in $B_i'$,
this would mean we have an edge-slope the slope of a meridian, which
can not happen in a minimal layered-triangulation distinct from the
two-tetrahedra $0/1$--layered-triangulation of the solid torus.
Hence, such a component of intersection must be trivial in $B_i'$.
But now, using standard techniques, we could modify the intersection
to reduce the number of components; this is a contradiction to our
having minimal intersection. If we have a component of intersection
that is an outermost, on $D_0$, spanning arc and if this is a
component common with, say $B_i'$, then we have a boundary
compression of $A_i$ into $B_i'$, or the arc of intersection is
inessential in $B_i'$. The latter would contradict that the
intersection of $D_0$ with the collection of complementary annuli is
minimal; hence, we must have a boundary compression of $A_i$ into
$B_i'$.

We now let $A_i$ and $B_i'$ play the role of $A$ and $B'$ in the
above argument and let $D_i'$ denote the outermost disk in $D_0$
that gives the boundary compression of $A_i$ into $B_i'$. Again, we
have a compressing disk, say $D_i$, that plays the role of the disk
$D_0$ above. Continuing in this way, we finally have as the only
possibility, that for some $A_j$, we have a compressing disk that
meets no complementary annuli. But as above, this means that two
edges must have the same slope. This contradiction completes the
argument.\end{proof}

We are now ready to characterize those normal or almost normal
surfaces in a minimal layered-triangulation of a lens space that are
Heegaard surfaces. First, we consider lens spaces distinct from
$S^3, \rppp, S^2\times S^1$, and $L(3,1)$; following this we will
consider these special cases.

\begin{thm}\label{samelevel} Suppose $L$ is a lens space distinct from $S^3, \rppp,
S^2\times S^1$, and $L(3,1)$. A normal or almost normal surface $S$
in a minimal layered-triangulation of $L$ is a Heegaard surface if
and only if $S$ has an almost normal tube along a thick edge or $S$
has an almost normal tube at the same level as a thin edge-linking
tube. In particular, there are no normal  or almost normal octagonal
Heegaard surfaces in a minimal layered-triangulation of
$L$.\end{thm}

\begin{proof} We begin by showing that either of these conditions is
sufficient for an almost normal surface to be a Heegaard surface.

Suppose $S$ is an almost normal surface with an almost normal tube
along a thick edge in a minimal layered-triangulation of $L$. Since
we have excluded $S^3, \rppp, S^2\times S^1$, and $L(3,1)$, the
layered-triangulation of $L$ has no creased $3$--cells and no
layering along a univalent edge. Hence, the layered-triangulation of
$L$ is formed by folding along an edge in the boundary of a minimal
(non-degenerate) layered-triangulation of a solid torus. We may
assume we have opened the layered-triangulation of the solid torus
so that the thick edge along which we have an almost normal tube is
in the one-tetrahedron solid torus at the bottom of the
layered-triangulation of the solid torus. If $S$ has more than one
tube, then the thin edge-linking tubes come in an order in the
layered solid torus, ordering from the bottom toward the top, where
we have opened the lens space triangulation. Beginning with the
first thin edge-linking tube from the bottom, there is a pinched
disk embedded in the layered solid torus between the edge the tube
is linking and a curve in the boundary of the one-tetrahedron solid
torus at the bottom of the solid torus. Thus this tube unwinds about
the thick edge and corresponds to a trivial handle attached to $S$.
We can remove this handle to get an almost normal surface $S_1$ with
an almost normal tube along a thick edge. In this way we see that
$S$ is a stabilization of the standard genus one Heegaard surface of
$L$ and therefore is a Heegaard surface.

Now, suppose $S$ is an almost normal surface with an almost normal
tube at the same level as a thin edge-linking tube. Then there is a
tetrahedron in our layering as shown in Figure \ref{f-samelevel}.
After a small isotopy of the almost normal tube, it and the thin
edge-linking tube provide a spine in the torus where they occur at
the same level. Hence, if there are any other tubes (necessarily
thin edge-linking), they come in an order and thus any one adjacent
to the level of the spine would unwind about it and form a trivial
handle. It follows that $S$ is obtained by adding trivial handles to
a surface $S'$, which is the boundary of a small neighborhood of the
spine of a level torus; so, $S'$ bounds a genus-$2$-handlebody, say
$H'$, on one side. However, the other side of $S'$ is two solid tori
attached by a $1$-handle determined by the complementary disk of
this spine on the level torus and therefore $S'$ is a genus $2$
Heegaard surface for $L$. The surface $S$ is a stabilization of $S'$
and thus $S$ is a Heegaard surface for $L$.

We have shown that an almost normal surface that has either an
almost normal tube along a thick edge or an almost normal tube at
the same level as a thin edge-linking tube, in a minimal
layered-triangulation of a lens space distinct from $S^3, \rppp,
S^2\times S^1$, and $L(3,1)$,  is sufficient for it to be a Heegaard
surface.

Conversely, suppose we have a normal or almost normal Heegaard
surface $S$ splitting $L$ into the handlebodies $H$ and $H'$. Then
$S$ is the vertex-linking $2$--sphere with tubes, possibly an almost
normal tube, or two vertex-linking $2$--spheres, one necessarily
with thin edge-linking tubes, that are connected by an almost normal
tube. In particular, $S$ can not be an almost normal octagonal
surface.

We first consider the case that we have a vertex-linking $2$--sphere
with tubes. If $S$ has only one tube, then for $S$ to be a Heegaard
surface, $S$ must be an almost normal surface with the almost normal
tube along a thick edge. So, we may assume $S$ has more than one
tube and thus must have a thin edge-linking tube. Since there can be
only one almost normal tube, there is a direction in the layering so
that the last tube is a thin edge-linking tube. If we split the
layering and the thin edge-linking tube at this level, we have two
solid tori, each with a layered-triangulation, $\bbb{T}$ and
$\bbb{T}'$, and the surface $S$ is split into a thin edge-linking
annulus with tubes, $A\subset\bbb{T}$, and a thin edge-linking
annulus, $A'\subset\bbb{T}'$. It is possible that $\bbb{T}'$ is
degenerate (a M\"obius band) and $A'$ is also degenerate (a simple
closed curve).

Recall that if $A$ has an almost normal tube along an edge and if it
does not have the almost normal tube along the thick edge, or at the
same level as a thin edge-linking tube, or at the same level as the
edge-linking annulus associated with $A$, then $A$ is isotopic,
keeping the boundary of $\bbb{T}$ fixed, to a \emph{normal} surface
that is a thin edge-linking annulus with all tubes thin
edge-linking. Thus if $A$ does not have an almost normal tube along
the thick edge or does not have an almost normal tube at the same
level as a thin edge-linking tube, including the level of the
edge-linking annulus associated with $A$, then we may assume $A$ is
normal.

On the other hand, if we let $B'$ denote the complementary annulus
to $A$ and have chosen notation so $B'\subset H'$, then $B'$ must
compress or $\bdy$--compress in $H'$. $B'$ can not compress, since
no edge has the slope of a meridian. Also, $B'$ can not
$\bdy$--compress into $\bbb{T}'$, since $A'$ can not be a
longitudinal. It follows that the only possibility is that $B'$
$\bdy$--compresses into $\bbb{T}$. But this is the same as $A$
$\bdy$--compressing into $\bbb{T}$. However, by Lemma
\ref{samelevel}, this is impossible for $A$ normal.

In particular, a normal surface can not be a Heegaard surface.

If $S$ were two vertex-linking $2$--spheres connected by an almost
normal tube, then we conclude that $S$ can not be a Heegaard
surface. For in this situation, if we compress the almost normal
tube, then we have two vertex-linking $2$--spheres, one necessarily
with tubes, say $S_1$ and $S_2$, where we may assume $S_1$ separates
$S_2$ from the vertex (they are nested).  If $S$ separates $L$ into
the handlebodies $H$ and $H'$, then notation may be chosen so that
$S_2$ bounds a handlebody that is a disk connected summand of $H'$
and $S_1$ is disjoint from this summand of $H'$. However, if this
were the case, then $S_2$ would necessarily bound a handlebody on
both sides and be, itself, a Heegaard surface. But the earlier
argument concludes that a normal surface can not be a Heegaard
surface.

We conclude that a Heegaard surface, in a lens space distinct from
$S^3, \rppp, S^2\times S^1$, and $L(3,1)$, is a vertex-linking
$2$--sphere with an almost normal tube along a thick edge or an
almost normal tube at the same level as a thin edge-linking tube.
\end{proof}

We now consider those cases excluded in Theorem \ref{samelevel};
i.e., the lens space $L$ is one of $S^3, \rppp, S^2\times S^1$, or
$L(3,1)$.

\vspace{.15 in}\noindent For $S^3$:

By Theorem \ref{surface-in-lens} a normal or almost normal surface
in the one-tetrahedron layered-triangulation of $S^3$ is: if normal,
then it is the vertex-linking $2$--sphere or the vertex-linking
$2$--sphere with a thin edge-linking tube about the edge $3$ (the
trefoil knot); if almost normal, then it is the vertex-linking
$2$--sphere with an almost normal tube along the thick edge, or the
vertex-linking $2$--sphere with a thin edge-linking tube along the
edge $3$ and an almost normal tube along the thick edge, or the
 almost normal octagonal $2$--sphere. There also is the
vertex-linking $2$--sphere with an almost normal tube along the edge
$3$ but it is isotopic to the vertex-linking $2$--sphere with a thin
edge-linking tube along the edge $3$. All of these but the
vertex-linking $2$--sphere with a thin edge-linking tube along the
edge $3$ is a Heegaard surface. It is curious that the genus $2$
example (the vertex-linking $2$--sphere with a thin edge-linking
tube along the edge $3$ and an almost normal tube along the thick
edge gives the classical example of an embedding of a genus $2$
handlebody in $S^3$ with one handle the trefoil knot and the other
unknotted but linking in such a way that the complement has closure
a handlebody. Just as in this example, the edge $3$ (the trefoil
knot) unwinds about the trivial handle; thus the genus $2$ example
is stabilized. Finally, the minimal layered-triangulation of $S^3$
is the only minimal layered-triangulation of a lens space with a
normal Heegaard surface (the vertex-linking $2$--sphere). In
addition, and as expect by the theory of sweep-outs, this
triangulation of $S^3$ \emph{must} have an  almost normal octagonal
$2$--sphere; it is the only octagonal Heegaard surface in a minimal
layered-triangulation of a lens space.

\vspace{.15 in}\noindent For $\rppp, S^2\times S^1$ and $L(3,1)$:

In each of these cases, we have two-tetrahedra
layered-triangulations. For $L(3,1)$ we use the minimal
(two-tetrahedra) $1/1$--layered-triangulation of the solid torus and
fold over along the thick edge.  Again, the argument is just a
matter of listing the possible normal and almost normal surfaces for
the lens spaces given in Theorem \ref{surface-in-lens}. In each case
we have the same characterization as in Theorem \ref{samelevel};
namely, a normal or almost normal surface is a Heegaard surface  if
and only if there is an almost normal tube at the same level as a
thin edge-linking tube or there is an almost normal tube along a
thick edge. There are no normal Heegaard surfaces.

\begin{cor} \label{spine-reduction} An almost normal Heegaard surface of genus $\ge 2$ in a minimal
layered-triangulation of a lens space is stabilized.
\end{cor}\begin{proof} If the genus is at least $2$, then by Theorem
\ref{samelevel} and the preceding paragraphs, there is an almost
normal tube along a thick edge or there is an almost normal tube at
the same level as a thin edge-linking tube. Furthermore, we showed
that an almost normal tube along a thick edge and any other tube
leads to a stabilization; and an almost normal tube along a thin
edge-linking tube and any other thin edge-linking tube, leads to a
stabilization. Hence, we only need to show that in the genus $2$
case, when we have an almost normal tube at the same level as a thin
edge-linking tube, we get a stabilization.

In this case, we may consider the thin edge-linking tube and the
almost normal tube in the same level-torus, sat $T$; thus we have a
spine, say $\Gamma$, for $T$ made up of the edge which the thin
edge-linking torus links and the core of the almost normal tube,
which may also be taken as an edge. Let $N(\Gamma)$ be a small
neighborhood of $\Gamma$ in $T$, which is a once-punctured torus,
and consider the small product neighborhood $N(\Gamma)\times I$.
$N(\Gamma)\times I$ is (up to isotopy) one of the handlebodies of
the genus $2$ Heegaard splitting. The other handlebody is the two
solid tori, say $\bbb{T}_1$ and $\bbb{T}_0$, complementary to the
interior of $T\times I$, which are  connected by the solid tube that
is the closure of $(T\setminus N(\Gamma))\times I$. Since $\Gamma$
is a spine for $T$, there is a meridional disk $D_0$ in the solid
torus $\bbb{T}_0$ with boundary in the interior of $N(\Gamma)\times
\{0\}$. On the other hand, since $N(\Gamma)\times \{0\}$ is a
punctured torus, there is an arc $\alpha$ in $N(\Gamma)\times \{0\}$
that crosses the boundary of $D_0$ precisely once. Let $D =
\alpha\times I$. Then $D_0$ and $D$ determine a stabilization for
the genus $2$ Heegaard splitting.

Note that following this stabilization, the genus $1$ Heegaard
splitting may not be normal or almost normal but after an isotopy,
it is a vertex-linking $2$--sphere with an almost normal tube along
a thick edge.\end{proof}

\begin{thm}\cite{wald-HS3, bon-otal} Up to isotopy a lens space has a
unique genus $g$ Heegaard splitting for every $g\ge 1$.\end{thm}

\begin{proof} It is sufficient to consider only irreducible Heegaard
surfaces and splittings with genus $\ge 1$. By \cite{cass-gord}, for
a lens space an irreducible Heegaard surface is strongly
irreducible; hence, by Rubinstein \cite{rubin-jac}, for any
triangulation such a Heegaard surface is isotopic to an almost
normal surface. In particular, this is true for a minimal
layered-triangulation of the lens space (for the lens space
$L(3,1)$, we can avoid the bad minimal layered-triangulation). Now,
from the above, if the Heegaard surface is irreducible, it must be
genus $g\ge 1$ and is the vertex-linking $2$--sphere with an almost
normal tube along a thick edge.\end{proof}

\subsection{Dehn fillings of knot and link manifolds} Our encounter with layered-triangulations
began with a problem on Dehn filling a knot-manifold (a compact,
orientable, irreducible $3$--manifold, with a single incompressible
torus boundary). In particular, we had constructed a one-vertex
triangulation, $\T$, of the knot-manifold, $X$, along with a
designated slope $\alpha$ on $\bdy X$.  We wanted  a one-vertex
triangulation of the Dehn filling of $X$ along $\alpha$ that
retained (extended) the triangulation $\T$ of $X$. In addition, we
needed to understand how the normal and almost normal surfaces in
the triangulation of the Dehn filling could meet the solid torus of
the filling.

So, suppose $X$ is a knot-manifold and  $\T$ is a one-vertex
triangulation of $X$. Let $\alpha$ be a slope on $\bdy X$. The
isotopy class of $\alpha$ in $\bdy X$ determines a unique triple of
intersection numbers with the three edges in the one-vertex
triangulation of the torus, $\bdy X$, say $\{p,q,p+q\}, 0\le p<q$.
Let $\bbb{T}$ denote a solid torus with the minimal
layered-triangulation extending the $p/q$-triangulation on its
boundary. We can attach $\bbb{T}$ to $X$ via a simplicial
isomorphism from $\bdy\bbb{T}\rightarrow \bdy X$ taking the edges
with intersection numbers $\{p,q,p+q\}$ in $\bdy\bbb{T}$ to those
with corresponding intersection numbers $\{p,q,p+q\}$ in $\bdy X$.
Hence, we get the closed $3$--manifold $X(\alpha)$ along with a
unique triangulation extending the triangulation $\T$ of $X$. We
denote this triangulation by $\T(\alpha)$ and say we have a {\it
triangulated Dehn filling} (of $X$ with respect to $\T$). Without
changing the triangulation $\T$, we get unique triangulated Dehn
fillings $\T(\alpha)$ of $X(\alpha)$ for all slopes $\alpha$. In
case the triple is either $\{0,1,1\}$ or $\{1,1,2\}$, we can use the
creased $3$--cell or the one-triangle M\"obius band to achieve a
triangulated Dehn filling. In these cases, there will be
identifications on the faces of $\T$ that are in the boundary torus.
One last remark. Typically, when working with triangulated Dehn
fillings, it is convenient to have that the frontier of a small
regular neighborhood of the layered solid torus is a normal torus.
This is equivalent to the frontier of a small neighborhood of the
boundary of $X$ being normal in $\T$. If this is the case, we say
{\it $\bdy X$ is normal in $\T$}. It is not necessary that the
boundary of a manifold be normal in a triangulation of the manifold.
However, we do show in \cite{jac-rub-blowup} that any compact
$3$--manifold with nonempty boundary, no component of which is a
$2$--sphere, admits a triangulation with normal boundary;
furthermore, we have such triangulations with all vertices in the
boundary and just one vertex in each boundary component. Hence,
normal boundary and minimal vertex triangulations are not
necessarily restrictions. Throughout our work we assume that any
triangulation of $X$, used in triangulated Dehn filling, has normal
boundary.

\vspace{.15 in} \noindent {\it Heegaard splittings of Dehn
fillings.}   The study of Heegaard splittings of Dehn fillings has
attracted much interest. A key to success in these efforts is an
understanding of how a strongly irreducible Heegaard surface can
meet a solid torus. This was done in \cite{mor-rubHS-curved},
followed by new proofs in \cite{shar-HS}, and has been studied
extensively by Y. Reick and E. Sedgwick
\cite{rieckHS,rieck-sedgHSstructure,rieck-sedgHSfinite}. A reader
familiar with these earlier results can easily see that
layered-triangulations of the solid torus provide an ideal
combinatorial environment for an investigation of how a normal or
almost normal, strongly irreducible, Heegaard surface can meet the
solid torus in a triangulated Dehn filling.

We give a combinatorial proof of a theorem that is in the spirit of
the main result in both the work of Moriah-Rubinstein
\cite{mor-rubHS-curved} and Scharlemann \cite{shar-HS}. We follow
the format of the statement in \cite{shar-HS}; however, here we give
a modest generalization and we note that possibility $3(c)$ in the
conclusion of Theorem \ref{HS-meets-torus} seems to have been missed
in the earlier versions. Following this result, we are able to give
a completely general analysis of how a strongly irreducible Heegaard
surface can meet a solid torus with a minimal layered-triangulation.
We extend our investigation and give examples in both
\cite{jac-rubfiniteHS} and \cite{jac-rub-blowup}.

Suppose $S$ is a normal or almost normal surface in a triangulated
$3$--manifold meeting a subcomplex  $\bbb{T}$ that is a solid torus
with a minimal layered-triangulation. Furthermore, suppose $S$ does
not meet $\bdy\bbb{T}$ in any trivial simple closed curves. From our
classification of normal and almost normal surfaces in a minimal
layered triangulation of a solid torus, it follows that $S$ can meet
$\bbb{T}$ in at most one of:\begin{enumerate}\item a collection of
meridional disks, possibly with an almost normal tube, or\item an
almost normal octagonal annulus (possibly) with thin edge-linking
tubes, or
\item a collection of normal edge-linking annuli (possibly) with
tubes, including the possibility of an almost normal tube.
\end{enumerate}

In case $(3)$  the components of $S\cap\bbb{T}$ nest in the sense
that the annuli must all be parallel (normally isotopic)
edge-linking annuli, each separates $\bbb{T}$ into two components,
one containing the vertex and, if $A_i$ and $A_j$ are components of
$S\cap\bbb{T}$ and $A_i$ separates $A_j$ from the vertex, then any
thin edge-linking tube for $A_i$ must run through a thin
edge-linking tube for $A_j$. Similarly, if there is an almost normal
tube between two normal surfaces and if both have thin edge-linking
tubes, then we also have the thin edge-linking tubes of one running
through the thin edge-linking tubes of the other. In either of these
situations we say we have {\it nested tubes} in $S\cap\bbb{T}$.

\begin{lem} \label{no-nesting} If in the above situation, $S$ is
a strongly irreducible Heegaard surface, then there are no nested
tubes in $S\cap\bbb{T}$.\end{lem}

\begin{proof}  If there were nested
tubes, then there is an edge $e$ and thin edge-linking tubes $t$ and
$t'$ about $e$ so that the meridian disk $d'$ for $t'$ contains the
meridian disk $d$ for $t$ and  only meets $S$ in $\bdy d$ (see
Figure \ref{f-no-nesting}). Let $C$ denote the annulus in the disk
$d'$ complementary to the interior of $d$.

\begin{figure}[h]

\psfrag{d}{\small $d$}\psfrag{e}{$e$}\psfrag{b}{\small $d'$}
\psfrag{c}{$C$}\psfrag{t}{$t'$}\psfrag{s}{$t$} {\epsfxsize = 2 in
\centerline{\epsfbox{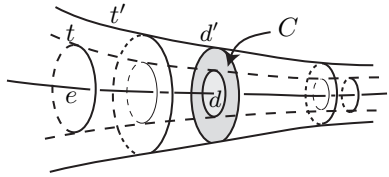}} } \caption{There are no
nested tubes.} \label{f-no-nesting}
\end{figure}

The annulus $C$ is contained in one of the handlebodies bounded by
$S$ and has both of its boundary components contained in $S$; so, it
must compress or $\bdy$-compress. If $C$ compresses, then using the
disk determined by the compression and the disk $d$, we have a
reduction of $S$, which is impossible. If $C$ $\bdy$-compresses,
then since $C$ does not separate, $C$ is not parallel into the
boundary of $S$ and so, we get a disk upon boundary compression that
along with the disk $d$ gives a weak reduction of $S$, which also is
impossible.\end{proof}

\begin{thm}\label{HS-meets-torus} Suppose $X$ is a knot-manifold with
a one-vertex triangulation $\T$,  $\alpha$ is a slope on $\bdy X$,
and $\T(\alpha)$ is the triangulated Dehn filling of $X$ with
respect to $\T$.  Let $\bbb{T}$ denote the solid torus of the Dehn
filling. If $S$ is a strongly irreducible Heegaard surface in
$X(\alpha)$ that is normal or almost normal in $\T(\alpha)$ and
$S\cap\bdy\bbb{T}$ has no components that are vertex-linking simple
closed curves, then $S$ meets $\bbb{T}$ in one of the
following:\begin{enumerate}\item a collection of meridional disks
along with possibly one other component that is two meridional disks
with an almost normal tube between them along an edge in
$\bdy\bbb{T}$, or
 \item a longitudinal, almost normal, octagonal annulus, or \item
a collection of edge-linking normal annuli and possibly one other
component that is  obtained either \begin{enumerate}\item by adding
an almost normal tube between two of the annuli, along an edge in
$\bdy\bbb{T}$, or\item by adding an almost normal tube to the
outermost annulus, along an edge in the complementary annulus, or
\item by adding an almost normal tube to the outermost annulus along
the thick edge.\end{enumerate}
\end{enumerate}\end{thm}

\begin{proof} First consider the situation where $\alpha$
does not have intersection numbers $\{0,1,1\}$ or
$\{1,1,2\}$.

Suppose $M = H\cup_S H'$ is a strongly irreducible Heegaard
splitting of $X(\alpha)$. We have the possibilities $(1), (2)$ and
$(3)$ from above, where in case $(3)$ we have no nesting. We analyze
each of these possibilities.

\vspace{.15 in}\noindent Case $(1)$ $S\cap\bbb{T}$ is a collection
of meridional disks, possibly with an almost normal tube.

Note that any almost normal tube from a meridional disk to itself or
between two meridional disks can be isotoped to an almost normal
tube along an edge in $\bdy\bbb{T}$. If such a tube were on a single
meridional disk, then we would have a weak reduction, which is a
contradiction. Hence, we are left only with possibility $(1)$ in the
statement of the theorem.

\vspace{.15 in}\noindent Case $(2)$ $S\cap\bbb{T}$ is an almost
normal octagonal annulus (possibly) with thin edge-linking tubes.

Let $B'$ be the complementary annulus to $A$ and suppose $B'\subset
H'$.  $B'$ must compress or $\bdy$--compress into $H'$. Since $B'$
has the slope of an edge, it is not trivial in $\bdy\bbb{T}$ nor can
it have the slope of a meridian. It follows that $B'$ can not
compress. We have that $B'$ must $\bdy$--compress.

If $B'$ $\bdy$--compresses to the outside of $\bbb{T}$, then we can
use this $\bdy$--compression to construct a properly embedded disk
$D$ in $H'$. If $D$ is not essential, then $S = A\cup A'$, where
$A'$ is an annulus parallel into $\bbb{T}$. But then we would have
$X$ in $H$ and $\bdy X$ incompressible in $H$; this is impossible.
So, if $B'$ $\bdy$--compresses outside of $\bbb{T}$, the disk $D$
must be essential. Note, however, by $S$ being strongly irreducible,
then we can have no thin edge-linking tubes for $A$, for otherwise,
we would get a weak reduction of $S$. It follows in this case that
$A$ is an almost normal octagonal annulus without thin edge-linking
tubes. Hence, $A$ is either longitudinal or is isotopic to a normal
edge-linking annulus. The former gives us conclusion $(2)$ of the
theorem and the latter gives us conclusion $(3)$.

Now, suppose $B'$ does not compress and does not $\bdy$--compress to
the outside of $\bbb{T}$. Thus $B'$ must $\bdy$--compress into
$\bbb{T}$. This is equivalent to saying that $A$ $\bdy$--compresses
into $B'$. But by Lemma \ref{no-comp}, this is impossible unless $A$
is a longitudinal almost normal octagonal annulus with no tubes.
Again, we have conclusion $(2)$ of the theorem.

\vspace{.15 in}\noindent Case $(3)$ $S\cap\bbb{T}$ is a collection
of normal edge-linking annuli (possibly) with tubes, including the
possibility of an almost normal tube.

Hence, the only possibility for $S\cap\bbb{T}$ is a number of copies
of the same annulus with  the outer most one (possibly) having thin
edge-linking tubes and (possibly) an almost normal tube either
between two of these annuli or attached to the outermost one. Notice
that any almost normal tube between any two of these annuli is
isotopic to an almost normal tube between them and along an edge in
$\bdy\bbb{T}$.

Let $A$ denote the outermost annulus (possibly) with tubes.  While
the argument here is similar to that above in Case $(2)$, we may
have an almost normal tube on $A$. However, if  $A$ does has an
almost normal tube, it may be isotoped to form a new thin
edge-linking tube on $A$ or it is at the same level as a thin
edge-linking tube, or along the thick edge, or, possibly, at the
same level as the thin edge-linking annulus associated with $A$.

So, we may assume we either have no almost normal tube or we have an
almost normal tube at the same level as a thin edge-linking tube, or
along the thick edge, or, possibly, at the same level as the thin
edge-linking annulus associated with $A$.

We shall consider, first, the case where we must have an almost
normal tube. If we have the almost normal tube along the thick edge,
then other thin edge-linking tubes could be unwound about the tube
along the thick edge and thus would give us a reduction. If we have
the almost normal tube at the level of a thin edge-linking tube or
at the level of the annulus associated with $A$, Then we have a
spine for a level torus at that level in $A$ and thus any other thin
edge-linking tubes could be unwound about this spine and thus would
give us a reduction. We conclude that $A$ is precisely one of:
\begin{enumerate}\item[(a)] an edge-linking annulus with a single
thin edge-linking tube and an almost normal tube at the same level
as this thin edge-linking tube, or
\item[(b)] an edge-linking annulus with an almost normal tube along
an edge in the complementary annulus, or \item[(c)] an edge-linking
annulus with an almost normal tube along the thick edge.
\end{enumerate}

But as in the proof of Corollary \ref{spine-reduction}, if we have
an almost normal tube at the same level as a thin edge-linking tube,
then we have a reduction. We conclude we can not have situation (a).
Situation (b) gives conclusion $3(b)$ and situation (c) gives
conclusion $3(c)$ of the theorem.

Finally, we consider the situation where we do not have an almost
normal tube; that is, the outermost edge-linking annulus $A$ does
not have an almost normal tube. Let $B'$ be the complementary
annulus to $A$ and suppose notation has been chosen so that
$B'\subset H'$. Then $B'$ must compress or $\bdy$--compress in $H'$.
As above, it can not compress; it therefore, must $\bdy$--compress.

If it $\bdy$--compresses to the outside, we get a properly embedded
disk $D$ in $H'$. If $D$ is essential, then $A$ has no tubes or we
would have a weak reduction of $S$; if it is inessential, then $S =
A'\cup A$, where $A'$ is an annulus parallel into $\bbb{T}$. Again,
as we saw above, this would place $\bdy\bbb{T}$ as an incompressible
torus in $H$, a contradiction. Hence, in this situation, we conclude
that $A$ has no tubes. This gives conclusion (3) of the theorem
where the collection has no component with tubes.

Since $A$ is normal, by Lemma \ref{no-comp}, $B'$ can not compress
or $\bdy$--compress into $\bbb{T}$.\end{proof}

While Theorem \ref{HS-meets-torus} is in a form analogous to those
in \cite{mor-rubHS-curved,shar-HS}, our version does not exclude
curves in the intersection between the Heegaard surface and the
boundary of the solid torus that have meridional slope. In fact, we
can extend our analysis to the general case, allowing all possible
curves of intersection between the Heegaard surface and the boundary
of the solid torus. This is a telling difference in using normal and
almost normal surfaces and layered-triangulations; we can address
situations that were totally intractable without these tools. We now
state the result in the general case; the earlier version provides a
rather direct comparison between our methods and the earlier
methods. Furthermore, by allowing inessential curves of
intersection, the proofs are not more difficult but there are more
possibilities for almost normal tubes, which makes the conclusions a
bit unwieldy.

\begin{thm}\label{HS-meets-torus-gen} Suppose $X$ is a knot-manifold with
a one-vertex triangulation $\T$,  $\alpha$ is a slope on $\bdy X$,
and $\T(\alpha)$ is the triangulated Dehn filling of $X$ with
respect to $\T$.  Let $\bbb{T}$ denote the solid torus of the Dehn
filling. If $S$ is a strongly irreducible Heegaard surface in
$X(\alpha)$ that is normal or almost normal in $\T(\alpha)$, then
$S$ meets $\bbb{T}$ in (possibly) a collection of vertex-linking
disks and one of the following:\begin{enumerate}\item a collection
of meridional disks along with possibly one other component that is
formed by adding an almost normal tube along an edge in
$\bdy\bbb{T}$ between two meridional disks, or between two
vertex-linking disks, or between a meridional disk and a
vertex-linking disk, or
 \item a longitudinal, almost normal, octagonal annulus, or \item
a collection of edge-linking normal annuli and possibly one other
component that is  obtained either \begin{enumerate}\item by adding
an almost normal tube  along an edge in $\bdy\bbb{T}$ between two of
the annuli, or between two vertex-linking disks, or between a
vertex-linking disk and the innermost edge-linking annulus, or\item
by adding an almost normal tube to the outermost annulus, along an
edge in the complementary annulus, or
\item by adding an almost normal tube to the outermost annulus along
the thick edge, or\end{enumerate}\item a vertex-linking disk with an
almost normal tube along the thick edge.
\end{enumerate}\end{thm}
\begin{proof} The really new part of the theorem is conclusion
$(4)$. We have this possibility since we now have that the
intersection between $S$ and $\bbb{T}$ might be a collection of
vertex-linking disks, possibly, with thin edge-linking tubes along
with a collection of edge-linking linking annuli, possibly, with
thin edge-linking tubes. However, we note that the nesting of tubes
in the case of just edge-linking annuli with tubes, extends to
having vertex-linking disks with tubes and edge-linking annuli with
tubes. But as before, there can be no nesting of tubes. Thus if we
have any edge-linking annuli, then the vertex-linking disk will not
have any thin edge-linking tubes and we are really back in the case
of considering edge-linking annuli.

So, we may assume the intersection between $S$ and $\bbb{T}$ is a
collection of vertex-linking disks with, possibly, the outermost one
having tubes and, possibly one other component that is formed by
adding an almost normal tube between two vertex-linking disks or
between a vertex-linking disk and the outermost vertex-linking disk
whit thin edge-linking tubes.

An almost normal tube on the outermost vertex-linking disk can be
isotoped to a thin edge-linking tube, which we also do should the
almost normal tube be along a thin edge in $\bdy\bbb{T}$,  or it is
along the thick edge, or at the same level as a thin edge-linking
tube. If we have an almost normal tube along a thick edge, then
there can be no other tubes or we would have a reduction. This gives
us conclusion $(4)$. If there is an almost normal tube at the same
level as a thin edge-linking tube, then in any case there must be a
reduction, which is a contradiction. So, we have conclusion $(4)$ or
we may assume that all the tubes from the outermost vertex-linking
disk are thin edge-linking.

We can order the thin edge-linking tubes from the top to the bottom.
Having done this, we can split the layered solid torus at the level
of the first tube. Using an initial segment of the
layered-triangulation of the solid torus, we have the Heegaard
surface $S$ meeting this solid torus in a collection of
vertex-linking disks and a thin edge-linking annulus, possibly with
thin edge-linking tubes and, possibly one other component formed by
adding an almost normal tube between two vertex-linking disks or
between this edge-linking annulus and a vertex-linking disk. Now, as
in earlier arguments the complementary annulus to the annulus we
have must compress or $\bdy$--compress into one side of the Heegaard
splitting. However, the slope of the boundary of this annulus can
not be meridional nor longitudinal; hence, the only possibility is
that it $\bdy$--compress into $\bbb{T}$. But by Lemma \ref{no-comp},
this is impossible.

We find that in this situation we have conclusion $(4)$ or we have a
collection of vertex-linking disks and, possibly, one other
component that is two vertex-linking disks with an almost normal
tube between them and along an edge in $\bdy\bbb{T}$.\end{proof}

We end this section with a summary which can be drawn from the
conclusions of the previous theorem. \begin{enumerate}\item The
intersection contains a meridional disk. There are only finitely
many slopes \cite{jac-sedg-dehn} for which this can happen. It is
typical to exclude these slopes; however, in some examples, the
normal surfaces in the knot space can be analyzed and often we can
make conclusions about such exceptional slopes. \item The
intersection contains a longitudinal, almost normal, octagonal
annulus.  The knot (core of the solid torus) is in the Heegaard
surface.\item The intersection includes any of the possibilities
with an edge-linking annulus. Then the Heegaard splitting is a
Heegaard splitting of the knot manifold $X$. In some of these
possibilities, we use the so-called Daisy Lemma of Rieck and
Sedgwick \cite{rieck-sedgHSstructure}.\item The intersection
contains the vertex-linking disk with an almost normal tube along
the thick edge. Then the Heegaard splitting is a Heegaard splitting
of the knot manifold $X$.

\end{enumerate}

\vspace{.15 in} \noindent {\it Efficient triangulations of Dehn
fillings.} The following result on efficiency of triangulated Dehn
fillings indicates some of the techniques for studying interesting
surfaces in Dehn fillings; however, as in the study of Heegaard
splittings of Dehn fillings, it requires knowledge of normal
surfaces in the knot-manifold that is being filled. We, again, point
out that triangulated Dehn fillings fix a triangulation of the
knot-manifold that remains for all Dehn fillings; so, in particular,
the collection of normal and almost normal surfaces and, hence, the
slopes of these surfaces in the knot-manifold are the same for all
fillings.

We consider when a $0$--efficient triangulation of a knot-manifold
carries over to a $0$--efficient triangulation of a Dehn filling of
the manifold. Recall from \cite{jac-rub0} that a triangulation of a
closed $3$--manifold is said to be {\it $0$--efficient} if and only
if the only normal $2$--spheres are vertex-linking; for a compact
$3$--manifold with boundary, a triangulation is {\it $0$--efficient}
if and only if the only normal disks are vertex-linking. The latter
condition guarantees that there are no normal $2$--spheres and the
manifold has incompressible boundary. In our consideration of Dehn
fillings, we shall assume we have a minimal triangulation of the
knot-manifold and consider triangulated Dehn fillings, extending
this triangulation. It follows from \cite{jac-rub0} that a minimal
triangulation of a compact, irreducible, $\bdy$--irreducible
$3$--manifold is $0$--efficient; so, a minimal triangulation of a
knot-manifold is $0$--efficient and  has just one vertex.

\begin{thm}\label{0-eff-dehn} Suppose $X$ is a knot-manifold and $\T$ is a minimal
triangulation of $X$. Then for all but finitely many slopes $\alpha$
the triangulated Dehn filling $\T(\alpha)$ is
$0$--efficient.\end{thm}
\begin{proof} Let $\bbb{T}(\alpha)$ denote the solid torus in the
triangulated Dehn filling $\T(\alpha)$. $\bbb{T}(\alpha)$ has a
minimal layered-triangulation determined by the slope $\alpha$ and
possibly is degenerate (a M\"obius band or creased $3$--cell).

Suppose $\Sigma$ is a normal $2$--sphere embedded in $\T(\alpha)$.
Since we have assumed $\T$ is minimal ($0$--efficient by
\cite{jac-rub0}), $\Sigma\cap\bbb{T}(\alpha)\ne\emptyset$. The
possibilities for the components of $\Sigma\cap\bbb{T}(\alpha)$
are:\begin{enumerate}\item [(1)]a collection of vertex-linking
disks, or\item[(2)] a collection of meridional disks, possibly, with
a collection of vertex-linking disks, or\item [(3)]a collection of
edge-linking annuli, possibly, with a collection of vertex-linking
disks.\end{enumerate}

Since the triangulation $\T$ for $X$ is fixed, there are only
finitely many slopes that are boundary slopes for normal (or almost
normal) surfaces in $\T$ \cite{jac-sedg-dehn}; hence, if we avoid
these slopes, then we may assume possibility (2) does not occur. In
our particular situation, the only slopes that need to be avoided
are slopes for normal planar surfaces in $\T$. These slopes can be
algorithmically determined from planar surfaces at the vertices of
the projective solution space for $\T$; and, in the case of
essential planar surfaces we need to avoid at most three slopes
\cite {gord-luecke-planar}. On the other hand, note that excluding
surgery slopes that are slopes of embedded normal surfaces in the
knot-manifold does not exclude possibility (3), as boundary slopes
of normal annuli in the solid torus come in several varieties and
{\it a priori} can match with slopes of normal surfaces embedded in
the knot-manifold, even if excluded as surgery slopes. However, we
see below that, almost magically, possibility (3) can not occur.

Assuming that (2) does not occur, we consider  disk components
(innermost disks) in
$\Sigma\setminus(\Sigma\cap\bdy\bbb{T}(\alpha))$. If such an
innermost disk were in $X$, then its boundary must be a
vertex-linking curve ($\T$ is $0$--efficient). But a component of
$\Sigma\cap X$ with a trivial boundary can only match up with a
vertex-linking disk in $\bbb{T}$, leaving only the possibility that
we are in (1) and $\Sigma$ is the vertex-linking $2$--sphere. Thus
all disk components in
$\Sigma\setminus(\Sigma\cap\bdy\bbb{T}(\alpha))$
 must be in $\bbb{T}(\alpha)$ and, therefore, vertex-linking.

If there is any curve in $\Sigma\cap\bdy\bbb{T}(\alpha)$ that does
not bound a disk component of
$\Sigma\setminus(\Sigma\cap\bdy\bbb{T}(\alpha))$, then there is a
component, say $P$, of
$\Sigma\setminus(\Sigma\cap\bdy\bbb{T}(\alpha))\subset X$, where all
curves in its boundary, except one, say $\delta$,  bound innermost
disks in $\Sigma$; we call $\delta$ an outermost boundary for $P$.
Each boundary curve of an innermost disk in $\Sigma$ is a
vertex-linking curve and thus $P$ can be capped off with disks in
$\bdy X$ leading to  $\delta$ bounding an embedded disk in $X$. But
$\bdy X$ is incompressible, it follows that $\delta$ must be
vertex-linking in $\bdy X$ as well. But then it too would bound an
innermost disk in $\bbb{T}(\alpha)$. We conclude that $\Sigma$ meets
$X$ in a planar surface all of whose boundary components bound
vertex-linking disks in $\bbb{T}(\alpha)$. It follows that
possibility (3) can not occur. We are left with only possibility
(1).

We note that in particular, our methods show that except for
possibly the finite number of excluded slopes, $X(\alpha)$ is
irreducible.

If we have (1), then we need to call upon deeper results from
\cite{jac-rub0}; namely, we show that if $\Sigma$ is not
vertex-linking and we are in (1), then  we can crush $\Sigma$ and
contradict the minimality of $\T$. There are a number of conditions
to verify before we can crush $\Sigma$. We refer the reader to
\cite{jac-rub0}, Section 4, and in particular, Theorem 4.1.

First, we must have $\Sigma$ bounding a $3$--cell that contains the
vertex. Since $X$ is irreducible and $\Sigma$ is isotopic into X,
$\Sigma$ bounds a $3$--cell. If the vertex is not in the $3$--cell
bounded by $\Sigma$, then we can engulf the vertex, getting a new
normal $2$--sphere bounding a $3$--cell, $B$, that contains the
vertex. (It is possible in the engulfing process we arrive at a
punctured $3$--sphere containing the vertex and having all boundary
components normal $2$--spheres; but as noted, $X(\alpha)$ is
irreducible and from above, normal $2$--spheres can not meet the
filling torus in other than vertex-linking disks. So we get the
desired conclusion.) We shall continue to call this $2$--sphere
$\Sigma$ and note that since we have excluded Dehn-fillings along
slopes of planar normal surfaces in $X$, then this new $2$--sphere
meets $\bbb{T}(\alpha)$ only in vertex-linking disks, and is not,
itself, vertex-linking.

Let $Y = X(\alpha)\setminus\open{B}$. Then $Y$ has a nice
cell-decomposition consisting of truncated-tetrahedra,
truncated-prisms, and triangular and quadrilateral product blocks.
We denote this cell-decomposition of $Y$ by $\C$.

We need to observe that there are not too many product blocks; i.e.,
$\P(\C)\ne Y$. However, since $Y$ does not meet $\bbb{T}(\alpha)$ in
only product blocks, we have that $\P(\C)\ne Y$. Next we need to
assure that we can fill product blocks in $Y$ to get trivial product
blocks, denoted $\P(Y)$. However, each region that needs to be
filled has a $2$--sphere boundary and is contained in $X$; since $X$
is irreducible, we can fill these blocks missing $\bbb{T}(\alpha)$.
So,  $\P(Y)\ne Y$.

Now, we must consider having too many truncated prisms; i.e.,
possible cycles of truncated-prisms that are not in $\P(Y)$. If we
have a cycle of truncated-prisms about a single edge, then we have a
compression  and obtain two new normal $2$--spheres by adding a
$2$--handle to $\Sigma$; neither of these $2$--spheres is
vertex-linking. Furthermore, one of these $2$--spheres bounds a
$3$--cell containing the $3$--cell $B$ and the vertex and meets
$\bbb{T}(\alpha)$ in a subcollection of the vertex-linking disks
from $\Sigma$. By Kneser's Finiteness Theorem \cite{kne}, this can
only happen a finite number of times. We shall continue to call the
normal, non vertex-linking $2$--sphere $\Sigma$. So, now we consider
having a cycle of truncated-prisms about different edges. In this
situation, we have a solid torus formed by the cycle of
truncated-prisms and we have either three disjoint annular bands on
$\Sigma$ meeting the boundary of this torus, each  a longitude of
the torus, or just one annular band on $\Sigma$ meeting the torus in
a curve that intersects the meridian of the torus three times. In
the first situation, the cycle of truncated-prisms must be in the
product region $\P(Y)$. In the second, by using a disk on $\Sigma$
complementary to the annular band, we have that the cycle of
truncated-prisms along with this disk determine a copy of $L(3,1)$.
However, this cycle of truncated-prisms does not meet
$\bbb{T}(\alpha)$, thus such a copy of $L(3,1)$ can be isotoped
into$X$ which contradicts $X$ being irreducible. We conclude there
are no cycles of truncated-prisms in $\C$ that are not in $\P(Y)$.

So, by \cite{jac-rub0}, we can crush $Y$ along $\Sigma$. Since
$\Sigma$ meets $\bbb{T}(\alpha)$ only in vertex-linking disks,
$\bbb{T}(\alpha)$ remains after crushing and thus the complement of
$\bbb{T}(\alpha)$ after crushing has closure $X$.  Since $\Sigma$ is
not vertex-linking it meets a tetrahedron of $\T$ in a
quadrilateral; such a tetrahedron disappears after crushing. Hence,
the new triangulation of $X$ would have fewer tetrahedra than $\T$.
But this contradicts $\T$ a minimal triangulation of $X$. This
completes the proof.\end{proof}

The only algorithm we know to decide if a given triangulation of a
$3$--manifold is minimal uses a solution to the homeomorphism
problem; whereas, there is a straight forward algorithm to decide if
a triangulation is $0$--efficient \cite{jac-rub0}. So, it seems a
bit unpleasant that we have to assume the triangulation is minimal.
But, note from the proof of Theorem \ref{0-eff-dehn} that the
problem we encountered is the existence of a normal planar surface
in $\T$ with all its boundary components vertex-linking but it is
not a vertex-linking disk (it has more than one boundary component).
The method of  proof in this situation, which comes from
\cite{jac-rub0}, shows that starting with $X$ irreducible and
$\bdy$--irreducible, then we can alter the triangulation $\T$ of $X$
to eliminate such planar surfaces. Having done this, then the above
proof works to show that all but finitely many triangulated Dehn
fillings of such a triangulation are $0$--efficient. We do not need
to know the given triangulation $\T$ is minimal; we can deform it so
that it is not only $0$--efficient but rid it of unwanted planar
normal surfaces with all boundary vertex-linking curves and make it
``super" $0$--efficient.

We have a similar result as to when we can obtain $1$--efficient
triangulated Dehn fillings. In Section 8, we defined a
$1$--efficient triangulation in the special case of
layered-triangulations of lens spaces. However, in working with
triangulated Dehn fillings, we need the full definition of a
$1$--efficient triangulation. In addition, as we needed methods from
\cite{jac-rub0} for the proof of Theorem \ref{0-eff-dehn}, we need
methods from \cite{jac-rub1} to establish this analogous result for
$1$--efficient triangulations. Thus, we have chosen to leave the
$1$--efficient case and include it in \cite{jac-rub1}.

\section{Layered-triangulations of Handlebodies}\label{lay-handle}
Having layered-triangulations of the solid torus and  genus one
manifolds, it is natural to question if there are analogous
triangulations for higher genera handlebodies and higher genera
manifolds. The answer is yes; there are direct generalizations. We
provide the definitions and a brief introduction here. While the
generalizations are straight forward, the level of complexity is
not. We have not given much time to the study of
layered-triangulations of handlebodies or layered-triangulations of
$3$--manifolds; however, we believe it reveals a very interesting
direction of study and we hope to follow-up on these ideas.

A solid torus has a minimal triangulation with just one tetrahedron.
In generalizing to handlebodies, we first determine the minimal
number of tetrahedra needed in a triangulation of a genus $g$
handlebody.

If we have a $3$--manifold with connected, genus $g$ boundary, then
for any triangulation, we have from Euler characteristic,
$$t = 3g+2+e_i-v_i,$$ where $t$ is the number of tetrahedra, $e_i$ is
the number of interior edges, and $v_i$ is the number of interior
vertices. Since, $e_i-v_i$  must be non negative, we have $3g-2$ as
a lower bound for the number of tetrahedra needed to triangulate a
$3$--manifold with connected boundary having genus $g$ (Euler
characteristic $1-g$). Notice if we assume the manifold is
irreducible and we have $e_i = 0$ (and hence, $v_i = 0$), then we
must have a handlebody.

For $g =1$, we layer a single tetrahedron onto the center (only
interior) edge of the one-triangle M\"obius band and we get a
one-tetrahedron solid torus. For higher genera, we note that a
one-vertex triangulation of a compact surface with one boundary
component and Euler characteristic $1-g$ has $3g-2$ interior edges.
Thus we are lead to layering $3g-2$ tetrahedra onto the interior
edges of such a surface. Specifically, let $S$ be a compact, surface
with one boundary component and Euler characteristic $1-g$.  Let
$\S$ be a one-vertex triangulation of $S$ and let
$e_1,\ldots,e_{3g-2}$ denote the $3g-2$ interior edges of $\S$. Let
$\mathbf{\Delta} = \{\td{\Delta}_1,\ldots,\td{\Delta}_{3g-2}\}$ be a
pairwise disjoint collection of oriented tetrahedra.  Let $K_1 =
S\cup_{e_1}\Delta_1$ denote the complex determined by layering
$\td{\Delta}_1$ onto $S$ along $e_1$ via an orientation reversing
attaching map, where $\Delta_1$ is the image of $\td{\Delta}_1$
after attaching (here the attaching map can take two faces of
$\td{\Delta}_1$ to the same face of $S$; so, the orientation
reversing condition may apply in this initial step, as well as later
in the attaching). Having defined $K_{i-1}$, $ 1< i<3g-2$, let
$K_i=K_{i-1}\cup_{e_i}\Delta_i$ be formed by layering
$\td{\Delta}_i$ onto $K_{i-1}$ along $e_i$ via an orientation
reversing attaching map.

Notice that by layering a tetrahedron on each interior edge, we have
a $3$--manifold at each point, except possibly at the vertex (there
is only one vertex). Since we have layered on all interior edges,
the vertex-linking triangles form a $2$--manifold and we can
determine its Euler characteristic.

To see this let $D$ denote the vertex-linking surface; let $\chi_D$
be its Euler characteristic, $f_D, e_D$ and $v_D$ be the number of
its faces, edges and vertices, respectively. $f_D = 4(3g-2)$, $e_D =
(3f_D-e_{\bdy D})/2 + e_{\bdy D}$, where $e_{\bdy D}$ is the number
of edges in the boundary of $D$. However, the number of edges in
$\bdy D$  is three times the number of faces in a genus $g$ surface
with a one-vertex triangulation, which is $12g-6$. Hence, $e_D =
6(3g-2)+6g-3$. Finally, $v_D$ is twice the number of edges in the
triangulation of the handlebody. In this case it is easy to see that
the number of edges in the triangulation, which are  the edges in
the once-punctured surface with Euler characteristic $1-g$, is
$3g-1$, and then we get a new edge for each tetrahedron we add;
hence, we have $(3g-1)+(3g-2)=6g-3$ edges and therefore $v_D =
12g-6$. (In general, if we have a compact $3$--manifold with
boundary and a one-vertex triangulation, having Euler characteristic
$1-g$, then $1-g = -t +f -e + 1$ and since the number of faces in
the boundary must be $4g-2$, the number of edges is $e = t+3g-1$.)
We conclude that
$$\chi_D = f_D-e_D+v_D=12g-8 - (24g -15) +12g-6 = 1.$$
So, $D$ is a disk and the underlying point set is a $3$--manifold.
We implicitly used that the layering is via orientation reversing
attaching maps to guarantee that the resulting identification space
is a manifold at the midpoints of the edges.

There are clearly no closed normal surfaces in such a triangulation;
so we have an irreducible $3$--manifold. Since the triangulation
collapses onto the once-punctured surface $S$, we have a handlebody
having Euler characteristic $1-g$ and thus has genus $g$. The
surface $S$ along with a fixed triangulation is called a \emph
{spine} or a \emph{$g$--spine}. Just as we considered the
one-triangle M\"obius band a degenerate layered-triangulation of the
solid torus, we shall consider a $g$--spine a degenerate
layered-triangulation of the genus-$g$-handlebody. We summarize in
the following proposition.

\begin{prop} A genus-$g$-handlebody has a one-vertex triangulation;
furthermore, a minimal triangulation is a one-vertex triangulation
having $3g-2$ tetrahedra.\end{prop}

\begin{figure}[htbp]
\psfrag{1}{\Large $\td{\Delta}_1$}
            \psfrag{2}{\Large $\td{\Delta}_2$}
            \psfrag{3}{\Large $\td{\Delta}_3$}
            \psfrag{4}{\Large $\td{\Delta}_4$}
            \psfrag{a}{\footnotesize $1$}
            \psfrag{b}{\footnotesize $2$}
            \psfrag{c}{\footnotesize $3$}
            \psfrag{d}{\footnotesize $4$}
            \psfrag{e}{\footnotesize $\td{1}$}
            \psfrag{f}{\footnotesize $\td{2}$}
            \psfrag{g}{\footnotesize $\td{3}$}
            \psfrag{h}{\footnotesize $\td{4}$}
             \psfrag{i}{\footnotesize $\delta$}

            \psfrag{D}{\small{$\partial$}}

            \psfrag{A}{\begin{tabular}{c}
            once punctured\\
   Klein bottle\end{tabular}}

        \vspace{0 in}
        \begin{center}
        \epsfxsize=3 in
        \epsfbox{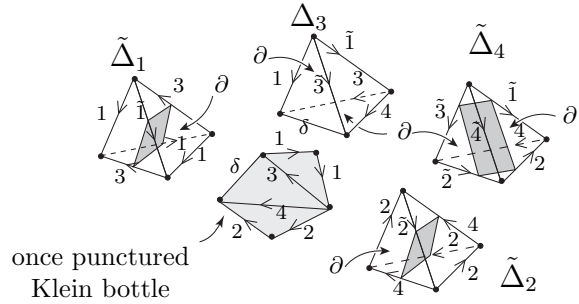}
        \caption{A minimal layered-triangulation of the genus 2 handlebody, having
$4$ tetrahedra and one vertex. The four tetrahedra are ``layered" on
a once-punctured Klein bottle. The six faces in the triangulation of
the genus $2$ boundary are indicated by the symbol $\bdy$. Normal
quadrilaterals $\{2,2\}, \{3,3\}$ and $\{\td{4},4\}$ are shown.}
        \label{f-genus2-layered}
        \end{center}

\end{figure}

See Figure \ref{f-genus2-layered} for an explicit example of a
minimal triangulation of a genus $2$ handlebody. In Figure
\ref{f-genus2-layered}, we have used the number $j$ for the edge
$e_j$ and the number $k$ with a tilde, $\td{k}$, for the edge
$\td{e}_k$, which eliminates the figure from being overwhelmed with
notation.

A triangulation of the genus-$g$-handlebody obtained by layering
$3g-2$ tetrahedra onto the interior edges of a $g$-spine will be
referred to as a {\it minimal layered-triangulation}. These are
analogous to the one-tetrahedron solid torus obtained by layering a
tetrahedron onto the interior edge of the one-triangle M\"obius
band. There is only one choice for layering in the genus one case.
However, for genus $g>1$, there are many choices. For an even genus
handlebody, we have the choice of layering on a nonorientable or an
orientable surface. Having made our choice of surface, we have the
many choices of a one-vertex triangulation of the surface; and we
have the various choices for the order in which we layer the $3g-2$
tetrahedra along the interior edges of the surface to get a minimal
layered-triangulation of the genus-$g$-handlebody. One is tempted to
place restrictions on the various choices but for the moment we have
no good basis for determining restrictions. In these higher genera
cases, there also are homotopy genus-$g$-handlebodies analogous to
the creased $3$--cell; and as the creased $3$--cell, they are not
manifolds. We do not know the role they might play here; so, at this
time we shall ignore these possibilities. Finally, we point out two
curious questions. The first should be easy to answer. Is a minimal
triangulation of a genus-$g$-handlebody a minimal
layered-triangulation? We suspect it is. The second is much more
formidable. How many distinct, up to isomorphism, minimal
layered-triangulations of the genus-$g$-handlebody are there? Recall
(see \cite{mosher-guide}) that, except for small genus, we do not
know the number of distinct triangulations, up to homeomorphism, of
the once-punctured surface of genus $g$. In our first case, genus
$2$, there is, up to isomorphism, $1$ one-vertex triangulation of
the once-punctured torus and 4 distinct one-vertex-triangulations of
the once-punctured Klein bottle. {\footnotesize\textsc{REGINA}}
\cite{burton-regina} has determined that there are $196$ non
isomorphic minimal triangulations (four tetrahedra) for the
genus-$2$-handlebody. This may seem horribly large. However, for the
once-punctured torus, alone, we have $384$ ways to layer four
tetrahedra; and all give a genus-$2$-handlebody. We have not
determined if all of the $196$ that {\footnotesize\textsc{REGINA}}
found are layered. We discuss these examples more below.

\subsection{ Layered-triangulations of a genus-$g$-handlebody} 
\begin{figure}[h]
\psfrag{b}{\scriptsize{$b$}}

            \psfrag{D}{$\td{\Delta}_{t}$}
            \psfrag{v}{\scriptsize{$\td{e}$}}\psfrag{u}{\scriptsize{$e$}}
\psfrag{x}{\scriptsize{$x$}}
            \psfrag{y}{\scriptsize{$y$}}\psfrag{z}{\scriptsize{$z$}}
            \psfrag{w}{\scriptsize{$w$}}
            \psfrag{3}{\scriptsize{$a$}}
            \psfrag{4}{\scriptsize{$b$}}\psfrag{5}{\scriptsize{$c$}}
            \psfrag{6}{\scriptsize{$d$}}
            \psfrag{L}{Layered}\psfrag{A}{\small{layering}}
            \psfrag{p}{\small{$2g+1)-gon$}}\psfrag{T}{\small{$T_1$}}
             \psfrag{S}{\small{$T_{t}$}}\psfrag{q}{\begin{tabular}{c}
            \small{once-punctured}\\
        \small{Klein bottle}\\\end{tabular}}
            \psfrag{s}{\begin{tabular}{c}
            \small{$t$}\\
       \small{layers}\\
            \end{tabular}}
            \psfrag{t}{\begin{tabular}{c}
            \small{$(t-1)$}\\
            \small{layers}\\
            \end{tabular}}

        \vspace{0 in}
        \begin{center}
\epsfxsize = 4 in \epsfbox{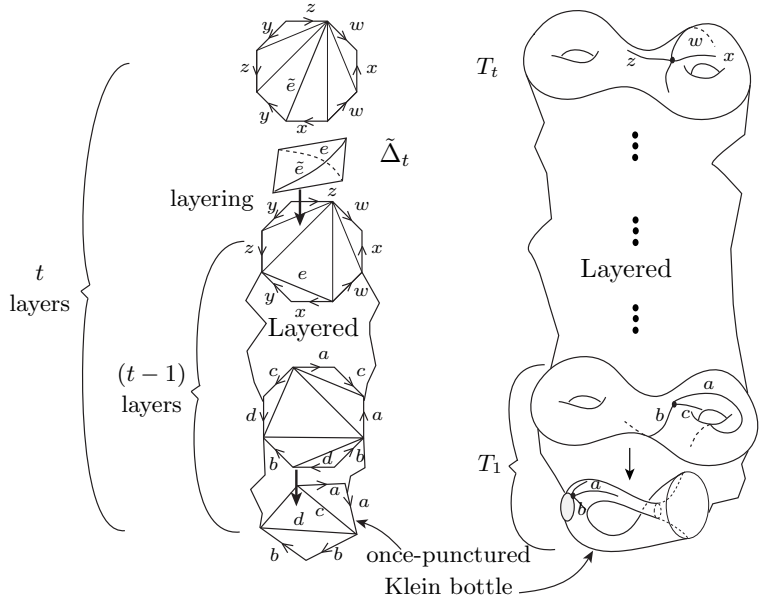} \caption{
Layered-triangulation of a handlebody. Shown for genus $2$.}
\label{f-layered-genusg-def}
\end{center}
\end{figure}

As in
the case of layered-triangulations of the solid torus, we
inductively define layered-triangulations of a handlebody. A
triangulation $\T_t$ is said to be  a {\it layered-triangulation of
the genus-$g$-handlebody}, with $t$--layers, if
\begin{enumerate}\item $\T_{-1}$ is a $g$-spine,
\item  $\T_0$ is a minimal layered-triangulation of the genus-$g$-handlebody
(obtained by layering $3g-2$ tetrahedra onto the interior edges of a
$g$-spine), \item $\T_{t} = \T_{t-1} \cup_e \td{\Delta}_t$ is a
layering along the edge $e$ of a layered-triangulation $\T_{t-1}$
having $t-1$ layers, $t \ge 1$. See Figure
\ref{f-layered-genusg-def}.
\end{enumerate}

If we have a genus-$g$-handlebody, there are infinitely many ways we
can place a one-vertex (minimal) triangulation of the genus $g$
surface on its boundary. If $\tau$ and $\tau'$ are one-vertex
triangulations on the boundary of a genus-$g$-handlebody, we shall
say that $\tau$ and $\tau'$ are {\it equivalent} if and only if
there is a homeomorphism of the handlebody taking $\tau$ onto
$\tau'$ via a simplicial isomorphism. This is an equivalence
relation on one-vertex triangulations on the boundary of a
genus-$g$-handlebody. We denote the equivalence class of $\tau$ by
$[\tau]$. For genus one, we have that the equivalence classes of
one-vertex-triangulations on the boundary of the solid torus are
parameterized by reduced fractions $p/q, 0< p<q$, along with $0/1$
and $1/1$. We do not have an analogous parametrization for higher
genera. However, we do have that these triangulations on the
boundary of the genus-$g$-handlebody can be extended to
layered-triangulations of the handlebody.

First, we give the following application of Lemma \ref{diag-flip},
which is very useful  to the theory of one-vertex triangulations of
$3$--manifolds.

\begin{lem}\label{pleated} Suppose $\T$ is a triangulation of
the compact $3$--manifold $M$ and $S$ is a component of $\bdy M$.
Furthermore, suppose $\mathcal{P}$ is the triangulation induced by
$\T$ on $S$ and $\P$ has one vertex. Then for any one-vertex
triangulation $\P'$ of $S$, there is a sequence
$\mathcal{T}=\mathcal{T}_0,\mathcal{T}_1,\ldots,\mathcal{T}_n = \T'$
of triangulations of $M$ and a sequence $\P = \P_0,\P_1\ldots,\P_n =
\P'$ of triangulations of $S$, where $\P_i$ is the triangulation of
$S$ induced by the triangulation $\T_i$ of $M$ and the triangulation
$\T_{i+1}$ is isomorphic to the triangulation obtained by layering
on $\T_i$ along an edge of $\P_i$. Hence, $\P_{i + 1}$ is isomorphic
to the triangulation obtained from $\P_i$ by a diagonal flip along
an edge of $\P_i, 0\le i\le n$. If $S'$ is a component of $\bdy M$
distinct from $S$, then the triangulation induced by $\T'$ on $S'$
is identical to that induced by $\T$ on $S'$.\end{lem}

\begin{proof} This is now a familiar argument in our work. The proof is by
induction on the length of the sequence $\P = \P_0,\P_1\ldots,\P_n =
\P'$. We have organized this result as a sequence of triangulations
on $M$ rather than triangulations on homeomorphic copies of $M$,
each obtained from the preceding by layering a tetrahedron along an
edge in the induced triangulation on the boundary.\end{proof}

We note in Lemma \ref{pleated} that in going from the triangulation
$\P$ on $S$ to the triangulation $\P'$ on $S$, we have embedded in
the triangulation $\T'$ of $M$ a sequence $S_0, S_1,\ldots,\S_{n-1}$
of subcomplexes of $\T'$, each is homeomorphic to $S$; furthermore,
the triangulation of $S_i$ is isomorphic to $\P_i$. We call such an
$S_i$ a {\it simplicial pleated surface} and think of this sequence
of pleated surfaces as {\it modifying the combinatorial structure on
$M$}, in a neighborhood of $S$, changing the induced triangulation
on $S$ from $\P$ to $\P'$.

Recall that if $\T_{\bdy}$ is a triangulation on the boundary of a
$3$--manifold $M$, we say the triangulation $\T$ of $M$ extends
$\T_{\bdy}$ if $\T$ restricted to $\bdy M$ is $\T_{\bdy}$ and the
vertices of $\T$ are precisely the vertices of $\T_{\bdy}$ (we do
not add vertices).

\begin{lem}\label{extend-one-vertex}  Suppose $\T_{\bdy}$ is a triangulation
on the boundary of a compact $3$--manifold $M$ and $\T_{\bdy}$ has
just one vertex in each component of $\bdy M$. Then $\T_{\bdy}$ can
be extended to a triangulation of $M$.\end{lem}

\begin{proof} For this, we first get a triangulation of $M$, ignoring $\T_{\bdy}$
for the moment,  that has all its vertices in the boundary and has
just one vertex in each boundary component. Such a triangulation can
be obtained starting with the observation of R.H.Bing \cite{bing1}
that any compact $3$--manifold with nonempty boundary admits a
triangulation with all of its vertices in the boundary. Having this,
one can use the operation of ``closing-the-book" (see
\cite{jac-rub0}) to get a triangulation of $M$  that not only has
all vertices in the boundary but has just one vertex in each
boundary component (no component of the boundary is a $2$--sphere or
projective plane, since by hypothesis each component of the boundary
admits a triangulation with just one vertex). This last step may
require adding tetrahedra by layering along edges in the boundary
but even if this is the case, at most two such tetrahedra are
required for each boundary component. Now, we have a triangulation,
say $\T$, of $M$ with all vertices in the boundary of $M$ and
precisely one vertex in each boundary component of $M$. The
triangulation $\T$ induces a one-vertex triangulation $\P$ on $\bdy
M$ and, up to isotopy, we may assume that $\P$ and $\T_{\bdy}$ have
the same vertex in each boundary component.

So, the problem is the triangulation $\T$ does not induce
$\T_{\bdy}$; hence, we must modify the PL structure on $M$ induced
by $\T$ in a neighborhood of each boundary component of $M$ from
$\P$ to $\T_{\bdy}$. By Lemma \ref{pleated}, we can do this, one
boundary component at a time; having modified the triangulation in a
neighborhood of a component $S'$ of $\bdy M$, we do not change the
triangulation induced on $S'$ when we move on to a boundary
component $S$ distinct from $S'$ and modify the triangulation in a
neighborhood of $S$. Thus, Lemma \ref{pleated} shows that we can,
through a  sequence of steps, modify the combinatorial structure in
a neighborhood of each boundary component of $M$ so as to eventually
include all boundary components.\end{proof}

While this result gives us immediately that any one-vertex
triangulation on the boundary of a genus-$g$-handlebody can be
extended to a triangulation of the handlebody, we want to sharpen
the conclusion so that not only have we extended the triangulation
but we have extended it to a layered-triangulation. This is easily
done by substituting a minimal layered-triangulation of the
genus-$g$-handlebody at the point in the preceding proof where we
used the triangulation coming from Bing's observation and
closing-the-book. The following corollary provides, yet, another way
to extend a one-vertex triangulation on the boundary of a solid
torus to a layered triangulation of the solid torus.

\begin{cor}\label{extend-handlebody} A one-vertex triangulation on the
boundary of a genus-$g$-handlebody can be extended to a
layered-triangulation of the handlebody.\end{cor}

\begin{proof} Let $M$ be a genus-$g$-handlebody and suppose $\P'$ is
a one-vertex triangulation on $\bdy M$. Now, $M$ can be triangulated
by any one of the minimal layered-triangulations, having $3g-2$
tetrahedra. Let such a triangulation be denoted $\T$ and let $\P$
denote the triangulation on $\bdy M$ induced by $\T$. Then we can
modify $\T$ in a neighborhood of $\bdy M$ via a sequence of pleated
surfaces to obtain a triangulation $\T'$ of $M$ that induces $\P'$
on $\bdy M$. The triangulation $\T'$ of $M$ is a
layered-triangulation because we started with $\T$ layered and each
subsequent step in pleating the PL structure is a layering along an
edge of the previous step.\end{proof}

What we provide here are existence proofs. One can, however,
construct a layered-triangulation extending a given
one-vertex-triangulation on the boundary of the
genus-$g$-handlebody. That is, given a one-vertex triangulation $\P$
of a closed surface $S$ and a complete system of closed curves for
the surface, then we can construct a layered-triangulation of a
genus-$g$-handlebody and a homeomorphism of the surface $S$ to the
boundary of the handlebody taking the given system of curves to
meridians for the handlebody.

From Corollary \ref{extend-handlebody}, if $\tau$ is a one-vertex
triangulation on the boundary of a genus-$g$-handlebody, $\tau$ can
be extended to a layered-triangulation of the handlebody. We call
such an extension a {\it $[\tau]$-layered-triangulation} of the
genus-$g$-handlebody. If in addition, such a
$[\tau]$-layered-triangulation of the genus-$g$-handlebody is
minimal among all $[\tau]$-layered-triangulations, we say it is a
{\it minimal $[\tau]$-layered-triangulation}. This should {\it not}
be confused with a minimal layered-triangulation of a
genus-$g$-handlebody, which has $3g-2$ tetrahedra. In the genus one
case, the minimal layered-triangulation (the one-tetrahedron torus)
is also the minimal $1/2$-layered-triangulation. In the higher
genera cases, we do not know what triangulations appear on the
boundaries of the minimal layered-triangulations. This is, itself,
an interesting question.  Below, we give examples of three different
triangulations of the genus $2$ surface that appear on the boundary
of minimal layered-triangulations of the genus-$2$-handlebody. These
were discovered by doing three examples by hand. There are nine
distinct (up to orientation preserving homeomorphism) triangulations
of the genus $2$ surface; and any one of these can be placed
infinitely many ways on the boundary of the handlebody. Here the
issue is which of these classes have a representative that can be
extended to a minimal layered-triangulation of the
genus-$2$-handlebody.  Ben Burton is adding code to
{\footnotesize\textsc{REGINA}} \cite{burton-regina} so such
questions can be aided by a computer search.

We have the following conjecture analogous to our conjecture in the
case of layered-triangulations of the solid torus.

\begin{conj} The minimal triangulation extending the one-vertex triangulation
$\tau$ on the boundary of a genus-$g$-handlebody is a minimal
$[\tau]$-layered-triangulation.\end{conj}

\vspace{.125 in}\noindent{\it The $L_g$ graph: Classification of
layered-triangulations of the genus-$g$--handlebody.} We construct a
graph to organize our study of one-vertex triangulations on the
boundary of a genus-$g$-handlebody and their extensions to
layered-triangulations of the handlebody. This is a preliminary
attempt for a genus $g$ analog  of the $L$--graph used above in the
study of one-vertex triangulations on the boundary of the solid
torus and their extensions to layered-triangulations of the solid
torus.

For genus $g$, we shall denote this graph by $L_g$. The $0$-cells of
$L_g$ are in one-one correspondence with equivalence classes of
one-vertex triangulations on the boundary of a genus-$g$-handlebody,
where equivalence is up to homeomorphism of the handlebody (we allow
orientation reversing as well as preserving homeomorphisms). The
$1$--cells are in one-one correspondence with diagonal flips on
edges of these triangulations. Specifically, there is a $1$--cell
between the $0$--cells corresponding to the equivalence classes
$[\tau]$ and $[\tau']$ if and only if some representative of
$[\tau']$ can be obtained from $\tau$ via a diagonal flip along an
edge of $\tau$ (hence, a representative of $[\tau]$ can be obtained
from $\tau'$ by a diagonal flip). We can perform diagonal flips on
any edge ($6g-3$ edges) in one of these triangulations, since an
edge must be adjacent to two distinct triangles. Note, the
$L_1$--graph is isomorphic to the $L$--graph above.

It follows from Proposition \ref{diag-flip} that $L_g$ is connected.
Also, by changing the order of independent diagonal flips, it is
easy to see that $L_g$, $g>1$, has cycles. We have the following
relationship between edge paths in $L_g$ and layered-triangulations
of the genus-$g$-handlebody.

\begin{prop} The $[\tau]$-layered-triangulations of a genus-$g$-handlebody
are in one-one correspondence with the paths in  $L_g$ that begin at
the vertex $[\t]$ and end at some $g$-spine.\end{prop}

A minimal $[\t]$--layered-triangulation corresponds to a minimal
edge path from the vertex $[\t]$ to a $g$--spine. There are various
possibilities for such a path; hence, unlike the genus one case,
there probably are a number of distinct minimal
$[\t]$--layered-triangulations. The answers to such questions will
come from the study and understanding of the properties of the $L_g$
graph and possibly its enlargement to a complex much like the
complex used by J.Harer \cite{harer} and others (see L. Mosher
\cite{mosher-guide}) to study the mapping class group.

\vspace{.125 in}\noindent{\it Normal surfaces in
layered-triangulations of handlebodies.} One of our goals is to
undertake an analysis of layered-triangulations of higher genera
handlebodies analogous to that for the solid torus. A major step
toward understanding and using layered-triangulations of the solid
torus is the classification of normal and almost normal surfaces and
their boundary slopes in minimal $p/q$--layered-triangulations of
the solid torus. While we have minimal
$[\tau]$--layered-triangulations in the higher genera case, we do
not expect to get a sharp classification like we had for the solid
torus. However, we do think there should be a useful organization of
the theory of normal surfaces in a layered-triangulation of a
genus-$g$-handlebody, as well as an understanding of their boundary
slopes.  For example, the layering serves as a height function and
the normal and almost normal surfaces are in a good position
relative to this height function. The local minima occur in the
$g$-spine; if critical values are introduced, they occur in a
tetrahedron upon layering and correspond to a saddle singularity
formed by a ``banding quadrilateral"; and the maxima are in the
boundary.  As a step toward what we might expect, we provide a small
example by classifying the normal surfaces in one of the minimal
layered-triangulations of the genus-$2$-handlebody. Recall that
while there are infinitely many incompressible distinct annuli and
non orientable surfaces embedded in a solid torus, a
layered-triangulation sees only finitely many of these as normal
surfaces. For higher genera handlebodies there are two constructions
of interesting families of embedded incompressible surfaces that are
consistent with the generation of normal surfaces. Namely ``stair
step surfaces" formed by adding several copies of an incompressible
normal surface to an incompressible normal annulus and ``spun
surfaces" formed by adding several copies of a incompressible normal
annulus to an incompressible normal surface (Dehn twisting about the
annulus). In normal surface parlance, surfaces of the form $kF + A$
and $F + kA$, where $F$ is an incompressible normal surface and $A$
is an incompressible normal annulus, both satisfying the same
quadrilateral conditions. In fact, we see in the simple example
below, a minimal layered-triangulation of a genus-$2$-handlebody, a
disk and annulus that generate an infinite family of ``spun" disks
by Dehn twisting a disk about the annulus.

\vspace{.125 in}\noindent{\bf Example:} Consider the minimal
layered-triangulation of the genus-$2$-handlebody given in Figure
\ref{f-genus2-layered}. To determine the normal surfaces in this
layered-triangulation, we use quadrilateral coordinates. There are
$12$ variables (the $12$ normal quadrilaterals in the four
tetrahedra) and no equations. (For quadrilateral coordinates, there
is one equation for each {\it interior} edge; in this example, there
are no interior edges.) To simplify notation, if $e$ and $f$ are
opposite edges in a tetrahedron, we shall denote the quadrilateral
separating these edges by $\{e,f\}$ and we shall use $j$ for the
edge $e_j$ and $\td{k}$ for the edge $\td{e}_k$.  In the example,
using this notation,  the three quadrilaterals in the tetrahedron
$\td{\Delta}_1$ are: $\{1,1\}, \{3,3\}$, and $\{\td{1},1\}$; in the
tetrahedron $\td{\Delta}_2$, the  quadrilaterals are: $\{2,2\},
\{4,4\}$ and $\{\td{2},2\}$; in the tetrahedron $\td{\Delta}_3$, the
quadrilaterals are: $\{1,4\},\ldots$, and so on. We  show the
quadrilaterals $\{2,2\}, \{3,3\}$ and $\{\td{4},4\}$ in Figure
\ref{f-genus2-layered}. Each quadrilateral is a quad solution and
the $12$ quads give a complete set of fundamental quadrilateral
solutions; in fact, they are vertex solutions. A quadrilateral
solution determines a unique normal solution, up to some multiple of
the vertex-linking disk. The  fundamental surfaces in the example
are:

\begin{tabbing}\=$\{1,1\} \Rightarrow$ \= meridional disk\hspace{.25
in}\=$\{3,3\} \Rightarrow$ \=M\"obius band\= \hspace{.25
in}\=$\{\td{1},1\}
\Rightarrow$ \=annulus\\
\>$\{2,2\} \Rightarrow$ \> meridional disk \>$\{4,4\} \Rightarrow$
\>M\"obius band\> \>$\{\td{2},2\}
\Rightarrow$ \>annulus\\
\>$\{1,4\} \Rightarrow$ \> meridional disk\>$\{\hspace{.02
in}\delta,\td{1}\} \Rightarrow$ \>annulus\> \>$\{\td{3},3\}
\Rightarrow$ \>meridional disk\\
\>$\{2,\td{3}\} \Rightarrow$ \> meridional disk\>$\{\td{1},\td{2}\}
\Rightarrow$ \>meridional disk\> \>$\{\td{4},4\}
\Rightarrow$ \>annulus.\\
\end{tabbing}

Hence, any normal surface in this minimal layered-triangulation of
the genus-$2$-handlebody is a disk, annulus, or M\"obius band. This
is as  expected; however, we do get some complexity  with just these
four tetrahedra that we did not see in layered-triangulations of the
solid torus. We have infinitely many distinct normal surfaces.  For
example, each sum $\{3,3\} + k\{\td{4},4\}, k=0,1,2,\ldots$ is a
meridional disk. These additions are the same as Dehn twisting the
disk determined by $\{3,3\}$ around the annulus determined by
$\{\td{4},4\}$; the latter is a thin edge-linking annulus about the
edge $e_{\td{4}}$ and the former is  a meridional disk orthogonal to
this edge. The disks formed in this way are isotopic but not
normally isotopic. We say more about this in the next section. There
is another interesting example here;  the sum of the two quad
solutions $\{1,1\} + \{1,4\}$ determines a meridional disk; the
normal solution for the disk is a fundamental solution (actually a
vertex solution) in normal solution space but its image in quad
space is not even fundamental. The sum  in normal solution space of
the disk determined by $\{1,1\}$ and the disk determined by
$\{1,4\}$  have the vertex-linking disk as a component.

\section{Layered-triangulations of $3$--manifolds}

Layered-triangulations for general $3$--manifolds are defined
analogously to those for lens spaces using genus $g$ Heegaard
splittings and layered-triangulations of genus-$g$-handlebodies.

We shall use compact surfaces having one boundary component and
Euler characteristic $1-g$ as degenerate genus-$g$-handlebodies,
just as we used the M\"obius band as a degenerate solid torus.
Furthermore, a minimal triangulation (one-vertex) of these surfaces
will be considered a (degenerate) layered-triangulation. Unlike the
one-triangle M\"obius band, there are many minimal triangulations of
these surfaces; as pointed out above, there are four distinct
minimal triangulations of the once-punctured Klein bottle. In the
genus $1$ case, we allowed the creased $3$--cell as a degenerate
solid torus and its one-tetrahedron triangulation as a degenerate
layered-triangulation. An analogous situation occurs for higher
genera. However, we have not investigated the role these degenerate
triangulations play; we suspect the role might be as for genus $1$,
where the creased $3$-cell was involved in minimal
layered-triangulations (Heegaard splittings) of small examples.

\vspace{.125 in}\noindent{\it $2$--symmetric triangulations of
closed orientable surfaces.} Recall that we defined a simplical
attachment between the boundaries of two solid tori with
layered-triangulations, incorporating an attachment in the cases
where one of the  factors was degenerate. A simplicial attachment of
the two-triangle torus to the one-triangle M\"obius band was a
folding along any one of the three edges of the two-triangle torus.
Thus the torus is a two-sheeted branched covering of the M\"obius
band, where the branching locus is the boundary of the M\"obius
band. The one-triangle triangulation of the M\"obius band lifts to
the two-triangle triangulation of the torus inducing a simplicial
involution on the torus fixed along an edge. There is such a
symmetry in the two-triangle torus about each edge.

\begin{figure}[h]

            \psfrag{a}{\footnotesize{$a$}}\psfrag{b}{\footnotesize{$b$}}
\psfrag{c}{\footnotesize{$a'$}}
            \psfrag{d}{\footnotesize{$b'$}}\psfrag{S}{$S_g$}
            \psfrag{B}{$B_g$}
            \psfrag{1}{$T_1$}
            \psfrag{2}{$T_2$}\psfrag{5}{$T_5$}
            \psfrag{3}{$T_3$}
            \psfrag{8}{$T_8$}
            \psfrag{I}{separating $2$--symmetry}
            \psfrag{J}{non separating $2$--symmetry}

        \vspace{0 in}
        \begin{center}
\epsfxsize = 4 in \epsfbox{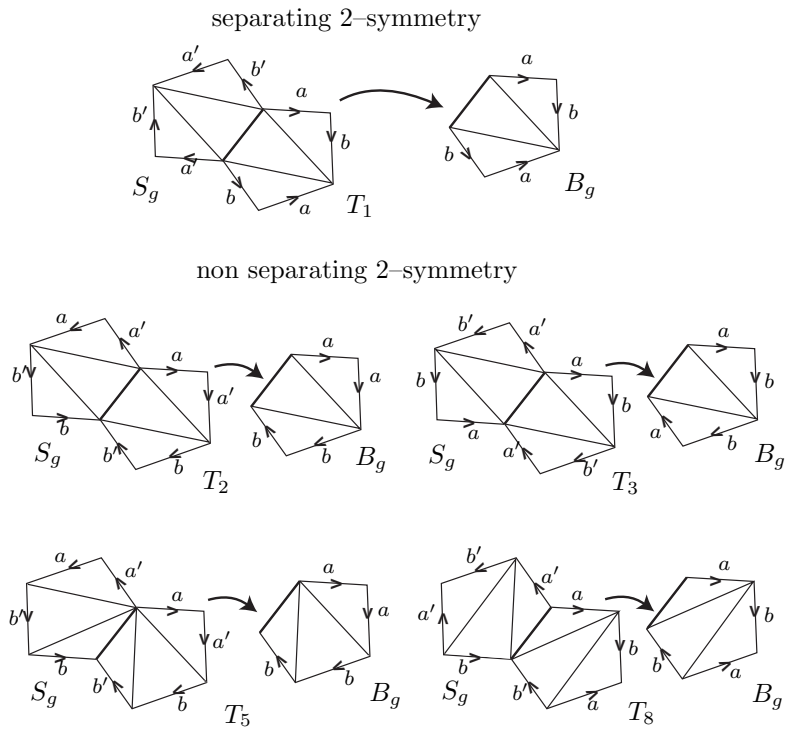} \caption{ The five
$2$--symmetric triangulations of a closed, orientable surface of
genus $2$ with orbit surfaces and quotient triangulations.}
\label{f-2-symmetry}
\end{center}
\end{figure}

There is an analogous situation for higher genera surfaces. If we
have a compact surface $B_g$ with one boundary component and Euler
characteristic $1-g$, then the closed orientable surface, $S_g$,
with Euler characteristic $2-2g$ is a branched double cover over
$B_g$ with branching locus the boundary of $B_g$. Hence, if we have
a minimal (one-vertex) triangulation of $B_g$ it lifts to a minimal
(one-vertex) triangulation of $S_g$; furthermore, such a
triangulation of $S_g$ has a simplicial involution fixed along an
edge, the edge over the branching locus ($\bdy B_g$).  We call such
a triangulation of $S_g$ {\it $2$--symmetric}; if the fixed edge,
under the simplicial involution, is separating, we say we have a
{\it separating $2$--symmetry}; and if it is non separating, we say
we have a {\it non separating $2$--symmetry}. The quotient manifold
in the case of a separating $2$--symmetry is orientable and in the
case of a non separating $2$--symmetry,  is non orientable. The
$2$--symmetric triangulation of $S_g$ induces a minimal
triangulation of $B_g$ so that the projection map is simplicial; and
conversely, a minimal triangulation of $B_g$ lifts to a
$2$--symmetric triangulation of $S_g$ so the projection map is
simplicial. In Figure \ref{f-2-symmetry}, we give the $2$--symmetric
triangulations for a genus $2$ closed, orientable surface. Chord
diagrams, such as those given in \cite{mosher-guide}, are very
helpful in recognizing $2$--symmetric triangulations; the chord
diagrams of a $2$--symmetric triangulation are symmetric about a
chord. For the examples in Figure \ref{f-2-symmetry}, we also use
the notation for the distinct minimal triangulations of the closed,
orientable surface of genus $2$ given in \cite{mosher-guide}. In the
case of genus $2$, where there are $9$ distinct minimal
triangulations, there are five $2$--symmetric triangulations: $T_1,
T_2, T_3, T_5$ and $T_8$. Since the $2$--symmetric triangulations of
a genus $g$ surface are in one-one correspondence with the minimal
triangulations of a compact surface with one boundary component and
Euler characteristic $1-g$, the ratio of $2$--symmetric
triangulations to the totality of minimal triangulations of a
closed, orientable surface of genus $g$ approaches zero as the genus
increases. More about this later.

\begin{figure}[htbp]
\psfrag{1}{$(1)$}\psfrag{2}{$(2)$}\psfrag{3}{$(3)$} \psfrag{0}{
$(0)$}

            \psfrag{5}{\Large {$T_5$}}
            \psfrag{4}{\Large {$T_1$}}
            \psfrag{a}{\tiny{$1$}}
            \psfrag{b}{\tiny{$2$}}
            \psfrag{c}{\tiny{$3$}}
            \psfrag{d}{\tiny{$4$}}
            \psfrag{e}{\tiny{$\td{1}$}}
            \psfrag{f}{\tiny{$\td{2}$}}
            \psfrag{g}{\tiny{$\td{3}$}}
            \psfrag{h}{\tiny{$\td{4}$}}
             \psfrag{i}{\tiny{$\delta$}}

            \psfrag{D}{$\partial$}

            \psfrag{A}{\begin{tabular}{c}
            \small{once-punctured}\\
   \small{Klein bottle}\end{tabular}}\psfrag{T}{\begin{tabular}{c}
            \small{once-punctured}\\
   \small{torus}\end{tabular}}\psfrag{F}{\small{folding}}

   \psfrag{B}{\small{boundary}}\psfrag{L}{\small{Layered}}
   \psfrag{S}{\Large{$S^2\times
   S^1$}}\psfrag{P}{\Large{$\rppp\#\rppp$}}

        \vspace{0 in}
        \begin{center}
        \epsfxsize=4 in
        \epsfbox{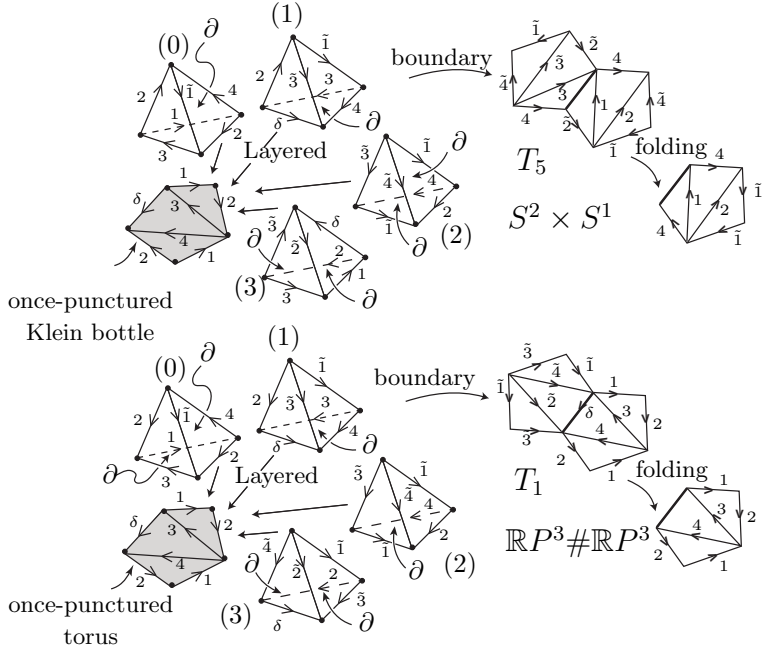}
        \caption{Two examples of triangulated Heegaard splittings of
genus $2$. The first gives a four-tetrahedron, layered-triangulation
of $S^2\times S^1$; and the second gives a four-tetrahedron,
(minimal) layered-triangulation of $\rppp\#\rppp$. For $S^2\times
S^1$, the Heegaard splitting is stabilized and for $\rppp\#\rppp$,
the Heegaard splitting is minimal genus but is a reducible
splitting.}
        \label{f-triang-HS}
        \end{center}

\end{figure}

\vspace{.15 in}\noindent {\it Definition.}  Suppose $M$ is a closed
$3$--manifold and $M = H\cup_{\bdy H\rightarrow\bdy H'}H'$ is a
genus $g$ Heegaard splitting of $M$. If $H$ and $H'$ can be endowed
with a layered-triangulation so that the attaching map $\bdy
H\rightarrow\bdy H'$ is simplicial, then the layered-triangulations
on $H$ and $H'$ give a triangulation of $M$. In this case, we call
the induced triangulation of $M$ a {\it layered-triangulation} and
more specifically a {\it genus $g$ layered-triangulation}. (Note
that the four tetrahedra triangulation of $S^2\times S^1$ shown in
Figure \ref{f-triang-HS} indicates that we should distinguish the
genus of a layered-triangulation; we typically think of a
layered-triangulation of this manifold being genus $1$ but clearly
any manifold with a genus $g$ layered-triangulation has a genus $g'$
layered-triangulation for all $g' >g$.) If
 one of the handlebodies is degenerate, then only one of $H$ or
$H'$ can be degenerate (say $H'$) and in this situation the induced
triangulation on $\bdy H$ must be $2$--symmetric and the attaching
map from $\bdy H\rightarrow\bdy H'$ is the projection of the
$2$--symmetric triangulation on $\bdy H$ to the $g$--spine $H'$. If
$e$ is the fixed edge under the $2$--symmetry of $\bdy H$, we call
the attaching map a {\it folding (of the triangulation of $\bdy H$)
along the edge $e$}. While the situation is more complex than in the
genus one case, we still have that a fixed layered-triangulation of
a $3$--manifold can be viewed in a number of different ways as a
union of two handlebodies, each with a layered triangulation.

\begin{thm} Every closed, orientable $3$--manifold admits a
layered-triangulation.\end{thm}
\begin{proof} Suppose $M$ is a closed, orientable $3$--manifold.
Then $M$ admits a  Heegaard splitting, $M = H\cup_S H'$.  Place any
one-vertex triangulation $\P$ on $S$. Then by Theorem
\ref{extend-to-layered}, $\P$ may be extended to a
layered-triangulation of $H$ and to a layered-triangulation of $H'$.
This gives a layered-triangulation of $M$.\end{proof}

In \cite{jac-rub0}, we provided a constructive proof that every
closed, orientable and irreducible $3$--manifold admits a one-vertex
triangulation. However, in \cite{jac-rub0}, we were interested in
$0$--efficient triangulations and placed additional conditions on
the triangulation; the methods there do not necessarily extend to
reducible $3$--manifolds. On the other hand, there are many ways to
arrive at one-vertex triangulations for $3$--manifolds. The previous
theorem provides one such method as a layered-triangulation has one
vertex. Thus as a corollary, we have:

\begin{cor} Every closed, orientable $3$--manifold admits a
one-vertex triangulation.\end{cor}

\begin{figure}[h]

            \psfrag{a}{{$a$}}\psfrag{b}{\small{$b$}}
\psfrag{c}{\small{$c$}}
            \psfrag{d}{\small{$d$}}\psfrag{T}{\Large{$\P''\approx\P'$}}
            \psfrag{S}{\Large{$\P = \P_0$}}\psfrag{R}{\Large{$\P_i$}}
            \psfrag{Q}{\large{$S_i$}}\psfrag{1}{\footnotesize{$1$}}
            \psfrag{2}{\footnotesize{$2$}}\psfrag{4}{\footnotesize{$\td{1}$}}
            \psfrag{3}{\footnotesize{$3$}}
\psfrag{P}{\begin{tabular}{c}simplicial\\
           pleated\\ surfaces\\\end{tabular}}
            \psfrag{L}{Layered}
            \psfrag{f}{$f:\bdy H\rightarrow\bdy H'$, simplicial}
\psfrag{G}{\Large{$H'$}}\psfrag{H}{\Large{$H$}}

        \vspace{0 in}
        \begin{center}
\epsfxsize = 3 in \epsfbox{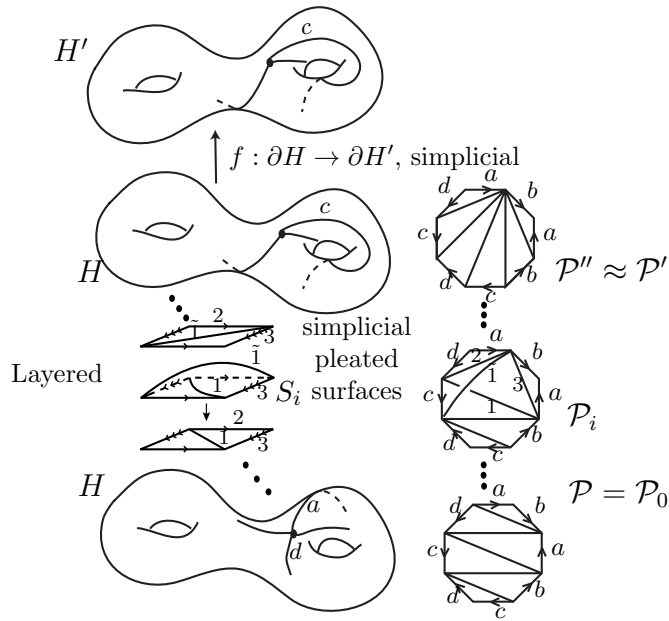} \caption{
Modification of combinatorial structure through a sequence of
simplicial pleated surfaces making the attaching map a simplicial
isomorphism.} \label{f-layered-general}
\end{center}
\end{figure}
There is a more constructive way to give a closed orientable
$3$--manifold a layered-triangulation, when it is given by a
Heegaard splitting. Specifically, suppose $M = H\cup_{f:\bdy
H\rightarrow\bdy H'} H'$ is a Heegaard splitting of $M$. Select any
minimal layered-triangulations $\T$ for $H$ and $\T'$ for $H'$. If
$g$ is the genus of the Heegaard splitting, then $\T$ and $\T'$ each
have $3g-2$ tetrahedra and, in fact, may be chosen to be isomorphic.
Let $\P$ and $\P'$ be the induced one-vertex triangulations on $\bdy
H$ and $\bdy H'$ respectively, and let $\P''= f^{-1}(P')$ be the
one-vertex triangulation on $\bdy H$ given by pulling the
triangulation $\P$ on $\bdy H'$ back to $\bdy H$. Then by Theorem
\ref{pleated} there is a sequence
$\mathcal{T}=\mathcal{T}_0,\mathcal{T}_1,\ldots,\mathcal{T}_n =
\T''$ of triangulations of $H$ and a sequence $\P =
\P_0,\P_1\ldots,\P_n = \P''$ of triangulations of $\bdy H$, where
$\P_i$ is the triangulation of $S$ induced by the triangulation
$\T_i$ of $M$ and the triangulation $\T_{i+1}$ is isomorphic to the
triangulation obtained by layering on $\T_i$ along an edge of
$\P_i$. Furthermore,  in going from the triangulation $\P$ on $\bdy
H$ to the triangulation $\P''$ on $\bdy H$, we have embedded in the
triangulation $\T''$ of $M$ a sequence $S_0, S_1,\ldots,\S_{n-1}$ of
subcomplexes of $\T''$, each is homeomorphic to $\bdy H$. The
triangulation of $S_i$ is isomorphic to $\P_i$. Hence, we can go
through a sequence of simplicial pleated surfaces building a {\it
new} layered-triangulation for $H$ that modifies the combinatorial
structure on $\bdy H$ so the attaching map $f$ is a simplicial
attachment. We provide a picture of this in Figure
\ref{f-layered-general}.

\vspace{.15 in}\noindent{\it Triangulated Heegaard splittings.} We
now give what seems to be the most interesting method of
representing a $3$--manifold via Heegaard splittings and
layered-triangulations. Let $H$ be a handlebody with a
layered-triangulation where the induced triangulation on the
boundary is $2$--symmetric along the edge $e$. If $M$ is the
$3$--manifold obtained by folding the boundary along the edge $e$,
then $M$ has a layered-triangulation and we say $M$ is presented by
a {\it triangulated Heegaard splitting}.

\begin{thm}\label{triang-HS} Every closed, orientable $3$--manifold can be presented
by a triangulated Heegaard splitting.\end{thm}
\begin{proof} Consider any layered-triangulation of $M$. This gives
$M = H\cup H'$, where $H$ and $H'$ both have layered-triangulations
and the attaching map is simplicial. If we open either $H$ or $H'$
where they are layered onto a $g$-spine, then we have $M$ presented
by a triangulated Heegaard splitting.\end{proof}

A triangulated Heegaard splitting of a lens space, where the
layered-triangulation of the solid torus is a minimal
$p/q$-layered-triangulation, is the method we used to present lens
spaces in the applications of layered-triangulations given in
Section 8.  Presenting a $3$--manifold by a triangulated Heegaard
splitting, where the layered-triangulation of the handlebody is a
minimal $[\tau]$--layered-triangulation, gives a compelling
direction for study of higher genera manifolds. There are numerous
questions that these presentations raise. For example, in the
generation of $3$--manifolds via triangulations, a major problem is
that the probability of getting a $3$--manifold by face pairings of
tetrahedra goes to $0$ as the number of tetrahedra gets large. So,
one might think that triangulated Heegaard splittings offer a nice
alternative to face pairings for the generation of triangulated
$3$--manifolds. However, the probability of getting $2$--symmetric
triangulations on the boundary of layered-triangulations of a
genus-$g$-handlebody goes to zero as the genus gets large. On the
other hand, each equivalence class of minimal triangulations of the
genus $g$ surface appears at infinitely many vertices of the
$L_g$-graph. An interesting question is the distribution of
$2$--symmetric triangulations in the $L_g$-graph. Each gives rise to
a closed, orientable $3$--manifold and by Theorem \ref{triang-HS},
every  closed, orientable $3$--manifold must appear.

For example, in the case of genus $2$, there are  $196$ distinct
minimal triangulations of the genus-$2$-handlebody. We do not know
if all of these are layered-triangulations. However, of those that
are layered-triangulations (and we suspect all are), we can
determine which of these have $2$--symmetric triangulations on their
boundaries; and for those that have a $2$-symmetric triangulation on
the boundary, we can determine the $3$--manifold  presented in this
way. Ben Burton is adding code to {\footnotesize\textsc{REGINA}} to
assist with the this exercise. We have done some examples beyond the
one given in Figure \ref{f-genus2-layered}. The example in Figure
\ref{f-genus2-layered} does not have a $2$--symmetric triangulation
on its boundary. One can easily check that the triangulation induced
on the boundary in this example is $T_6$, using the notation in
\cite{mosher-guide}. Having discovered this triangulation is not
$2$--symmetric, we constructed three additional examples by hand,
two of these had $2$--symmetric triangulations induced on their
boundaries. We give these examples in Figure \ref{f-triang-HS}. The
first is a triangulated Heegaard splitting of $S^2\times S^1$, of
genus $2$ with four tetrahedra; the second is a triangulated
Heegaard splitting of $\rppp\#\rppp$, also of genus $2$ with four
tetrahedra. The first splitting is stabilized. However, the second
is not stabilized but must be reducible as the manifold is
reducible; it is a minimal triangulation of $\rppp\#\rppp$. There
are $15$ manifolds having minimal triangulations with $4$ tetrahedra
and few with more than one such minimal triangulation. Only $5$ of
the $15$ manifolds requiring $4$ tetrahedra are manifolds having
genus $2$. Many of the $196$ minimal triangulations of the
genus-$2$-handlebody probably do not have $2$--symmetric
triangulations on their boundaries and many of those that do have
$2$--symmetric triangulations on their boundaries probably do not
give minimal triangulations. {\footnotesize\textsc{REGINA}} will
give us this answer.

Finally, suppose $M$ is a closed $3$--manifold with a triangulated
Heegaard splitting. Then the number of tetrahedra gives a measure of
complexity of the Heegaard splitting. There are many questions
related to the number of tetrahedra in such a triangulated Heegaard
splitting relative to the genus of the splitting and distance of the
Heegaard splitting as defined in \cite{hempel-distance}. Also, if we
consider all possible triangulated Heegaard splittings of $M$, then
the number of tetrahedra in a minimal triangulated Heegaard
splitting of $M$ is an invariant of $M$. It may be the minimal
number of tetrahedra needed to triangulate $M$. Indeed, we have a
question similar to the conjectures given above.

\begin{ques} With possibly a few exceptions, is a minimal
triangulation of a $3$--manifold a layered-triangulation?\end{ques}

Hopefully with the aid of {\footnotesize\textsc{REGINA}}, we can
explore these questions and gain a better understanding toward the
use of layered-triangulations and triangulated Heegaard splittings
in the study and understanding of $3$--manifolds of genus larger
than one.

\vspace{.15 in}\noindent{\it Layered-triangulations and foliations.}
There is an interesting connection between layered triangulations
and foliations, following Lackenby's view of how to construct taut
ideal triangulations \cite{lack-taut}. We sketch this procedure in
the closed case and will discuss it further in a later paper.  Given
a tetrahedron $\Delta$, choose a pair of opposite edges $e,e'$. Now
we can build a ``compressed" foliation on $\Delta$ with leaves being
a family of quadrilateral disks which all have common boundary being
the 4 edges of $\Delta$ other than $e,e'$. Note that if these four
edges are removed from $\Delta$, we see a genuine foliation. In
Lackenby's construction, a taut foliation can be constructed by
gluing together such foliations on all the ideal simplices and then
slightly uncompressing the leaves along the edges. There is a simple
compatibility condition required to make this work. Note that in our
construction of layered triangulations of handlebodies and closed
manifolds like lens spaces, similar compressed foliations can be
built on all the simplices. The opposite edges $e,e'$ chosen always
to include the one being layered on. It is easy to see that this
foliation can be uncompressed to be the standard product foliation,
when the layered manifold is homeomorphic to a surface $\times I$.
When a closed manifold is obtained from a triangulated Heegaard
splitting by folding a $2$--symmetric triangulation of the boundary
of a handlebody along an edge, then a spine of the surface which is
the $g$-spine can be chosen as the singular limit of the leaves of
the foliation (for example, the center line of a M\"obius band in
case of a solid torus).

An interesting point about this is that normal and almost normal
surfaces behave nicely relative to this singular foliation, when we
use layered triangulations. So a basic idea is that normal surfaces
can have saddle singularities and death singularities relative to
this foliation, but birth singularities can only occur at the limit
leaf. (A birth is the appearance of an inessential simple closed
curve and a death is the disappearance of one. Also, here we are
following normal surface intersections with the leaves in from the
boundary). This is what gives the strong control of normal and
almost normal surface theory in layered solid tori and partial
control in layered-triangulations of more general handlebodies.


\begin{thebibliography}{10}

\bibitem{bing1}
R.~H. Bing.
\newblock {S}ome aspects of the {T}opology of $3$--{M}anifolds {R}elated to the
  {P}oincar{\'{e}} {C}onjecture.
\newblock In {\em Lectures on Modern Mathematics, Vol. II}, pages 92--128. John
  Wiley \& Sons, New York, 1964.

\bibitem{bon-otal}
F.~Bonahon and J.~P. Otal.
\newblock Scindements de {H}eegaard des espaces lenticulaires.
\newblock {\em Ann. Scient. \'Ec. Norm Sup.}, 16:451--466, 1983.

\bibitem{bredon-wood}
Glen~E. Bredon and John~W. Wood.
\newblock Non-orientable surfaces in orientable $3$--manifolds.
\newblock {\em Invent. Math.}, 7:83--110, 1969.

\bibitem{burton-regina}
B.A. Burton.
\newblock \textsc{REGINA}.
\newblock {http://regina.sourceforge.net/}, 2001.

\bibitem{burton-thesis}
B.A. Burton.
\newblock Minimal triangulations and normal surfaces.
\newblock {\em Thesis, University of Melbourne}, pages 1--233, 2003.

\bibitem{casler}
B.G. Casler.
\newblock {An embedding theorem for connected $3$--manifolds with boundary }.
\newblock {\em Proc. Amer. Math. Soc.}, 16:559--556, 1965.

\bibitem{cass-gord}
A.~Casson and C.~McA. Gordon.
\newblock Reducing {H}eegaard splittings.
\newblock {\em Topology Appl.}, 27:275--283, 1987.

\bibitem{fom}
E.~A. Fominykh.
\newblock A complete description of normal surfaces for infinite series of
  $3$--manifolds.
\newblock {\em Siberian Math. J.}, 43(6):1112--1123, 2002.

\bibitem{gord-luecke-planar}
C.~Gordon and J.~Luecke.
\newblock Reducible manifolds and {D}ehn surgery.
\newblock {\em Topology}, 35:385--410, 1996.

\bibitem{harer}
J.~Harer.
\newblock The virtual cohomological dimension of the mapping class group of an
  oriented surface.
\newblock {\em Invent. Math.}, 84:157--176, 1986.

\bibitem{hempel-distance}
John Hempel.
\newblock $3$--manifolds viewed from the curve complex.
\newblock {\em Topology}, 40(3):631--657, 2001.

\bibitem{jac-rubfiniteHS}
W.~Jaco and J.~H. Rubinstein.
\newblock Finiteness of {H}eegaard splittings for $3$--manifolds.
\newblock (in preparation).

\bibitem{jac-rub1}
W.~Jaco and J.~H. Rubinstein.
\newblock $1$--efficient triangulations of 3-manifolds.
\newblock 2002.
\newblock (in preparation).

\bibitem{jac-rub0}
W.~Jaco and J.~H. Rubinstein.
\newblock $0$--efficient triangulations of 3-manifolds.
\newblock {\em J. Diff. Geom.}, 65:61--168, 2003.

\bibitem{jac-sedg-dehn}
W.~Jaco and E.~Sedgwick.
\newblock Decision problems in the space of {D}ehn fillings.
\newblock {\em Topology}, 42:845--906, 2003.

\bibitem{jac-rub-blowup}
William Jaco and J.~H. Rubinstein.
\newblock Blow-ups of ideal triangulations.
\newblock (in preparation).

\bibitem{jac-rub-sSFS}
William Jaco and J.~H. Rubinstein.
\newblock Canonical triangulations and {H}eegaard splittings of small seifert
  fiber spaces.
\newblock (in preparation).

\bibitem{kne}
H.~Kneser.
\newblock {\it Geschlossene {F}l\"achen in dreidimensionalen
  {M}annigfaltigkeiten}.
\newblock {\em Jahresbericht der Deut. Math. Verein.}, 38:248--260, 1929.

\bibitem{lack-taut}
Marc Lackenby.
\newblock {Taut ideal triangulations of 3-manifolds}.
\newblock {\em Geom. Topol.}, 4:369--395, 2000.

\bibitem{italians-3}
B.~Martelli and B.~Petronio.
\newblock Three-manifolds having complexity at most $9$.
\newblock {\em Experimental Math.}, 10:207--236, 2001.

\bibitem{matveev1}
S.V. Matveev.
\newblock {Complexity Theory of three-dimensional manifolds}.
\newblock {\em Acta Appl. Math.}, 19(2):101--130, 1990.

\bibitem{matveev2}
S.V. Matveev.
\newblock Tables of $3$--manifolds up to complexity $6$.
\newblock {\em Max-Planck Institute Preprint}, MPI 1998(67):1--50, 1998.

\bibitem{matveev3}
S.V. Matveev.
\newblock Tables of spines and $3$--manifolds up to complexity $7$.
\newblock {\em Max-Planck Institute Preprint}, MPI 2002(71):1--65, 2002.

\bibitem{matveev-fom}
S.V. Matveev and E.~A. Fominykh.
\newblock Normal surfaces in $3$--manifolds.
\newblock {\em Doklady Math.}, 65(3):727--730, 2002.

\bibitem{mor-rubHS-curved}
Y.~Moriah and J.H. Rubinstein.
\newblock Heegaard structures of negatively curved manifolds.
\newblock {\em Comm.in Geom. and Ann.}, 5(3):375--412, 1997.

\bibitem{mosher-guide}
L.~Mosher.
\newblock A user's guide to the mapping class group: once punctured surfaces.
\newblock In {\em MSRI 1995-001. Geometric and computational perspectives on
  infinite groups}, volume~25 of {\em DIMACS Ser. Discrete Math. Theoret.
  Comput. Sci.}, pages 1--73. Amer. Math. Soc., Providence, RI, 1996.

\bibitem{pachner1}
U.~Pachner.
\newblock P.{L}. homeomorphic manifolds are equivalent by elementary shellings.
\newblock {\em European J. Combin.}, 12(2):129--145, 1991.

\bibitem{perelman}
G.~Perelman.
\newblock The entropy formula for the ricci flow and its geometric
  applications.

\bibitem{rieckHS}
Y.~Rieck.
\newblock Heegaard structures of manifolds in the {D}ehn filling space.
\newblock {\em Topology}, 39:619--641, 2000.

\bibitem{rieck-sedgHSstructure}
Y.~Rieck and E.~Sedgwick.
\newblock Persistence of {H}eegaard structure under {D}ehn filling.

\bibitem{rieck-sedgHSfinite}
Y.~Rieck and E.~Sedgwick.
\newblock Finiteness results for {H}eegaard surfaces in surgered manifolds.
\newblock {\em Topology and its Applications}, 109:41--53, 2001.

\bibitem{rubin-jac}
J.H. Rubinstein and William Jaco.
\newblock Sweep-outs in three-manifolds.

\bibitem{shar-HS}
M.~Scharlemann.
\newblock {H}eegaard splittings.
\newblock In R.J. Daverman and R.B. Sher, editors, {\em Handbook of Geometric
  Topology}, pages 921--955. Elsevier Science, The Netherlands, 2002.

\bibitem{wald-HS3}
F.~Waldhausen.
\newblock Heegaard-zerlegungen der $3$--sph\"are.
\newblock 7:195--203, 1968.

\end{thebibliography}
\end{document}